\documentclass[letterpaper]{article}
\usepackage{authblk}
\usepackage{url}
\usepackage[margin=1in,footskip=0.25in]{geometry}
\usepackage{amssymb}
\usepackage[labelfont=bf]{caption}
\usepackage{graphicx}
\usepackage{amsmath}
\usepackage{upgreek}
\usepackage{bbm}
\usepackage{algorithmic} 
\usepackage{algorithm}
\usepackage{float}
\usepackage{lipsum}
\usepackage{subcaption}
\usepackage{wrapfig}
\usepackage{color}
\usepackage{easylist}
\usepackage{lineno}
\usepackage{stmaryrd}
\usepackage{booktabs}
\usepackage{hyperref}
\usepackage{cleveref}
\usepackage{multirow}
\usepackage{hhline}
\usepackage{soul}
\usepackage{blindtext}
\usepackage{threeparttable}
\usepackage[bbgreekl]{mathbbol}
\usepackage[toc,page]{appendix}
\usepackage[bbgreekl]{mathbbol}

\usepackage{bm}
\usepackage{amssymb}
\usepackage{movie15}
\usepackage{xcolor}
\usepackage{ae,aecompl}
\usepackage{hyperref}

\usepackage[subrefformat=parens,labelformat=parens]{subcaption} 
\usepackage{tabularx} 
\usepackage{multirow}

\usepackage{enumitem}

\usepackage{cite}

\global\long\def\V#1{\boldsymbol{#1}}

\global\long\def\grad{\boldsymbol{\nabla}}

\global\long\def\J{\sM J}
\global\long\def\S{\sM S}

\renewcommand{\vec}[1]{\bm{\mathrm{#1}}}

\def \dt{\Delta t}
\def \b{\vec{b}}
\def \x{\vec{x}}
\def \xp{\x'}
\def \X{\vec{X}}
\def \Xn{\X_m}
\def \Vn{V_m}
\def \Xl{\X_l}

\def \Xp{\X'}

\def \u{\vec{u}}

\def \up{\u'}
\def\Fubar{\vec{F}_{\text{PD}}}

\def \Yubar{\underbar{\vec{Y}}}
\def \tubar{\underbar{t}}
\def \Tubar{\underbar{\vec{T}}}
\def \Mubar{\underbar{\vec{M}}}
\def \horizon{\vec{\mathcal{H}_X}}
\def \dhorizon{\vec{\mathcal{H}}}
\def \PP{\vec{\mathbb{P}}}
\def \FF{\vec{\mathbb{F}}}

\def \BB{\vec{\mathbb{K}}}
\def \Xi{\vec{\xi}}
\def \Eta{\vec{\eta}}

\def \pdrho{\rho_0}
\def \frho{\rho}
\def \fmu{\mu}

\def \v{\vec{v}}
\def \V{\vec{V}}
\def \f{\vec{f}}
\def \F{\vec{F}}

\def \Omegas{\Omega^\text{s}}
\def \vchi{\vec{\chi}}

\def \N{\boldsymbol{\mathcal{N}}}
\def \S{\boldsymbol{\mathcal{S}}}
\def \J{\boldsymbol{\mathcal{J}}}

\def \horizonsize{\epsilon}

\def \sc{s_\text{c}}

\def \dV{ \mathrm{d} {\X}}
\def \dVp{ \mathrm{d} {\Xp}}

\def \omegah{\hat{\omega}}

\title{An immersed peridynamics model of fluid-structure interaction accounting for material damage and failure}
\author[1]{Keon Ho Kim}
\author[2]{Amneet P. S. Bhalla}
\author[3,4,5,6]{Boyce E. Griffith}
\affil[1]{Department of Mathematics, University of North Carolina, Chapel Hill, NC, USA}
\affil[2]{Department of Mechanical Engineering, San Diego State University, San Diego, CA, USA}
\affil[3]{Departments of Mathematics, Applied Physical Sciences, and Biomedical Engineering,
University of North Carolina, Chapel Hill, NC, USA}
\affil[4]{Carolina Center for Interdisciplinary Applied Mathematics, University of North Carolina,
Chapel Hill, NC, USA}
\affil[5]{Computational Medicine Program, University of North Carolina, Chapel Hill, NC, USA}
\affil[6]{McAllister Heart Institute, University of North Carolina, Chapel Hill, NC, USA}
\affil[ ]{\texttt{keonho@email.unc.edu}, \texttt{asbhalla@sdsu.edu}, and \texttt{boyceg@email.unc.edu}}

\date{August 13, 2023}

\begin{document}

\maketitle

\begin{abstract}

This paper develops and benchmarks an immersed peridynamics method to simulate the deformation, damage, and failure of hyperelastic materials within a fluid-structure interaction framework. 
The immersed peridynamics method describes an incompressible structure immersed in a viscous incompressible fluid. 
It expresses the momentum equation and incompressibility constraint in Eulerian form, and it describes the structural motion and resultant forces in Lagrangian form.
Coupling between Eulerian and Lagrangian variables is achieved by integral transforms with Dirac delta function kernels, as in standard immersed boundary methods. 
The major difference between our approach and conventional immersed boundary methods is that we use peridynamics, instead of classical continuum mechanics, to determine the structural forces. 
We focus on non-ordinary state-based peridynamic material descriptions that allow us to use a constitutive correspondence framework that can leverage well-characterized nonlinear constitutive models of soft materials. 
The convergence and accuracy of our approach are compared to both conventional and immersed finite element methods using widely used benchmark problems of nonlinear incompressible elasticity.
We demonstrate that the immersed peridynamics method yields comparable accuracy with similar numbers of structural degrees of freedom for several choices of the size of the peridynamic horizon. 
We also demonstrate that the method can generate grid-converged simulations of fluid-driven material damage growth, crack formation and propagation, and rupture under large deformations.

\end{abstract}


\section{Introduction}
Under extreme loading conditions, structural deformations can cause material damage and, ultimately, failure.
Conventional continuum mechanics models are formulated in terms of \textit{local} definitions of strain, which are evaluated using the derivatives of the deformation field, but such derivatives become ill-defined along propagating crack surfaces, where the deformation field becomes \textit{discontinuous}.
These discontinuities make it challenging to simulate structural mechanics in the face of material failure using conventional continuum approaches, particularly if the failure event involves an evolving crack front that is not known in advance but instead must be determined by the model.
Nonetheless, many continuum approaches to describing material failure have been developed, including linear elastic fracture mechanics  \cite{griffith1921vi} using the finite element (FE) method, and substantial work has created numerical methods for simulating fracture, such as the cohesive zone model \cite{hillerborg1976analysis} and the eXtended finite element method  \cite{moes1999finite}, although these formulations can present challenges themselves in terms of implementation and efficient computation \cite{Sukumar2015}. 



As an alternative to local approaches to fracture mechanics, Silling introduced a \textit{non-local} formulation of solid mechanics called \textit{peridynamics} \cite{silling2000reformulation}, which avoids the use of derivatives in determining strains and instead considers interactions among all material points within a finite \textit{horizon}.
Such non-local models can be better suited for simulating material failure than conventional continuum mechanics formulations because they can tolerate discontinuous deformation fields.
Silling's first theory of peridynamics is now known as \textit{bond-based peridynamics} (BB-PD).
BB-PD asserts that the forces acting between two interacting particles in a material body are completely determined by the relative positions of the particles in the reference and current configurations (i.e., the reference and deformed \textit{bonds}).
This formulation further requires that the resulting \textit{bond forces} are parallel to the deformed bonds in the current configuration.
The assumptions of BB-PD impose strong restrictions on the types of materials that can be modeled.
Specifically, the BB-PD theory can only describe isotropic materials with a Poisson's ratio of $\nu = \frac{1}{3}$ for plane stress and $\nu = \frac{1}{4}$ for plane strain \cite{gerstle2005peridynamic}.
Trageser and Seleson \cite{Trageser2020} demonstrated that these restrictions are the consequence of Cauchy’s relations for isotropic materials in the BB-PD model.

Silling et al.~\cite{silling2007peridynamic} subsequently introduced \textit{state-based peridynamics}, which overcomes many of the limitations of bond-based formulations.
State-based peridynamics adopts the concept of \textit{peridynamic states}.
Examples of peridynamic states include the  force vector state and the deformation vector state, which are detailed in Sec.~\ref{s:PD_states}.
Whereas a constitutive model in a standard continuum mechanics formulation relates the (local) stress and the (local) strain, a constitutive model in the state-based peridynamic theory is a relationship between the force vector state and the deformation vector state. 
In ordinary state-based peridynamics (OSB-PD), the force vectors between any two particles may differ in magnitude but are always parallel to the deformed bond. 
In contrast, the force vectors in non-ordinary state-based peridynamics (NOSB-PD) can differ both in magnitude and direction, which is far more general than BB-PD or OSB-PD and is helpful for developing peridynamic models for real materials. 
In addition, the force vectors can be determined using non-local analogues of the stress and strain tensors in classical continuum mechanics, resulting in so-called constitutive correspondence models, which allow these methods to generate non-local generalizations of existing continuum material models.

%
%

This paper integrates NOSB-PD \cite{silling2007peridynamic, warren2009non} with an immersed fluid-structure interaction (FSI) framework \cite{griffith2020immersed} to simulate the deformations of incompressible hyperelastic materials under FSI. 
We focus on the constitutive correspondence model in NOSB-PD because it enables simulations using well-characterized constitutive models of soft materials within the PD framework.
There has been relatively limited previous work to couple PD models with FSI frameworks to simulate failure in both stiff \cite{BEHZADINASAB2021100045,Shende_2022} and flexible \cite{DallaBarba2020,GAO2020126,ZHANG2021110267,DALLABARBA2022115210} materials under fluid-driven loading conditions.
These studies use BB-PD \cite{DallaBarba2020,ZHANG2021110267} or OSB-PD \cite{BEHZADINASAB2021100045,Shende_2022,GAO2020126,DALLABARBA2022115210}, however, and we are unaware of prior work using NOSB-PD within an FSI framework.

To treat FSI, we use the framework of the immersed boundary (IB) method \cite{peskin_2002}. 
The IB method for FSI was introduced by Peskin in the 1970's  to simulate the dynamics of heart valves \cite{peskin1972flow,peskin1977numerical}, and it has been subsequently used in a broad range of applications, including cardiovascular dynamics \cite{PESKIN1989372, MCQUEEN1989289, McQueen_2000,  McQueen_book, griffith_2009, griffith2012immersed, Lee_2020, LEE202160, Choi_2014, Bailoor_2021, KOLAHDOUZ2020108854}, esophageal transport \cite{Kou_2015, Kou_2017}, aquatic locomotion  \cite{Bhalla_2013, BHALLA2013446, HERSCHLAG201184, Koumoutsako_2008, BORAZJANI20087587, Borazjani_2008}, and insect flight \cite{fluids3030045, Wang_2005}. 
The IB method uses Lagrangian variables for the deformations, stresses, and resultant forces of the immersed structure and Eulerian variables for the momentum, viscosity, and incompressibility of the coupled fluid-structure system. 
Coupling between Lagrangian and Eulerian variables is mediated by integral transforms with Dirac delta function kernels in the continuous formulation. 
In discretized IB formulations, the singular delta function is replaced by a regularized delta function \cite{peskin_2002}. 
This coupling approach enables nonconforming discretizations of the fluid and immersed structures \cite{griffith2020immersed,griffith2017hybrid}. 
Conventional IB methods use the framework of nonlinear continuum mechanics to compute the elastic body forces of the immersed structure \cite{griffith2017hybrid, BOFFI20082210, ZHANG20042051, Wells:2021uh}. 
In this work, the constitutive correspondence model of NOSB-PD is used to compute structural forces instead of continuum mechanics.

The key contribution of this work is that it develops an immersed peridynamics (IPD) method using NOSB-PD to enable simulations including nonlinear hyperelastic material responses, even with discontinuities in the displacement field (i.e., material damage and failure), which is challenging with continuum-based IB methods. 
A second contribution of this study is that it investigates the performance of NOSB-PD models of hyperelastic materials.
To date, there has been a relatively limited amount of research on the use of PD to simulate such rubber-like materials \cite{Diehl2022}.
The IPD method developed here considers a hyperelastic structure immersed in the surrounding fluid. 
We consider several standard benchmark cases that were first introduced to test solvers for conventional formulations of large-deformation elasticity through a variety of solid mechanics benchmarks in both two and three spatial dimensions, including the compressed block \cite{reese1999new}, Cook's membrane \cite{cook1974improved}, and the torsion test \cite{bonet1997nonlinear}, which make it straightforward to compare the quasi-static results under FSI to classical quasi-static solid mechanics benchmarks.
All benchmarks show that the accuracy of the structural deformations obtained by the proposed method is comparable to that yielded by a stabilized finite element method \cite{reese1999new} and a stabilized immersed finite element-finite difference method \cite{vadala2020stabilization} for incompressible nonlinear elasticity.
Consequently, this work advances our understanding of the fidelity of hyperelastic structural responses within the PD framework.
Material responses under fluid traction boundary conditions are also investigated using both static and dynamic versions of an elastic band test \cite{vadala2020stabilization}. 
A modified elastic band with and without an initial crack is used to simulate purely fluid-driven material failure.
Some test cases involve interactions between flexible and fixed structures. 

\section{Non-ordinary state-based peridynamics}
\label{s:NOSB-PD}
This section presents a NOSB-PD formulation \cite{silling2007peridynamic, warren2009non} that can be used with material models that are characterized by strain energy functionals, such as those often used with continuum mechanics descriptions of hyperelastic materials. 
We describe the continuous formulation of the peridynamic model that we use in our IPD method for FSI

\subsection{Peridynamic states}

\label{s:PD_states}
\begin{figure}[]
	\centering
	\includegraphics[width = 0.7\textwidth]{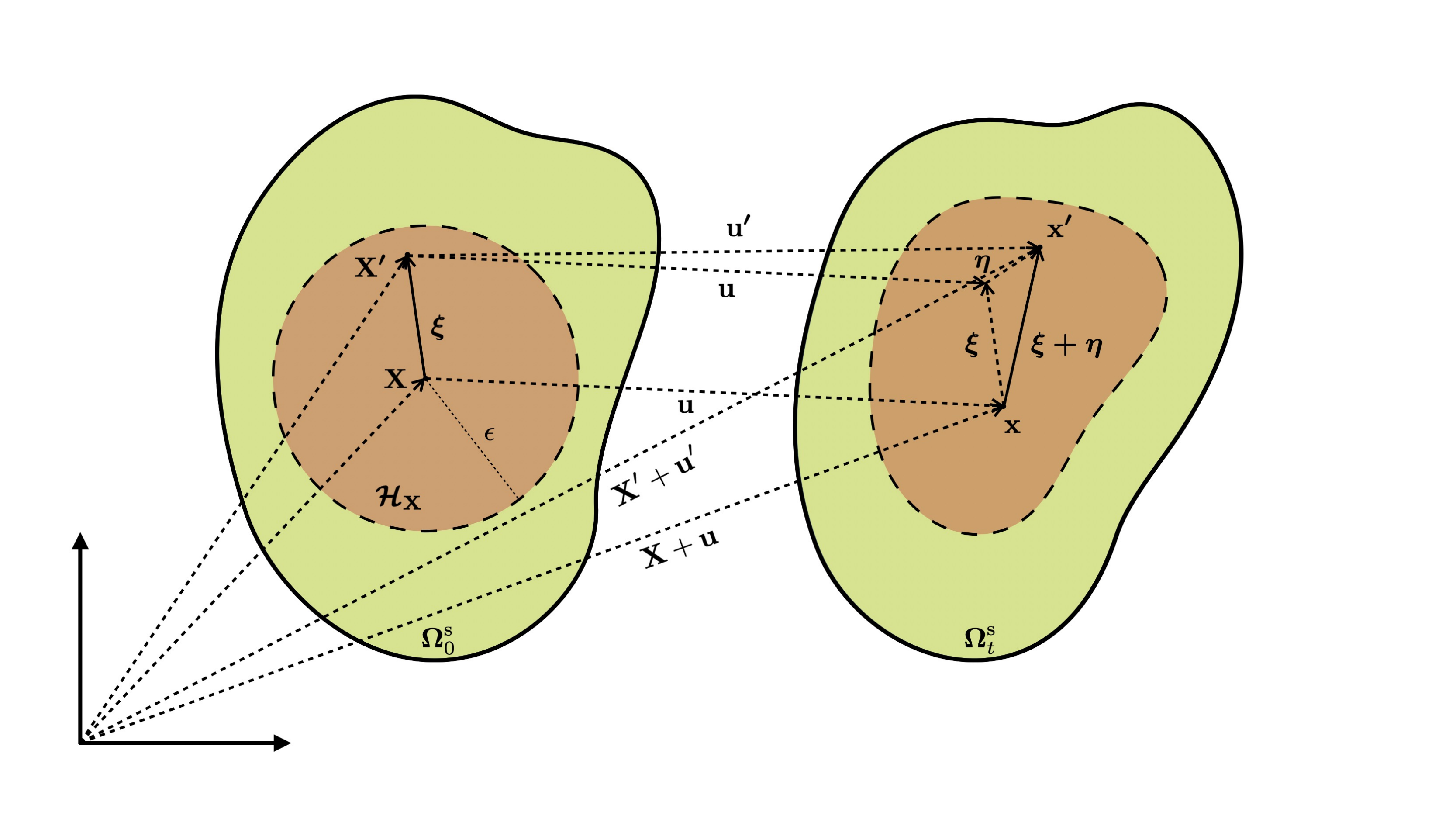}
	\caption{Peridynamic material points, bonds, and horizons in the reference and current configurations.}
	\label{f:PD_coordinates}
\end{figure}

The dynamics of an elastic body occupying the region $\Omegas_t$ at time $t$ are described using material (reference) coordinates $\X \in \Omegas_0$, and the deformation mapping $\vchi : (\Omegas_0,t) \mapsto \Omegas_t$ relates the reference and deformed coordinate systems, so that $\x = \vchi(\X,t)$ is the physical position of the material point $\X$ at time $t$.
Each material point $\X$ in peridynamics interacts with all material points within a finite region around $\X$ called the \textit{horizon} and denoted as $\horizon$.
Here, we choose $\horizon$ to be the sphere of radius $\horizonsize >0$ centered at $\X$, but other choices of the horizon shape have also been used in practice (e.g., a cube of side length $2\horizonsize$)\footnote{In the peridynamics literature, the horizon size is commonly denoted by $\delta$, but here we reserve that symbol for the Dirac delta function.}.
A material point $\X$ interacts with each point $\Xp$ in its horizon $\horizon$ through a \textit{bond}, which is defined in the reference frame by $\Xi = \Xp - \X$, and which deforms to $\Xi + \Eta = \xp - \x$ in the deformed frame, with $\xp = \vchi(\Xp,t)$ and $\x = \vchi(\X,t)$.
See Fig.~\ref{f:PD_coordinates}. 
$\u$ and $\up$ represent the displacements of the material points labeled by $\X$ and $\Xp$ from the reference configuration to the deformed configuration, respectively, so that $\x = \X + \u$ and $\xp = \Xp + \up$. 
The deformed bond is $\xp - \x = (\Xp + \up) - (\X + \u) = ( \Xp -\X) + ( \up - \u ) =  \Xi + \Eta$. The definition of \textit{deformation vector state} $\Yubar$ is motivated by this consideration:
\begin{align}\label{deformation_vec_state}
    \Yubar = \Yubar[\X,t]\langle \Xi \rangle = \vchi(\X + \Xi,t) - \vchi(\X,t) = \xp - \x.
\end{align}
Thus, $ \Yubar[\X,t]\langle \Xi \rangle$ is the deformation of the bond $\Xi$ associated with material point $\X$ at time $t$.

The \textit{force vector state} associated with a bond $\Xi$ acting on material point $\X$ at time $t$ is denoted by
\begin{align}
\Tubar = \Tubar[\X,t]\langle\Xi\rangle.
\end{align} 
Constitutive models in state-based peridynamics provide the values for the force vector state field based on the deformation state and possibly other variables:
\begin{align}
   \Tubar = \hat{\Tubar}(\Yubar,\Lambda),
\end{align}
in which $\Lambda$ denotes all variables that determine the force state other than the deformation state.  There are several restrictions on the force vectors based on different material models.
\begin{figure}[]
	\centering
	\includegraphics[width = 0.9\textwidth]{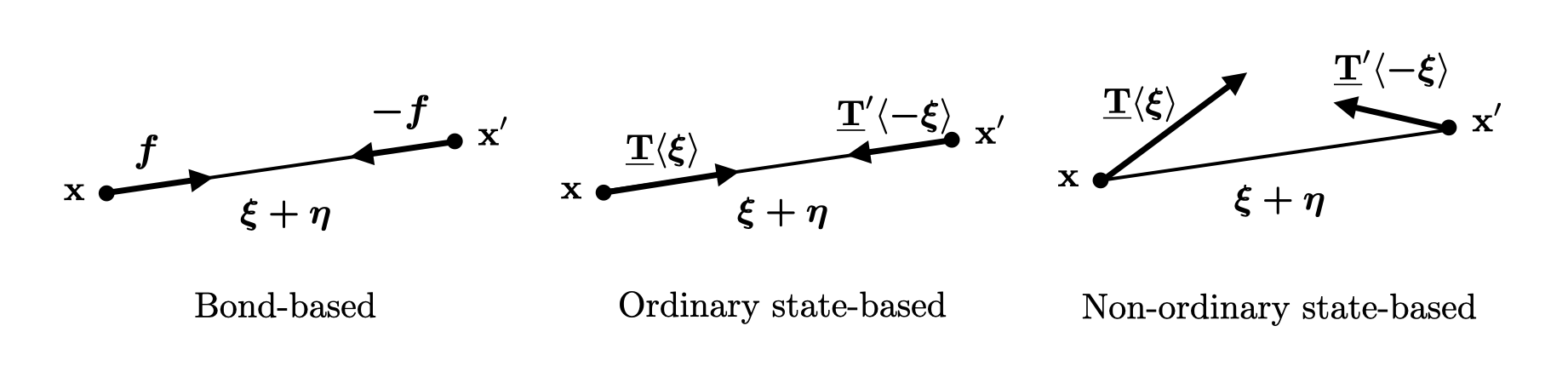}
	\caption{Force vector states for bond-based, ordinary state-based, and non-ordinary state-based peridynamics.}
	\label{f:PD_thoeries}
\end{figure}
For example, the force vector state of an \emph{ordinary} material is defined by
\begin{align}\label{deformed_dir_vec}
\Tubar = \tubar \ \Mubar,
\end{align}
in which $\Mubar$ is the deformed direction vector state that is the unit directional vector of the deformed vector state $\Yubar$ and $\tubar$ is the scalar force state. 
In ordinary peridynamic material models, the force vectors are parallel to its deformed bond. 
If the force vectors are not parallel to its corresponding bond in the deformed configuration, the material description is called non-ordinary. 
Here, we primarily focus on NOSB-PD for nonlinear material responses because it allows us to leverage existing hyperelastic constitutive laws for soft materials. 
Fig.~\ref{f:PD_thoeries} illustrates possible force vector states for bond-based, ordinary state-based, and non-ordinary state-based peridynamics. 
Silling et al.~\cite{silling2007peridynamic} provide additional details about peridynamic states.

\subsection{Kinematics and balance laws}
In a state-based peridynamic model, the basic equation of motion is
\begin{align}
    \pdrho \frac{\partial^2 \u}{\partial t^2}(\X,t) = \int_{\horizon} \left( \Tubar[\X,t]\langle\Xi\rangle - \Tubar[\Xp,t]\langle-\Xi\rangle \right) \dVp + \b(\X,t),
\end{align}
in which $\pdrho$ is the mass density of the material and $\b(\X,t)$ is an external force density. 
The pairwise bond force function $\Fubar(\Xp,\X,t)$ is
\begin{align}
\Fubar(\Xp,\X,t) &  = \Tubar[\X,t]\langle\Xi\rangle - \Tubar[\Xp,t]\langle-\Xi\rangle,
 \label{eqn_f_pd}            
\end{align}
which accounts for contributions of the material model at both $\X$ and $\Xp$ to balance linear momentum. 
Note that force balance is automatically satisfied for arbitrary force vector states: 
\begin{equation}
\Fubar(\X,\Xp,t) = \Tubar[\X,t]\langle\Xi\rangle - \Tubar[\Xp,t]\langle-\Xi\rangle =- \left(\Tubar[\Xp,t]\langle-\Xi\rangle - \Tubar[\X,t]\langle\Xi\rangle \right)  = - \Fubar(\Xp,\X,t),
\end{equation}
which is consistent with Newton's third law of motion.
Conservation of linear momentum is guaranteed in all three peridynamic models for any definition of the force vector states:
\begin{equation}
\begin{aligned}
\label{balance_of_linear_momentum}
&\int_{\Omegas_0} \int_{\horizon} \left(  \Tubar[\X,t]\langle\Xi\rangle - \Tubar[\Xp,t]\langle-\Xi\rangle \right)  \dVp \dV,\\
&= \int_{\Omegas_0} \int_{\Omegas_0} \left(  \Tubar[\X,t]\langle\Xi\rangle - \Tubar[\Xp,t]\langle-\Xi\rangle \right)  \dVp \dV,\\
&= \int_{\Omegas_0} \int_{\Omegas_0} \Tubar[\X,t]\langle\Xi\rangle   \dVp \dV -  \int_{\Omegas_0} \int_{\Omegas_0} \Tubar[\Xp,t]\langle-\Xi\rangle  \dVp \dV,  \\
&= \int_{\Omegas_0} \int_{\Omegas_0} \Tubar[\X,t]\langle\Xi\rangle  \dVp \dV-  \int_{\Omegas_0} \int_{\Omegas_0} \Tubar[\X,t]\langle\Xi\rangle  \dV \dVp,\\
&= \int_{\Omegas_0} \int_{\Omegas_0} \Tubar[\X,t]\langle\Xi\rangle   \dVp \dV -  \int_{\Omegas_0} \int_{\Omegas_0} \Tubar[\X,t]\langle\Xi\rangle  \dVp \dV,\\
&=  \vec{0}.
\end{aligned}
\end{equation}
Note that  $\Tubar[\X,t]\langle\Xi\rangle = \Tubar[\Xp,t]\langle-\Xi\rangle= \vec{0}$ for all $ \Xp \notin \horizon$, so that the integration domain can be taken as the entire solid domain $\Omegas_0$. 

Balance of angular momentum holds if
\begin{align}\label{balance_of_angular_momentum}
\int_{\horizon} \left( \Yubar[\X,t]\langle \Xi \rangle \times \Tubar[\X,t]\langle\Xi\rangle \right) \dVp = \vec{0}.
\end{align}
Eq.~\eqref{balance_of_angular_momentum} is automatically satisfied if the force vector $\Tubar$ is aligned with its corresponding deformed bond $\Yubar$; however, the force vector in NOSB-PD is not parallel to its relative position vector in general, so there is a need to define the force vector state in a way that satisfies the conservation of angular momentum.

\subsection{Constitutive correspondence}
In classical continuum mechanics, the \textit{local deformation gradient tensor} is 
\begin{align}\label{deformation_gradient}
\FF_{\text{local}} = \frac{\partial\vchi}{\partial \X},
\end{align}
which determines the local spatial deformation of a point. 
The local stress and strain tensors are computed using the local deformation gradient tensor with various choices of constitutive models. 
For a hyperelastic material model, the first Piola-Kirchhoff stress tensor is determined from a strain energy functional $\Psi$ by
\begin{align}\label{piola_kirchhoff}
\PP = \frac{\partial \Psi \left( \FF_{\text{local}}\right)}{\partial \FF_{\text{local}}}.
\end{align} 


Developing a peridynamic model that corresponds to classical nonlinear hyperelasticity in cases that do not involve material failure requires constructing non-local analogues of the deformation gradient tensor and stress. 
The \textit{non-local deformation gradient tensor} is 
\begin{equation}\label{nonlocal_deformation_gradient}
    \FF_{\text{non-local}}  = \left[ \int_{\horizon} \omega (| \Xi |) \; \left( \Yubar\langle \Xi \rangle \otimes \Xi \right) \; \dVp  \right] \BB^{-1},
\end{equation}	
in which $\omega (| \Xi |)$ is a non-negative scalar valued function called the \textit{influence function}, which controls the influence of peridynamic points away from the current point, and $\BB$ is the \textit{shape tensor}, 
\begin{equation}\label{shape_tensor}
	\BB =  \int_{\horizon} \omega (| \Xi |) \; \left( \Xi \otimes \Xi \right) \; \dVp.
\end{equation} 
If the deformation mapping $\vchi$ is continuously differentiable and $\horizon \subset \Omegas_0$, then, by a Taylor expansion, Eq.~\eqref{nonlocal_deformation_gradient} can be written as
\begin{equation}
\begin{aligned}
\FF_{\text{non-local}} =& \left[ \int_{\horizon} \omega (| \Xi |) \; \left( \left(\vchi(\X + \Xi,t) - \vchi(\X,t) \right) \otimes \Xi \right) \; \dVp  \right] \BB^{-1},\\
=& \left[ \int_{\horizon} \omega (| \Xi |) \; \left( \left( \frac{\partial \vchi}{\partial \X} \Xi + O (|\Xi|^2) \right) \otimes \Xi \right) \; \dVp  \right] \BB^{-1},\\
=& \left[ \int_{\horizon} \omega (| \Xi |) \; \left( \Xi  \otimes \Xi \right) \; \dVp  \right] \BB^{-1} \frac{\partial \vchi}{\partial \X}   + O (\horizonsize^2),\\
=& \ \FF_{\text{local}}  + O(\horizonsize^2).\label{nonlocal-local}
\end{aligned}
\end{equation}
Consequently, away from the boundary of the structural domain, the non-local deformation gradient tensor is a second-order accurate approximation to the local deformation gradient tensor used in the conventional continuum theory. Near the boundary of the domain, the shape of the peridynamic horizon or the influence function is not symmetric, however, and the non-local deformation gradient tensor is only a first-order approximation to the local deformation gradient tensor \cite{Hillman_2019}.
Notice that the non-local deformation gradient tensor is not defined only on a particular bond. 
Rather, it provides an averaged description of all interacting bonds in the horizon of a material point. 
The integral operators in the peridynamic theory allow for descriptions of discontinuities in materials, such as cracks and fractures. 
To simplify notation for the remainder of the paper, we denote the non-local deformation gradient tensor in NOSB-PD by $\FF$.

In the NOSB-PD constitutive correspondence introduced by Silling et al.~\cite{silling2007peridynamic}, the force vector state is defined by 
\begin{align}\label{general_force_vector}
    \Tubar[\X,t]\langle\Xi\rangle = \omega (| \Xi |) \PP \BB^{-1} \Xi,
\end{align}
in which the first Piola-Kirchhoff stress tensor $\PP$ is computed using the non-local deformation gradient tensor instead of the local deformation gradient tensor in Eq.~\eqref{piola_kirchhoff}.
Note that Eq.~\eqref{general_force_vector} guarantees the balance of angular momentum \cite{silling2007peridynamic,  warren2009non, madenci2014peridynamic}. 
In brief, once the non-local deformation gradient $\FF$ is determined, the first Piola-Kirchhoff stress tensor $\PP$ can be obtained from a classical constitutive material model, which thereby allows the peridynamic force vector to be determined.
This relation between the conventional stress tensor in continuum mechanics and the peridynamic force vector state is the so-called constitutive correspondence \cite{silling2007peridynamic}.

\subsection{Failure and damage}
\label{s:PD_failure_damage}
In peridynamics, the formation and propagation of a crack occurs as bonds break between the Lagrangian material points.
In this study, a critical stretch criterion is used to determine if a bond breaks. 
Other failure criteria have also been suggested, including strain-based criteria \cite{silling2005meshfree} and energy-based criteria \cite{foster2011energy}. 
Bond breakage is modeled as an irreversible process: once a bond breaks, it cannot be reformed. The bond stretch is
\begin{align}\label{bond_stretch}
s = \frac{|\Xi + \Eta|}{|\Xi|}.
\end{align}
If the bond stretch exceeds its critical value $\sc$, there is no longer interaction between the two material points connected by the bond, i.e, the bond breaks. 
To track the connectivity between two material points that are initially connected by the bond $\Xi$ under deformations, we use an indicator function $I(\Xi,t)$ \cite{behera2020peridynamic}:
\begin{align}\label{fail_paramter}
I \left( \Xi, t  \right) = 
\begin{cases} 
1, \ s \le \sc, \\
0, \ s> \sc.
\end{cases}
\end{align}
We use a modified influence function  $\omegah (|\Xi|,t) =  \omega (|\Xi| ) I \left( \Xi, t  \right)$ that takes the value of 0 if the bond $\Xi$ breaks.
This also implies that the deformation gradient tensor and force vector state, Eqs.~\eqref{nonlocal_deformation_gradient} and \eqref{general_force_vector}, must be modified after a bond breaks in the material horizon.  

The local damage at a material point $\X$ at time $t$ can be computed by the bond connectivity within the peridynamic horizon $\horizon$ as \cite{behera2020peridynamic}
\begin{align}\label{damage}
\varphi(\X,t) = 1 - \frac{\int_{\horizon} I (\Xi,t) \dVp}{\int_{\horizon} \dVp},
\end{align}
which is a volume-weighted ratio of the number of eliminated bonds to the number of initial bonds at a material point within the horizon. 
Note that the local damage is equal to $0$ if its initial bonds are all active and its value is equal to $1$ if all bonds are broken.

\section{Immersed peridynamics method}
This section presents the continuous and discrete IPD formulations for simulating FSI with and without material damage and failure. 

\subsection{Continuous formulation}
\begin{figure}[]
	\centering
	\includegraphics[width = .9\textwidth]{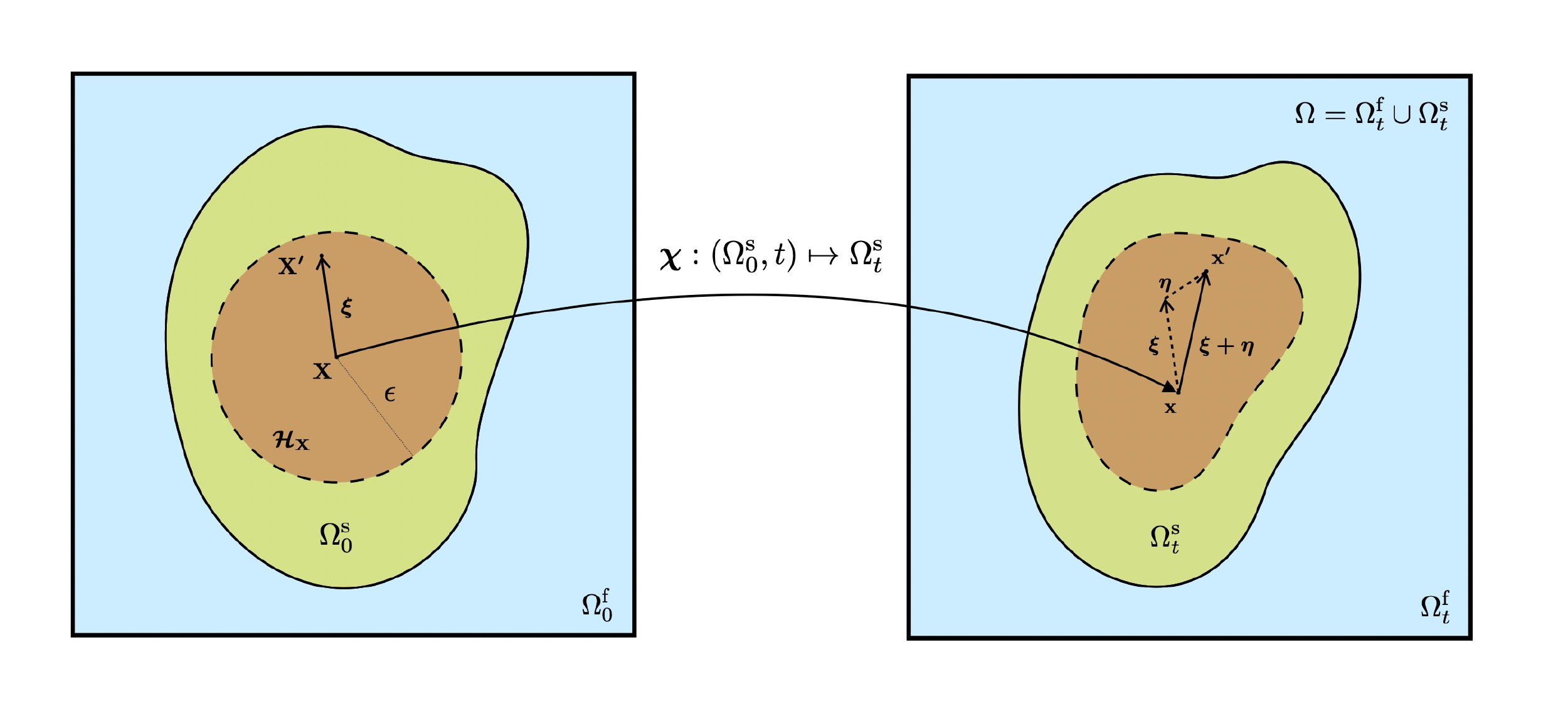}
	\caption{Lagrangian and Eulerian coordinate systems in the immersed peridynamics formulation with peridynamic material points, bonds, and horizons. The computational domain is divided into the solid and fluid subdomains, $\Omega_t^{\text{s}} $ and $ \Omega_t^{\text{f}}$ at time $t$, respectively. The peridynamic point $\X \in \Omega_0^{\text{s}}$ interacts with its neighborhood in the finite region called horizon, which is denoted by $\horizon \subset \Omega_0^{\text{s}}$.}
	\label{f:IPD_mapping}
\end{figure}
The continuous IPD formulation considers an Eulerian computational domain $\Omega$ that consists of a time-dependent fluid subdomain $ \Omega_t^{\text{f}}$ and solid subdomain $\Omega_t^{\text{s}} $ that are indexed by time $t$. 
We use both fixed spatial coordinates $\x \in \Omega$ and reference coordinates $\X \in\Omega_0^{\text{s}}$, with $\Omega_0^{\text{s}}$ indicating the region occupied by the solid structure at time $t = 0$. 
The dynamics of the fluid-structure system are described by
\begin{align}
\frho \frac{\mathrm{D} \v}{\mathrm{D} t} \left( \x,t \right)  &= - \grad p\left( \x,t \right) + \fmu \grad^2 \v\left( \x,t \right) + \f \left( \x,t \right), \label{navier_stokes} \\
\grad \cdot \v \left( \x,t \right) & = 0, \label{incompressibility} \\
\f\left( \x,t \right) &= \int_{\Omega_0^{\text{s}}} \F (\X , t)\, \delta \left(\x - \vchi(\X,t)\right) \mathrm{d}\X, \label{force_spreading}\\
\frac{\partial \vchi}{\partial t} (\X,t) &= \V (\X,t) =  \int_{\Omega}  \v (\X,t)\, \delta \left(\x - \vchi(\X,t)\right) \mathrm{d}\x = \v \left(\vchi(\X,t),t\right), \label{velocity}
\end{align}
in which $\frho$ is the mass density, $\fmu$ is the viscosity, $\v (\x,t)$ is the Eulerian velocity, $\V (\X,t)$ is the velocity of the structure, $p (\x,t)$ is the Eulerian pressure field, $\f (\x,t)$ is the Eulerian structural force density generated by the deformations of the structure, and $\F (\X,t)$ is the Lagrangian structural force density. 
The operators $\grad$, $\grad^2$, and $\grad \cdot$ are with respect to spatial (current) coordinates, and $\frac{\mathrm{D}}{\mathrm{D}t} = \frac{\partial}{\partial t} + \v \cdot \grad$ is the material time derivative in current coordinates.  
Eqs.~\eqref{navier_stokes} and \eqref{incompressibility} are the incompressible Navier-Stokes equations. 
The Eulerian and Lagrangian variables are coupled by integral transforms with Dirac delta function kernels in Eqs.~\eqref{force_spreading} and \eqref{velocity}. 
Notice that Eq.~\eqref{velocity} implies that the no-slip boundary condition holds along the fluid-solid interface. 
Eq.~\eqref{velocity} also implies that the structure is exactly incompressible, because its motion is determined by $\v\left(\x,t\right)$, which satisfies the incompressible constraint, Eq.~\eqref{force_spreading}, throughout the entire computational domain $\Omega$. 
See Fig.~\ref{f:IPD_mapping}.

Some cases in the present study include both fixed and elastic parts. 
We approximately impose the zero-displacement constraint $\frac{\partial \vchi}{\partial t} (\X,t) = \vec{0}$ by an approximate Lagrangian multiplier force $\F_{\text{c}}$ of the form:
\begin{align}\label{Lagrangian_multiplier}
\F_{\text{c}} (\X,t) = \kappa \left( \X - \vchi\left( \X,t \right) \right) - \eta \V(\X,t),
\end{align}
in which $\kappa \geq 0$ is a stiffness penalty parameter and $\eta \geq 0$ is a damping penalty parameter. 
Note that as $\kappa \rightarrow \infty$, $\vchi(\X,t) \rightarrow \X$, and $\frac{\partial \vchi}{\partial t} (\X,t) \rightarrow \vec{0}$, and we recover the exactly stationary model.
Adding a damping term reduces numerical oscillations that can occur with finite values of $\kappa$. 
This approach corresponds to an approximate Lagrange multiplier approach like those developed by Glowinski et al. \cite{GLOWINSKI1999755,PATANKAR20001509,GLOWINSKI2001363}.

\subsection{Discrete formulation}
\subsubsection{Eulerian discretization}
The incompressible Navier-Stokes equations are discretized in space using the second-order finite difference scheme on a staggered Cartesian grid \cite{griffith2009accurate}.
For simplicity, we describe the spatial discretization in two spatial dimensions. 
The extension to three spatial dimensions is straightforward. 
The computational domain $\Omega = [0, L]^2$ is discretized by an $N \times N$ Cartesian grid with a uniform grid spacing $h = L/N$ in the $x$- and $y$-directions. 
Let $(i,j)$ label the Cartesian grid cell for integer values of $i$ and $j$ and $0 \le i,j < N$. 
The discrete Eulerian velocity $\v = (v_1,v_2)$ is defined by vector components that are normal to the edges of the Cartesian grid cells at positions, $\x_{i-\frac{1}{2},j} = \left(ih, \left( j + \frac{1}{2}h \right) \right)$ and $\x_{i,j-\frac{1}{2}} = \left(\left( i + \frac{1}{2}h \right), jh \right)$. 
In three spatial dimensions, the discrete velocity is determined by the components of the velocity at the faces of the grid cells.
The components of the discretized elastic body force density $\f = \left( f_1, f_2 \right)$ are defined at the same locations as the velocity. 
The pressure $p$ is approximated at the centers of the Cartesian grid cells at positions $\x_{i,j} = \left( \left(i + \frac{1}{2}\right)h,\left(j + \frac{1}{2}\right)h \right)$.  
The nonlinear advection term $\v \cdot \grad \v$ is computed using a version of the piecewise parabolic method \cite{colella1984piecewise}.

\subsubsection{Lagrangian discretization}
For simplicity, we discretize the structure using a uniformly distributed point cloud, or lattice.
In such a description, we assign uniform discrete volume elements to each material point.
Let $\dhorizon_{\Xl} \subset \Omega_0^{\text{s}}$ be the set of interacting neighborhoods of radius $\horizonsize$ centered at the Lagrangian marker $\Xl \in \Omega_0^{\text{s}}$.
Then it is natural that the spatial integrals in the continuous NOSB-PD formulation are discretized as volume weighted sums. 
The volumes associated with interacting PD points near the boundary of $\dhorizon_{\Xl}$ are partially located inside of the $\epsilon-$ball. 
However, the contributions of interacting PD nodes are added up to calculate overall quantities, such as PD net body forces, which can cause a larger amount of the discrete PD volumetric force to act on a PD node in the discrete IPD formulation compared to the exact evaluation of the PD force.
The accuracy of numerical solutions is improved by using a volume correction method \cite{hu2010numerical, Seleson2019}:
\begin{align}
\Vn^{(l)} = \begin{cases}
\Vn \ &\text{if} \ |\Xi| \le \horizonsize - \frac{\Delta X}{2}, \\
\frac{1}{\Delta X} \left[ \horizonsize - \left(|\Xi| - \frac{\Delta X}{2}\right) \right] \Vn \ &\text{if} \ |\Xi| \le \horizonsize, \\
0 \ &\text{otherwise},
\end{cases}
\end{align}
in which $\Xi$ is a bond connecting material points $\Xl$ and $\Xn$ in the reference configuration and $\Vn$ is the volume occupied by material point $\Xn$.
In our computations, we use the same corrected volumes in both two and three spatial dimensions.

Now we focus on computing the Lagrangian force density $\F$, Eq.~\eqref{force_spreading}, in the IPD formulation, which uses the NOSB-PD constitutive correspondence model to obtain the internal elastic body force. 
The non-local deformation gradient tensor $\FF$, Eq.~\eqref{nonlocal_deformation_gradient}, and shape tensor $\BB$, Eq.~\eqref{shape_tensor}, at material point $\Xl$ are discretized into finite sums,
\begin{align}
\FF_l &= \sum_{\Xn \in \dhorizon_{\Xl}}  \ \omega \left(| \Xn - \Xl |\right) \Yubar \langle \Xn - \Xl \rangle \otimes \left( \Xn - \Xl \right) \BB_l^{-1} \, \Vn^{(l)} , \label{discrete_nonlocal_deformation_gradient}\\
\BB_l &= \sum_{\Xn \in \dhorizon_{\Xl}} \ \omega \left(| \Xn - \Xl |\right) \left( \Xn - \Xl \right) \otimes \left( \Xn - \Xl\right) \, \Vn^{(l)} , \label{discrete_shape_tensor}
\end{align}
in which $\Xn$ is a neighborhood of $\Xl$ in the peridynamic horizon $\dhorizon_{\Xl}$. 
Likewise, the discretized force vector state of material point $\Xl$ is 
\begin{align}\label{discrete_force_vector}
\Tubar[\Xl,t]\langle\Xn - \Xl\rangle = \omega (| \Xn - \Xl |) \PP_l \BB_l^{-1} \left( \Xn - \Xl \right),
\end{align}
and the pairwise bond force function is 
\begin{align}\label{discrete_pairwise_bond_f}
\Fubar \left( \Xl, \Xn, t \right)  = \omega (| \Xn - \Xl |) \left( \PP_l \BB_l^{-1} +  \PP_m \BB_m^{-1} \right) \left( \Xn - \Xl \right),
\end{align}
in which $\PP_m$ and $\BB_m$ are the discretized first Piola-Kirchhoff stress tensor and shape tensor of particle $\Xn$, respectively. The discretized first Piola-Kirchhoff stress tensor is computed by the classical constitutive relations, such as Saint-Venant or neo-Hookean material models, but using the discretized non-local deformation gradient tensor here instead of the classical deformation gradient tensor. Consequently, the net internal body force density at material point $\Xl$ is 
\begin{align}\label{discrete_pd_force}
\F (\Xl,t) =  \sum_{\Xn \in \dhorizon_{\Xl}} \Fubar(\Xl,\Xn,t) \, \Vn^{(l)}.
\end{align}
This peridynamic net bond force is used as an elasticity model for the immersed structure at material point $\Xl$ at time $t$. 


To control the contribution of each Lagrangian point in the peridynamic horizon, we use the influence function defined by
\begin{align}
\omega(r) = 
\begin{cases}
C \left( \frac{2}{3} - r^2 + \frac{r^3}{2} \right) \ &\text{if} \ r < 1.0, \\
C\frac{(2-r)^3}{6} \ &\text{if} \ r \le 2.0, \\
0  \ &\text{otherwise},
\end{cases}
\end{align}
in which $r = \frac{2|\Xi|}{\horizonsize}$ and $C = \frac{15}{7\pi}$ in the two spatial dimensions or $C = \frac{3}{2\pi}$ in the three spatial dimensions. 
Seleson et al.~\cite{Seleson_2011} provide a detailed discussion about the role of influence functions in the peridynamic theory.

Simulating a failure process during the deformation requires the modification of the discretized non-local deformation tensor and peridynamic force vectors based on the connectivity of internal bonds in an immersed structure at each time. 
Therefore, for models that involve material failure in our IPD simulations, we replace the influence function in the integral equations, Eqs.~\eqref{discrete_nonlocal_deformation_gradient}--\eqref{discrete_pairwise_bond_f}, to the modified influence function $\omegah = \omega I$ as explained in Sec.~\ref{s:PD_failure_damage}.

\subsubsection{Lagrangian-Eulerian coupling}
In the continuous equations, coupling between Eulerian and Lagrangian variables is achieved by integral transforms with Dirac delta function kernels as Eqs.~\eqref{force_spreading}--\eqref{velocity}. 
In the discrete formulation, the singular delta function is replaced by a regularized delta function $\delta_h$, which is formed as a tensor product of one-dimensional kernel functions,
\begin{align}\label{regularized_delta_func}
\delta_h (\x) = \Pi_{i=1}^2 \delta_h(x_i) = \frac{1}{h^2} \phi \left( \frac{x_1}{h}\right) \phi \left( \frac{x_2}{h}\right),
\end{align}
in which $\phi(r)$ is a basic one-dimensional kernel function \cite{peskin_2002}. 
We use the four-point IB kernel function introduced by Peskin \cite{peskin_2002} unless otherwise mentioned. 

The immersed body $\Omega_t^{\text{s}}$ is discretized as a collection of Lagrangian points. Then the volume integral, Eq.~\eqref{force_spreading}, is approximated by 
\begin{align}
\left( f_1 \right)_{i - \frac{1}{2},j} &= \sum\limits_{l}  F_{l,1} \, \delta_h \left( \x_{i - \frac{1}{2},j} - \vchi \left( \Xl,t \right) \right) h^2,\\
\left( f_2 \right )_{i,j - \frac{1}{2}} &= \sum\limits_{l} F_{l,2} \, \delta_h \left( \x_{i,j - \frac{1}{2}} - \vchi \left( \Xl,t \right) \right) h^2,
\end{align}
in which $\F_l = \left( F_{l,1}, F_{l,2} \right)$ is the Lagrangian force density at a Lagrangian marker of index $l$. Note that the Lagrangian force density is computed by NOSB-PD in the proposed method.
We use the notation
\begin{align}
\f = \S\left[ \vchi \left(\cdot,t\right)\right] \F,
\end{align}
in which $\S\left[\vchi \left(\cdot,t\right)\right]$ is the discrete force-spreading operator. The structural body interacts with the surrounding fluid by spreading the force to the Eulerian grid and moves with the local fluid velocity. Similarly, Lagrangian and Eulerian velocities are related by
\begin{align}
V_{l,1}  \left(\X,t\right) &= \sum\limits_{i,j} (v_1)_{i-\frac{1}{2},j} \, \delta_h \left( \x_{i - \frac{1}{2},j} - \vchi \left( \X,t \right) \right) h^2, \label{velocity_interpolation_1}\\
V_{l,2}  \left(\X,t\right) &= \sum\limits_{i,j} (v_2)_{i,j-\frac{1}{2}} \, \delta_h \left( \x_{i,j-\frac{1}{2}} - \vchi \left( \X,t \right) \right) h^2, \label{velocity_interpolation_2}
\end{align}
in which $\V_l = \left( V_{l,1}, V_{l,2} \right)$ is the Lagrangian velocity at a Lagrangian marker $\X_l$. We use the notation
\begin{align}
\V = \J \left[ \vchi \left(\cdot,t\right) \right] \v,
\end{align}
in which $\J$ is the discrete velocity restriction or interpolation operator. 
In the present formulation, $\S$ and $\J$ are adjoint operators if evaluated using the same structural configurations \cite{griffith2017hybrid}.

\subsection{Computational algorithm}
We now briefly outline the key steps of the implementation of the IPD method used in the computational examples in Sec.~\ref{s:benchmarks}. Let $\v^n$ and $\vchi^n$ be the fluid velocity and discrete deformation at time $t^n = n \Delta t$, respectively, in which $\Delta t$ is the time step size. We use a second-order time stepping scheme \cite{griffith2017hybrid}, as follows:
\begin{align}
\frac{\vchi^{n+\frac{1}{2}} -\vchi^n}{\dt / 2} &= \J \left[ \vchi^n \right] \v^n, \label{intermediate_lag} \\ 
\frho \left( \frac{\v^{n+1} -\v^n}{\dt} + \N^{\left(n + \frac{1}{2}\right)}\right) &= - \grad_h p^{n + \frac{1}{2}} + \fmu \grad_h^2  \left( \frac{\v^{n+1} + \v^n}{2}\right) + \f^{n + \frac{1}{2}}, \label{discretized_navier_stokes}\\
\grad_h \cdot \v^{n+1} &= 0, \label{discretized_incompressibility} \\
\f^{n + \frac{1}{2}} &=  \S \left[ \vchi^{n + \frac{1}{2}} \right] \F^{n + \frac{1}{2}}, \label{discretized_force_spreading}\\
\frac{\vchi^{n+1} -\vchi^n}{\dt} &= \J \left[ \vchi^{n + \frac{1}{2}} \right] \left( \frac{\v^{n+1} + \v^n}{2} \right) \label{update_lag}, 
\end{align}
in which $\N^{\left(n + \frac{1}{2}\right)} = \frac{3}{2}\v^{n} \cdot \grad_h \v^{n} - \frac{1}{2} \v^{n -1} \cdot \grad_h \v^{n -1}$ is an explicit approximation to the nonlinear advection term. 

Given the initial positions of the Lagrangian markers, the discrete shape tensor is determined by the bond connections between material points in the reference configuration. Then, as the structure moves, the discretized non-local deformation gradient tensor and pairwise bond force function are computed in each time step to account for structural deformations and changes in connectivity (bond breakage). The net bond force is determined by classical constitutive relations with the discrete peridynamic tensors. The discrete Lagrangian force density $\F^{n + \frac{1}{2}}$ is computed at each Lagrangian marker via Eq.~\eqref{discrete_pd_force}. The computational algorithm is performed as follows:
\begin{itemize}[leftmargin=.5in]
\item[Step 1.] Given $\v^n$ and $\vchi^n$. Update $\vchi^{n + \frac{1}{2}}$ using Eq.~\eqref{intermediate_lag}.
\item[Step 2.] Update the non-local deformation gradient tensor and force vector state using Eqs.~\eqref{discrete_nonlocal_deformation_gradient}-\eqref{discrete_pairwise_bond_f}.
\item[Step 2.] Evaluate $\F^{n + \frac{1}{2}}$ using Eq.~\eqref{discrete_pd_force}.
\item[Step 3.] Spread the intermediate Lagrangian force density to the Eulerian grid using Eq.~\eqref{discretized_force_spreading}.
\item[Step 4.] Solve Eqs.~\eqref{discretized_navier_stokes}-\eqref{discretized_incompressibility} for $\v^{n+1}$ and $p^{n+\frac{1}{2}}$.
\item[Step 5.] Update $\vchi^{n+1}$ using Eq.~\eqref{update_lag}.
\end{itemize}

\section{Benchmarks}
\label{s:benchmarks}
We first investigate standard benchmark problems for hyperelastic materials in the conventional solid mechanics literature, with the structure bodies embedded in an incompressible Newtonian fluid for FSI. 
Numerical tests detailed in Sec.~\ref{s:non-failure} exhibit clear constitutive correspondence between the IPD formulation and benchmark FE and FE based IB results.
Afterwards, we focus on fluid-driven material failure through several numerical experiments, as detailed in Sec.~\ref{s:failure}. 

The IBAMR software \cite{griffith2007adaptive, ibamr} is used for all IPD simulations. IBAMR is a distributed-memory parallel implementation of the IB method with support for Cartesian grid adaptive mesh refinement (AMR). 
We compare the results obtained using our new IPD method against an immersed finite element-finite difference (IFED) method \cite{vadala2020stabilization} that is also implemented in IBAMR, and also to a stabilized FE method for incompressible nonlinear elasticity \cite{reese1999new}.

It is well known that conventional IB-type methods can suffer from poor volume conservation \cite{PESKIN199333, griffith_2012}. 
The exact incompressibility condition in the continuous IB formulation, Eq.~\eqref{incompressibility}, can be lost under the spatial discretization, time stepping errors, and the use of regularized delta function kernels. 
To improve volume conservation, we adopt a modified neo-Hookean model unless otherwise mentioned, which was previously shown to improve the volume conservation of the IFED method by Vadala-Roth et al. \cite{vadala2020stabilization}.
In the modified neo-Hookean model, the strain energy and elastic stress are additively decomposed into two parts, isochoric and volumetric,
\begin{align}\label{neo_hookean}
\Psi &= \frac{G}{2}\left( J^{-2/3} \text{tr}\left(\mathbb{C}\right) - 3 \right) + \frac{\kappa_{\text{stab}}}{2} \left( \ln J \right)^2,\\
\PP &= G J^{-2/3} \left( \FF - \frac{\text{tr}\left(\mathbb{C}\right)}{3} \FF^{-T} \right) + \kappa_\mathrm{stab} \ln \left(J\right) \FF^{-T},
\end{align}
in which $G$ is the shear modulus, $J$ is the determinant of non-local deformation gradient tensor, $J = \det \left( \FF \right)$, $\mathbb{C} = \FF^T \FF$ is the right Cauchy-Green strain, and $\kappa_{\mathrm{stab}}$ is the numerical bulk modulus. 
The numerical Poisson's ratio, $\nu_{\mathrm{stab}}$, can be used to define the numerical bulk modulus via
\begin{align}\label{numerical_bulk_modulus}
\kappa_{\mathrm{stab}} = \frac{2G\left(1+\nu_{\mathrm{stab}}\right)}{3\left(1 - 2 \nu_{\mathrm{stab}}\right)}.
\end{align}
This volumetric term reinforces the discrete incompressibility of the immersed structure. 
As in prior work that uses an immersed finite element structural description \cite{vadala2020stabilization}, we test several values of the numerical Poisson's ratio in our benchmarking studies. 

Unless otherwise noted, the density and viscosity of the fluid are respectively set to $\rho = 1.0 \, \frac{\text{g}}{\text{cm}^3}$ and  $\mu = 0.01 \,  \frac{\text{dyn$\cdot$s}}{\text{cm}^2}$, corresponding to water. 
For simplicity, we use the same density for both the structure and fluid. 
The computational domain is $\Omega = [0 , L]^d$, in which $d$ is the spatial dimension and $L$ is the length of domain. 
We define the mesh factor ratio $M_{\mathrm{FAC}} = \frac{\Delta X}{\Delta x}$, in which $\Delta X$  and $\Delta x$ are the Lagrangian and Eulerian grid spacings, respectively. 
The Eulerian grid size is $\Delta x = \frac{L}{N}$, in which $N$ is the number of Cartesian grid cells in one spatial direction. 
We use $M_{\mathrm{FAC}} = 0.5$ in our IPD simulations, so that the structure discretization is twice as fine as the background Cartesian grid.

In our numerical simulations, both static and dynamic versions of benchmarks are considered. 
To efficiently obtain numerical solutions at steady states, the maximum amount of load is applied to the immersed structure using the polynomial $q(t) = -2 \left( \frac{t}{T_{\text{l}}} \right)^3 + 3 \left( \frac{t}{T_{\text{l}}} \right)^2 $, in which $T_{\text{l}} = \alpha T_{\text{f}}$ with $\alpha \in (0,1)$ is a loading time and $T_{\text{f}}$ is a final simulation time. 
In static benchmarks, the final simulation time $T_{\text{f}}$ is determined when the velocity $\V$ is approximately zero.
In addition, viscous damping force is used in the solid region to dampen oscillations and accelerate reaching steady states. Viscous damping is applied to the immersed structure by adding a damping force $-\eta \V$ to the Lagrangian body force $\F$, as in Eq.~\eqref{Lagrangian_multiplier}, in which $\eta > 0$ is the damping coefficient. 

The effect of peridynamic horizon is also investigated in the IPD simulations. 
For the simplicity, a uniform $\horizonsize$-ball is used for the peridynamic horizon. 
The peridynamic horizon size $\horizonsize$ is always taken to be a constant multiple of the Lagrangian mesh spacing $\Delta X$ in the reference configuration, which is commonly used to define the $\horizonsize$-ball in the PD literature \cite{madenci2014peridynamic,behera2020peridynamic,Wang2023}.
Consequently, our grid refinement studies consider the $\horizonsize$-convergence of the IPD formulation\footnote{We remark that the notion of convergence that we call $\horizonsize$-convergence is more commonly called $\delta$-convergence in the peridynamics literature \cite{silling2008convergence,seleson2016convergence}, in which $\delta$ defines the horizon size. As mentioned previously, we avoid using $\delta$ to describe the horizon size since we use $\delta$ to denote the Dirac delta function and $\delta_h$ to denote the regularized delta function.}. 
A larger horizon size implies more interactions between Lagrangian points, and it requires more computations compared to a smaller horizon size. 
Therefore, finding an optimum horizon size is important for optimizing the computational performance of the method. 
Our simulations examine different peridynamic horizon sizes for the constitutive correspondence: $\horizonsize = 1.015  \Delta X, \ 2.015  \Delta X, \ 3.015 \Delta X$.

\subsection{Non-failure benchmarks}
\label{s:non-failure}
This section presents non-failure benchmarks, including standard benchmark problems in solid mechanics literature, using the IPD method to demonstrate the constitutive correspondence to the classical continuum based theory.

\subsubsection{Compression test}
\label{s:Benchmark_Compression}
\begin{figure}[]
\centering
    \includegraphics[width=.45\textwidth]{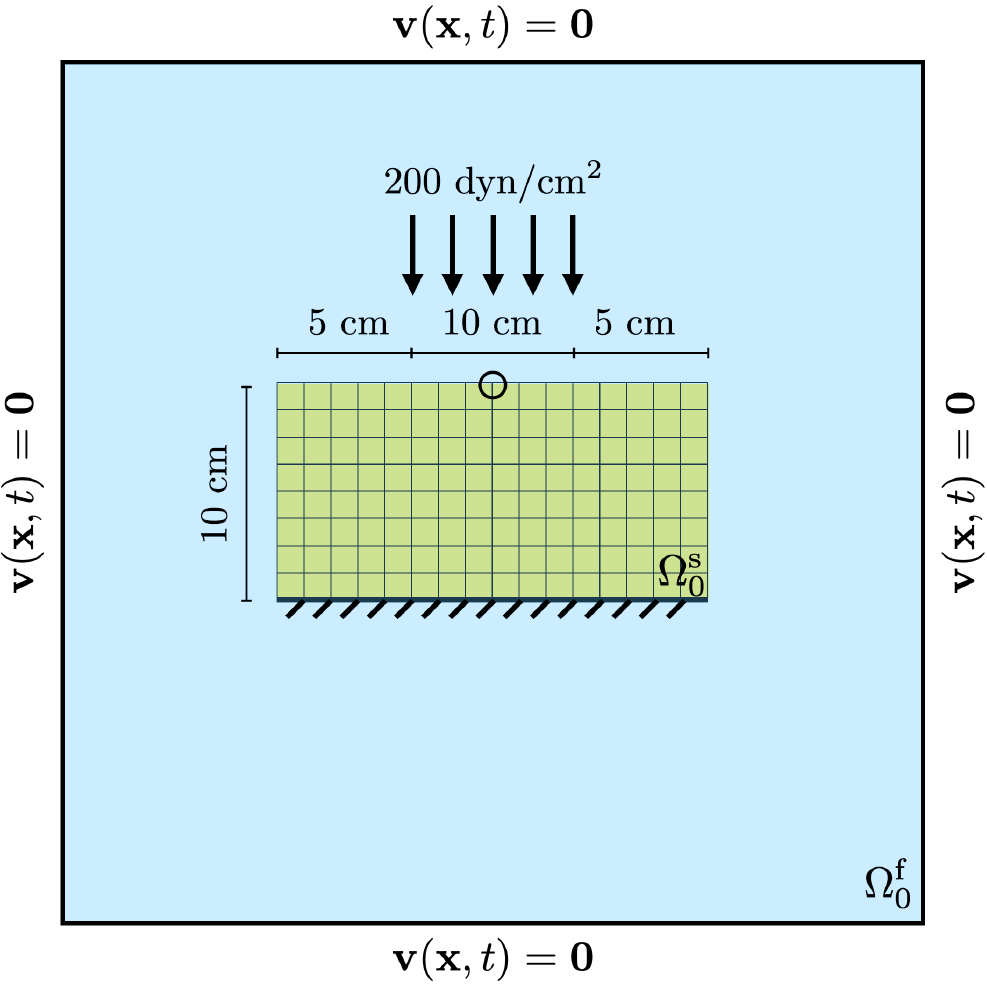}
    \caption{Schematic diagram for the compression test (Sec.~\ref{s:Benchmark_Compression}). The initial configurations of the immersed structure and a fluid are denoted by $\Omega_0^{\text{s}}$ and $\Omega_0^{\text{f}}$, respectively. The entire computational domain is $\Omega = \Omega_0^{\text{s}} \cup \Omega_0^{\text{f}}$. Zero fluid velocity is enforced on the outer boundaries of the computational domain.}
    \label{f:Compression_schematics}
\end{figure}
 We compress a rectangular block to demonstrate a hyperelastic material response under plane strain.
 The computational domain is $\Omega = [0, L]^2$, with $L = 40 \, \text{cm}$. 
 A downward uniaxial traction is loaded in the center of the top of the block. 
 Zero horizontal and vertical displacements are respectively applied to the top and bottom boundaries of the block and all other boundaries have zero traction. 
 This test was introduced by Reese et al.~\cite{reese1999new} to test a stabilization technique for low-order finite elements.
 Fig.~\ref{f:Compression_schematics} provides a schematic of this test case. 
 A shear modulus of $G = 80.194 \, \frac{\text{dyn}}{\text{cm}^2}$ is used for the incompressible neo-Hookean hyperelastic material, and the downward traction is set to $200 \, \frac{\text{dyn}}{\text{cm}^2}$. 
 The load time is $T_{\text{l}} = 100 \, \text{s}$, and the final time is $T_{\text{f}} = 500 \, \text{s}$. 
 An additional damping is set to $\eta = 4.0097 \, \frac{\text{g}}{\text{s}}$.  
 To verify the correspondence to benchmark FE results, material failure (i.e., bond breakage) is not allowed.

 \begin{figure}[]
\centering
	\begin{tabular}{cc}
	 \includegraphics[width=.95\textwidth]{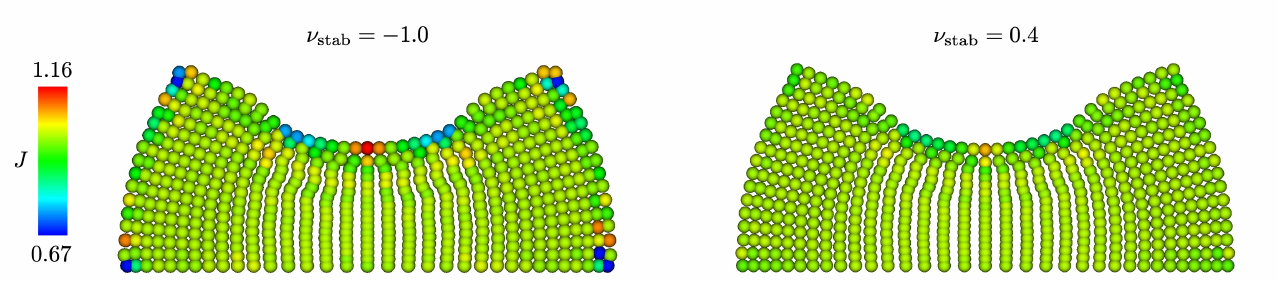}
	\end{tabular}
	    \caption{Deformations of the hyperelastic block along with the values of $J$ at material points using the neo-Hookean material model with $G = 80.194 \, \frac{\text{dyn}}{\text{cm}^2}$. The deformations are computed using 561 solid degrees of freedom (DoF) and $\horizonsize = 2.015 \Delta X$. The left panel shows the deformation obtained using $\nu_{\mathrm{stab}}  = -1.0$, and the right panel shows the result for $\nu_{\mathrm{stab}}  = 0.4$.}
    \label{f:Compression_deformation}
\end{figure}

\begin{figure}[]
\centering
    \begin{tabular}{cc}
        \begin{subfigure}{.3\textwidth}
          		\includegraphics[width=\textwidth]{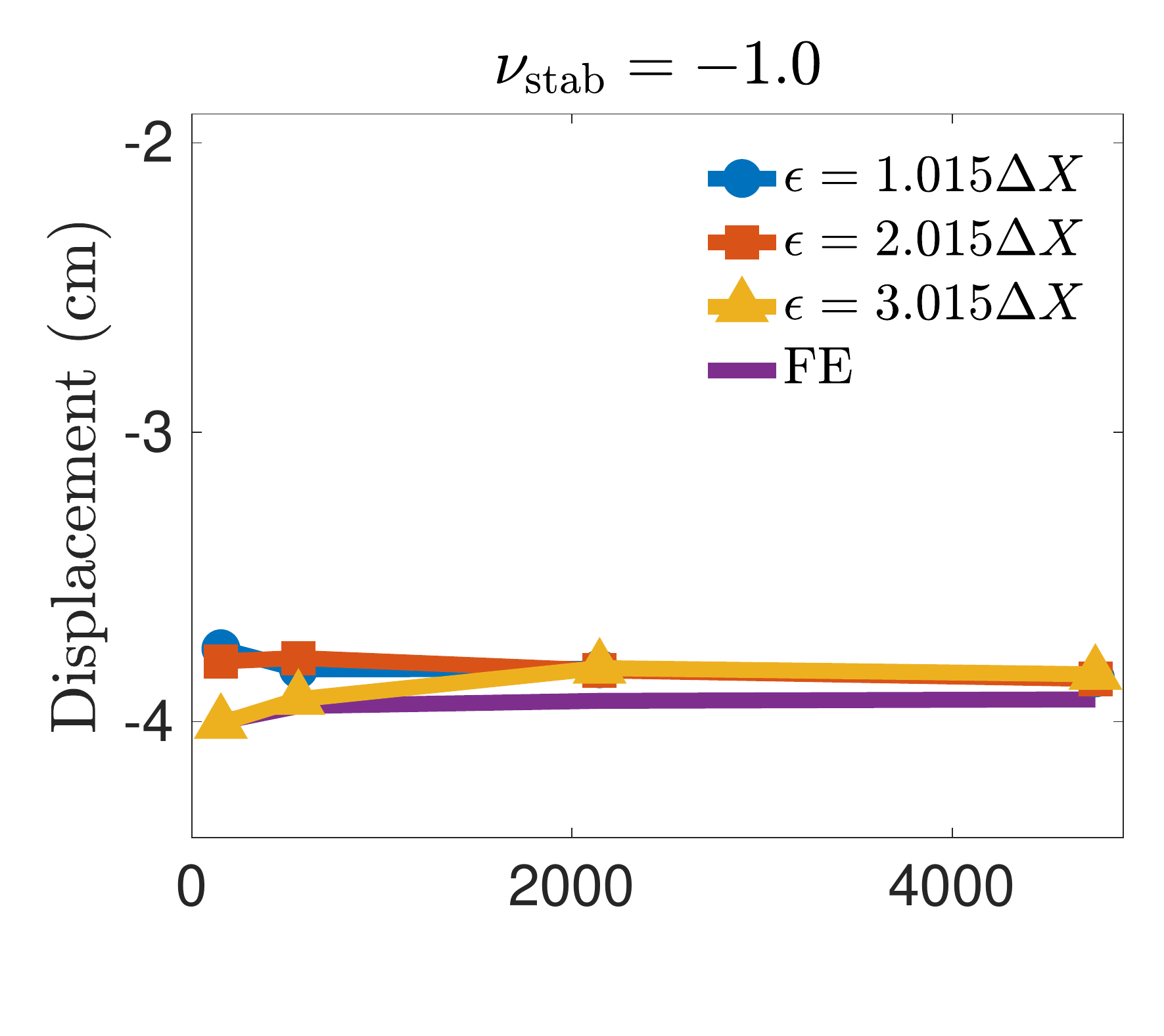}
        \end{subfigure} 
        \begin{subfigure}{.3\textwidth}
                \includegraphics[width=\textwidth]{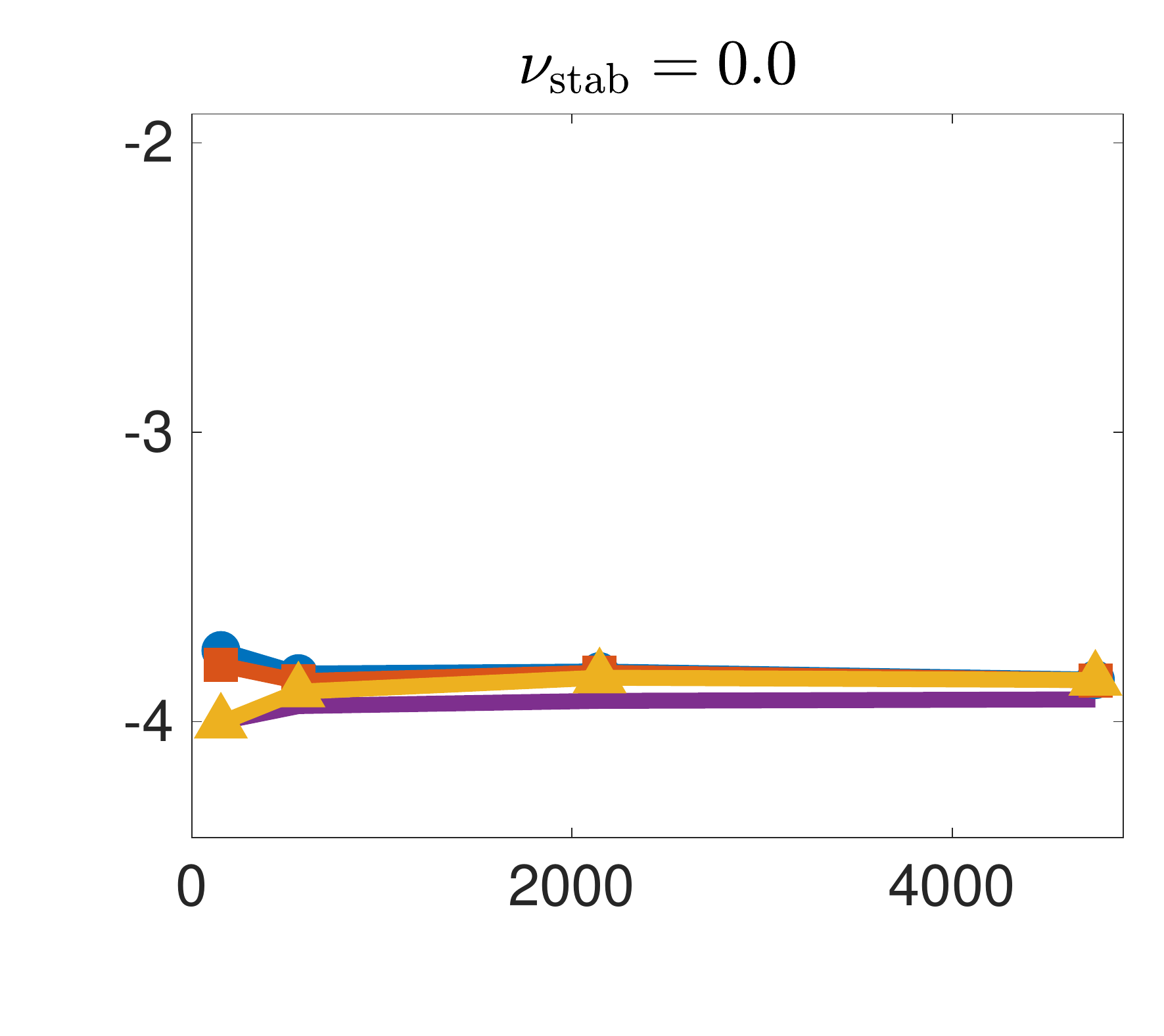}
        \end{subfigure}
        \begin{subfigure}{.3\textwidth}
                \includegraphics[width=\textwidth]{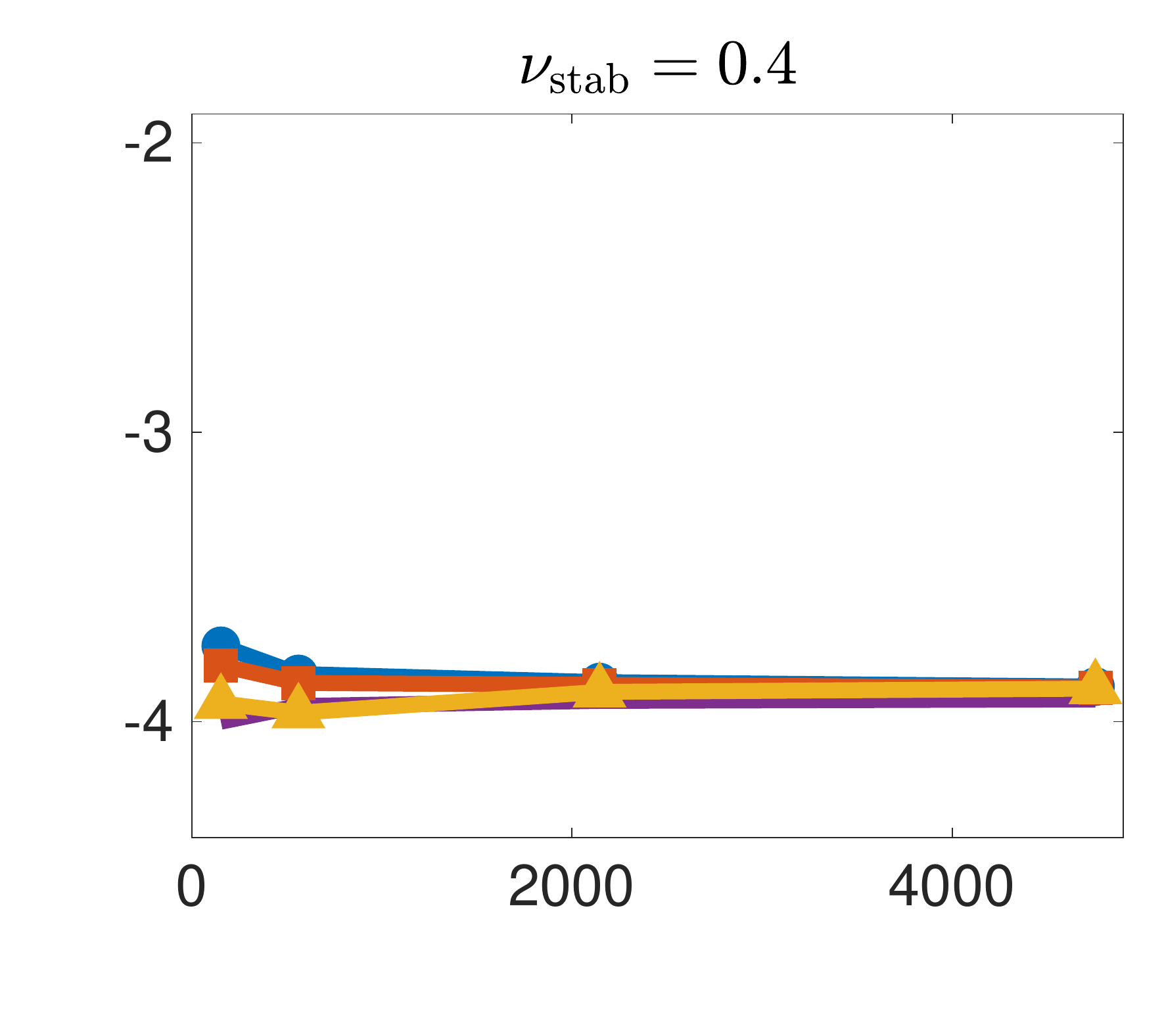}
        \end{subfigure} \\
         \begin{subfigure}{.3\textwidth}
          		\includegraphics[width=\textwidth]{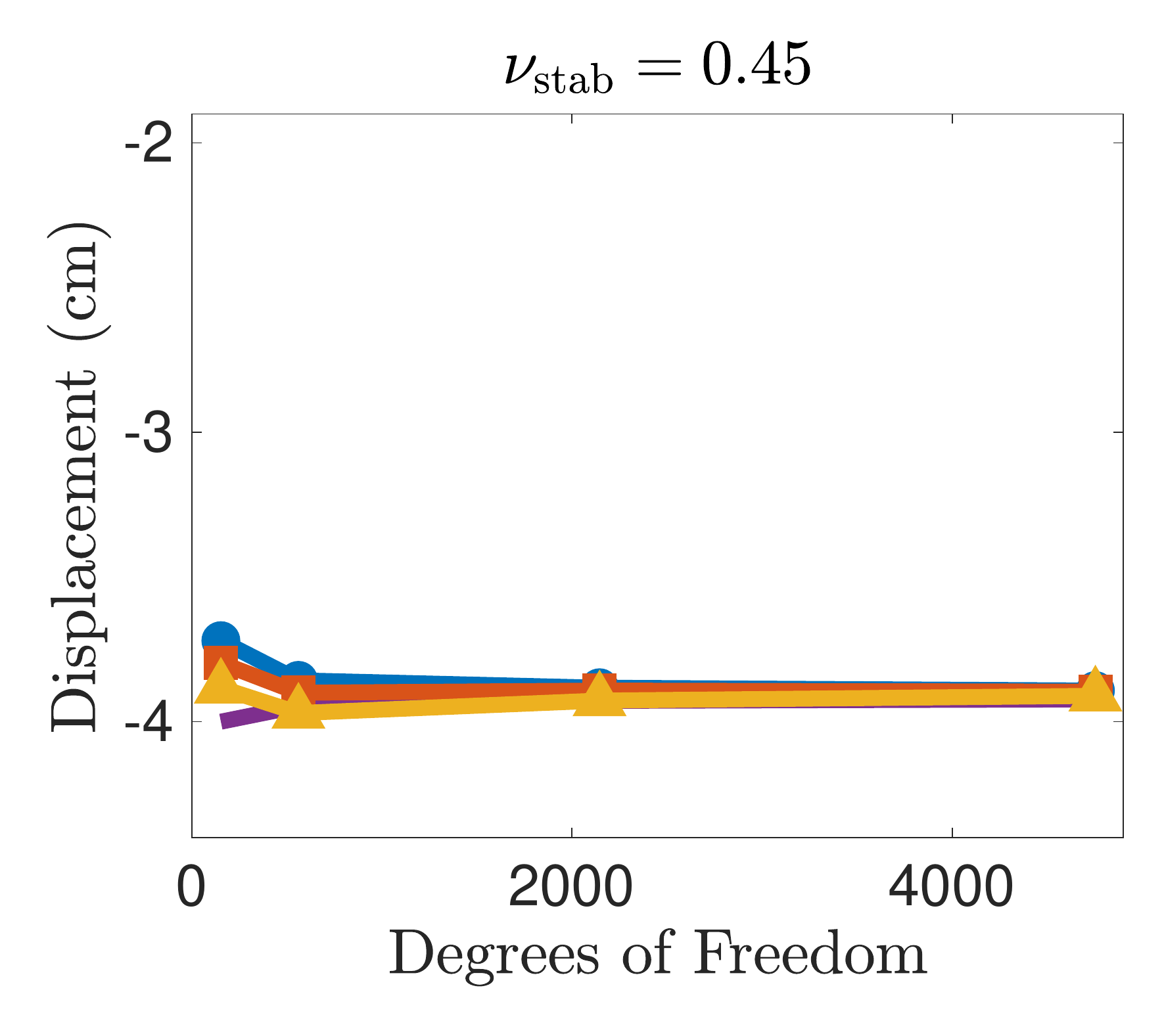}
        \end{subfigure} 
        \begin{subfigure}{.3\textwidth}
                \includegraphics[width=\textwidth]{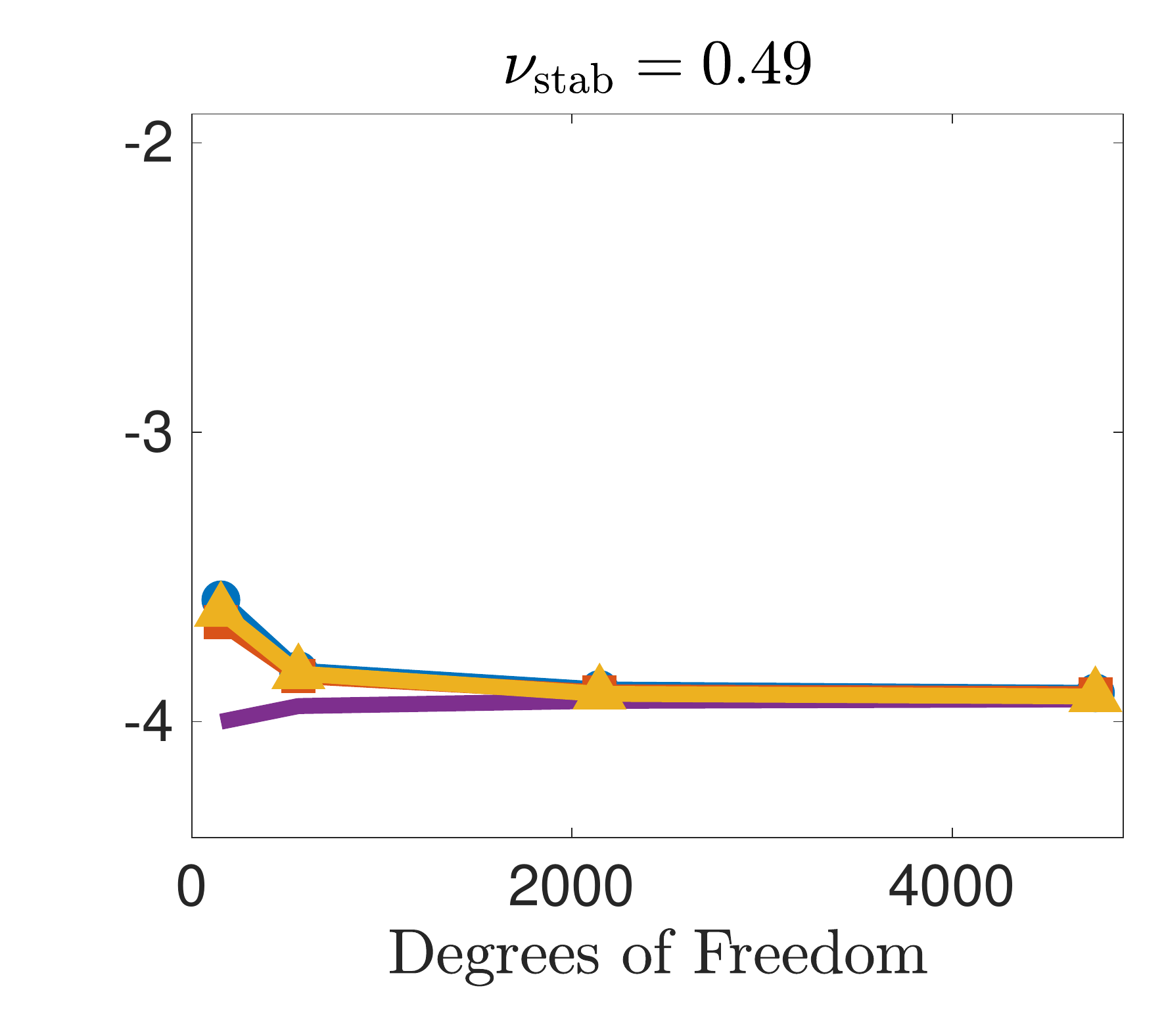}
        \end{subfigure}
        \begin{subfigure}{.3\textwidth}
                \includegraphics[width=\textwidth]{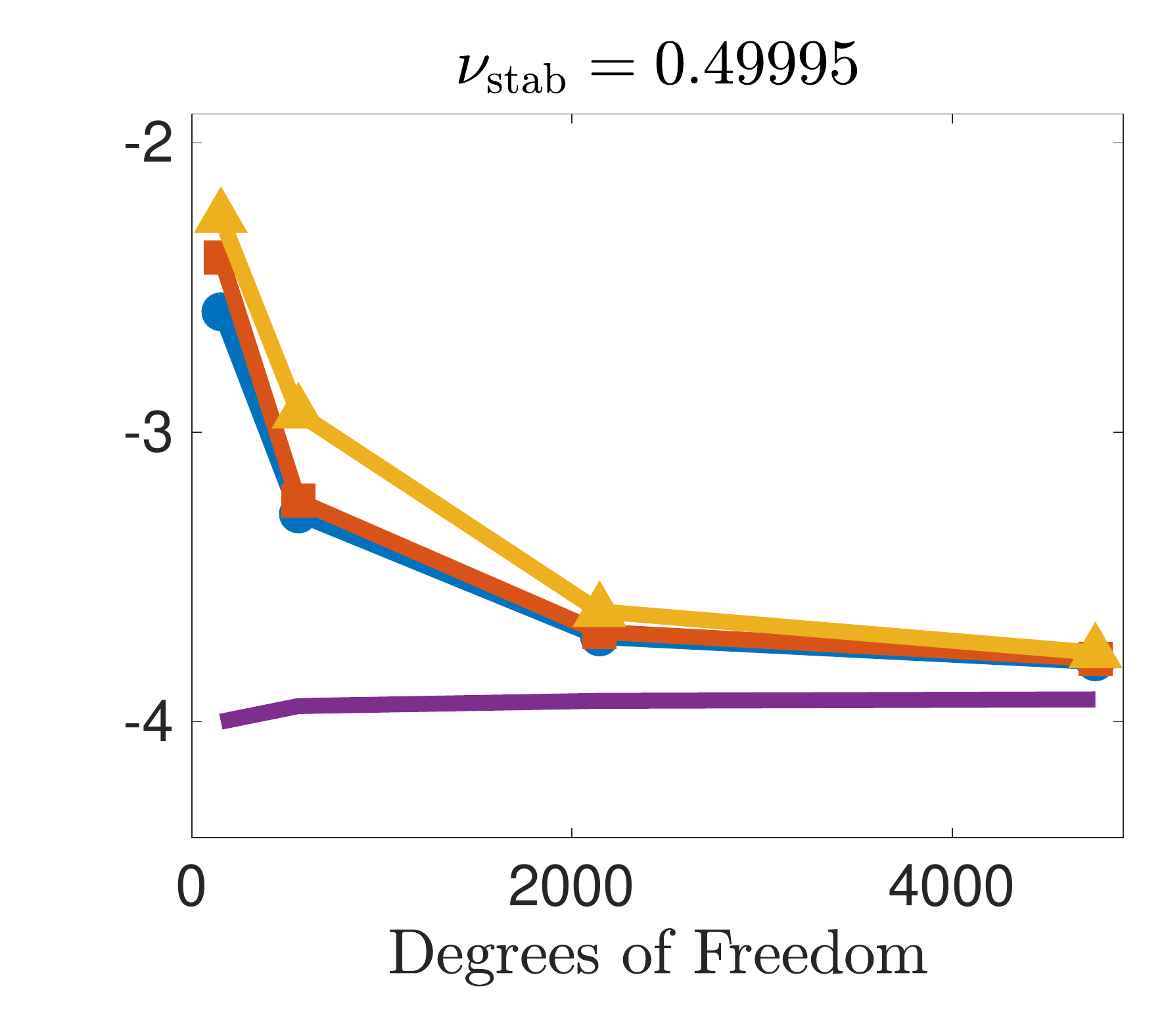}
        \end{subfigure} \\
    \end{tabular}
    \caption{Vertical displacements of the top center point of the compressed block, highlighted in Fig.~\ref{f:Compression_schematics}, for different choices of peridynamic horizon size $\horizonsize$ and numerical Poisson's ratio $\nu_{\mathrm{stab}}$. The solid DoF range from $153$ to $4753$. Note that locking clearly occurs for $\nu_{\mathrm{stab}} = 0.49995$. As in standard computational mechanics approaches, however, the IPD formulation ultimately converges under grid refinement even with high (but fixed) levels of volumetric penalization.}
    \label{f:Compression_disp}
\end{figure}

Fig.~\ref{f:Compression_deformation} illustrates the material body after the deformation along with pointwise values of the non-local Jacobian determinant $J$, which is evaluated from the non-local deformation gradient tensor.
Fig.~\ref{f:Compression_disp} shows the vertical displacements of the top center material point, highlighted in Fig.~\ref{f:Compression_schematics}, for various numerical Poisson's ratios $\nu_{\mathrm{stab}}$ and peridynamic horizon sizes $\horizonsize$ under grid refinement. 
The maximum displacement of the point obtained using IPD method is in excellent agreement with that obtained using the standard FE method, and it converges under grid refinement  to approximately $3.92 \, \text{cm}$.
The maximum displacement of the point of interest is relatively small (between $2.25 \, \text{cm}$ and $3.30 \, \text{cm}$) at low grid resolutions if a larger value of $\nu_{\mathrm{stab}}  = 0.49995$ is used. 
In the computational mechanics literature, this issue is referred to as volumetric locking and can occur with large values of the volumetric penalty parameters.
Note that $\kappa_{\mathrm{stab}}  \rightarrow \infty$ as $\nu_{\mathrm{stab}} \rightarrow 0.5$. 
Under grid refinement, we ultimately recover accurate deformations for fixed finite values of $\kappa_{\text{stab}}$, as in standard methods for nearly incompressible elasticity.

\begin{figure}[]
\centering
   \begin{tabular}{cc}
        \begin{subfigure}{.3\textwidth}
          		\includegraphics[width=\textwidth]{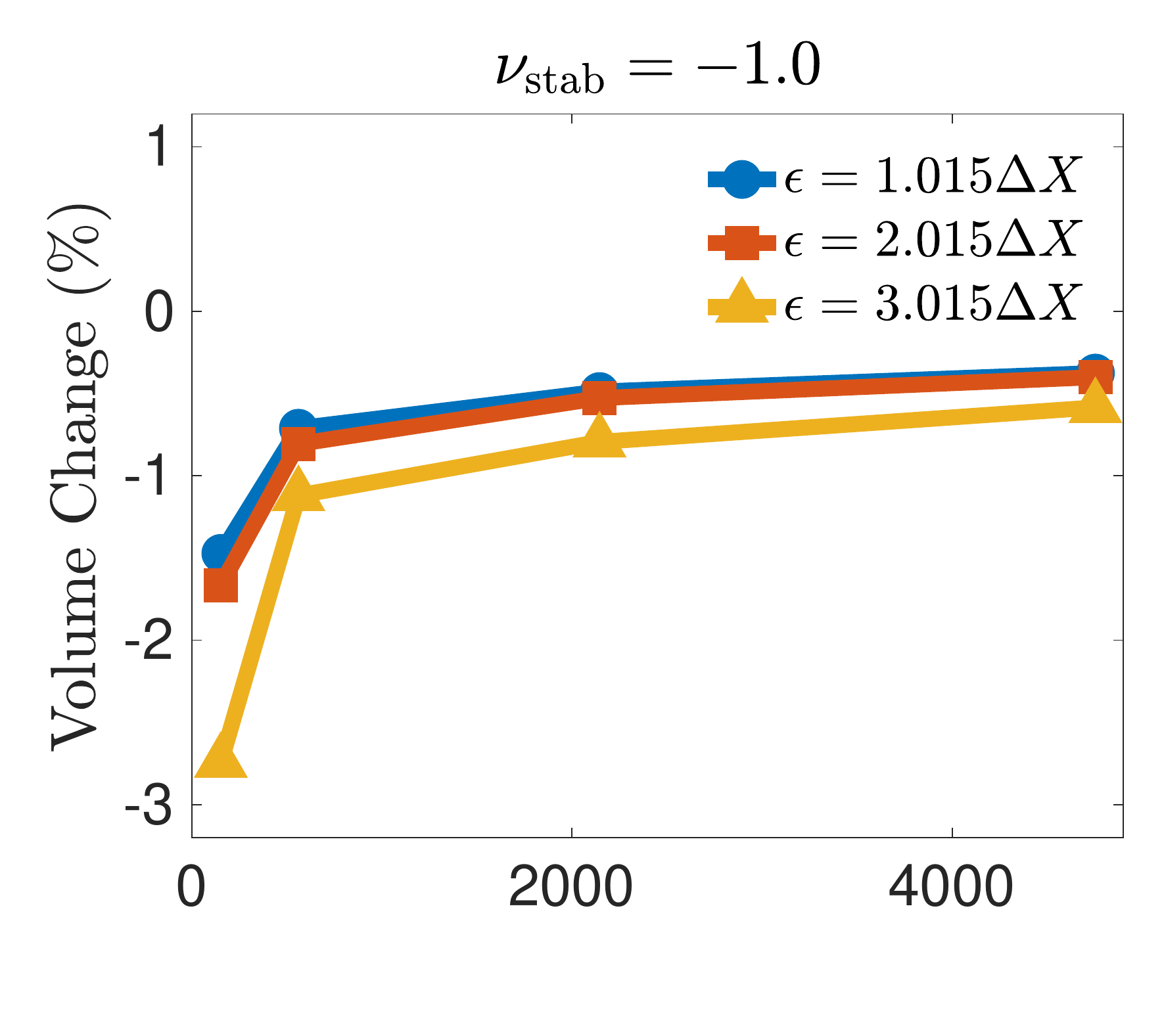}
        \end{subfigure} 
        \begin{subfigure}{.3\textwidth}
                \includegraphics[width=\textwidth]{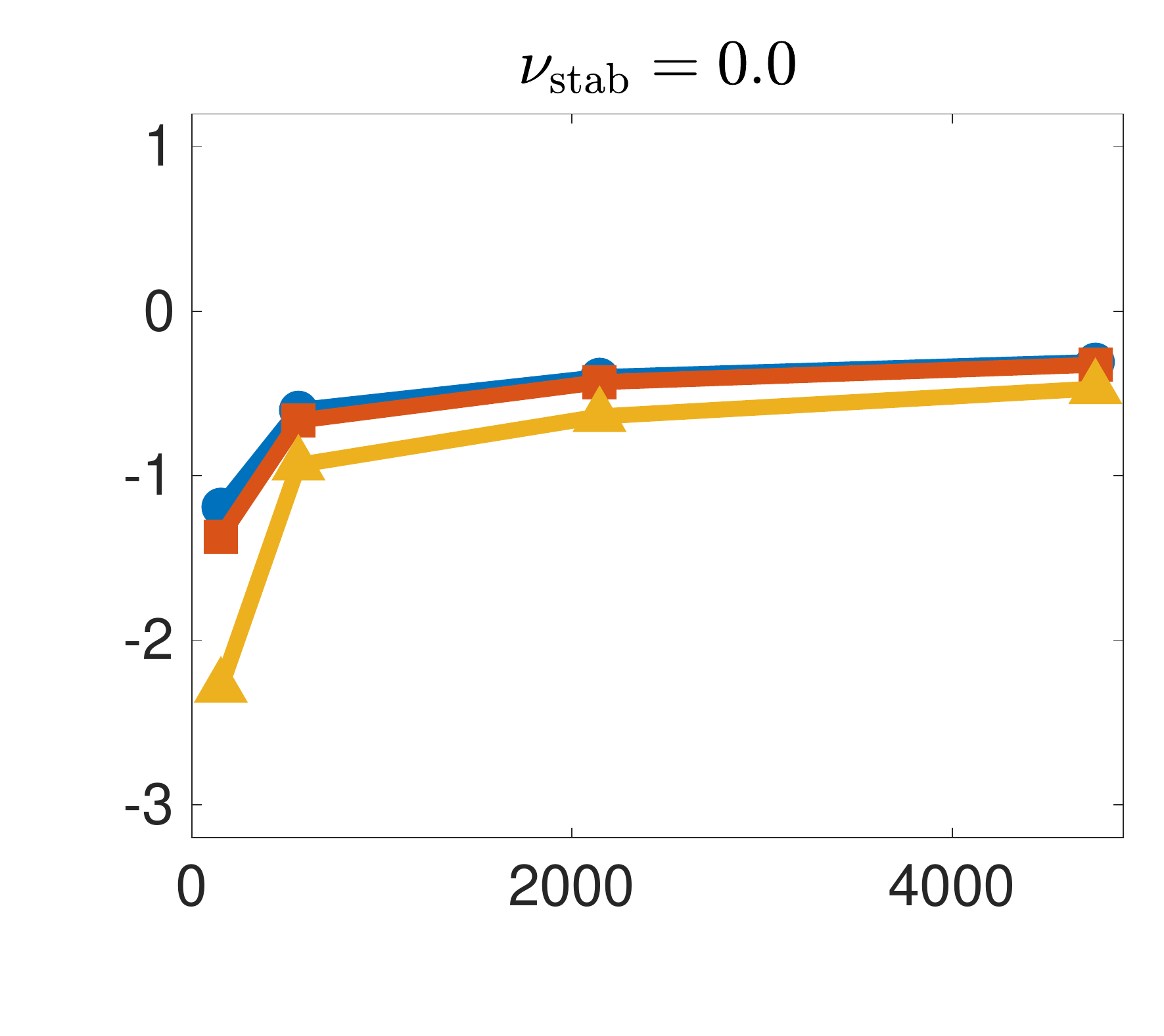}
        \end{subfigure}
        \begin{subfigure}{.3\textwidth}
                \includegraphics[width=\textwidth]{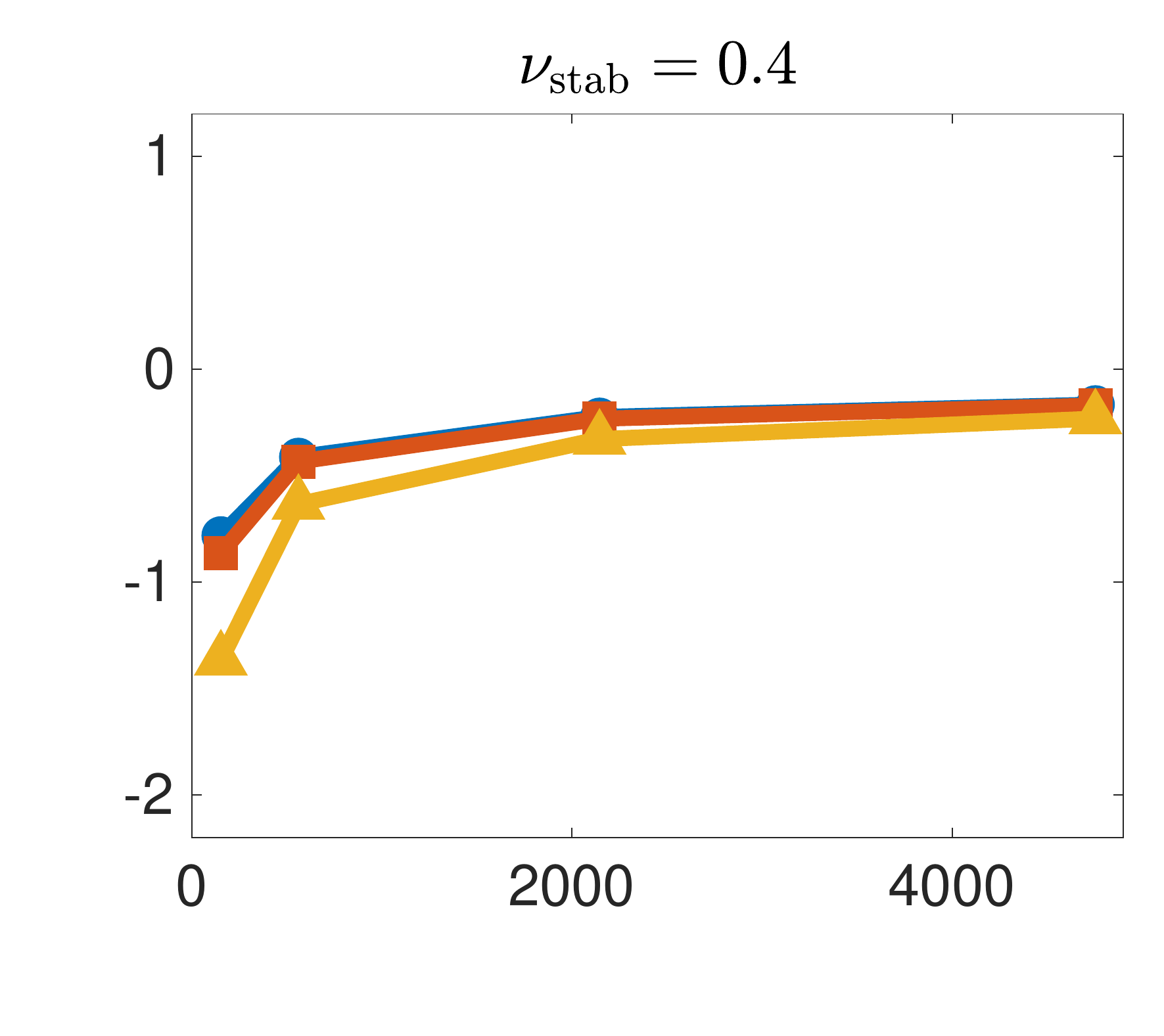}
        \end{subfigure} \\
         \begin{subfigure}{.3\textwidth}
          		\includegraphics[width=\textwidth]{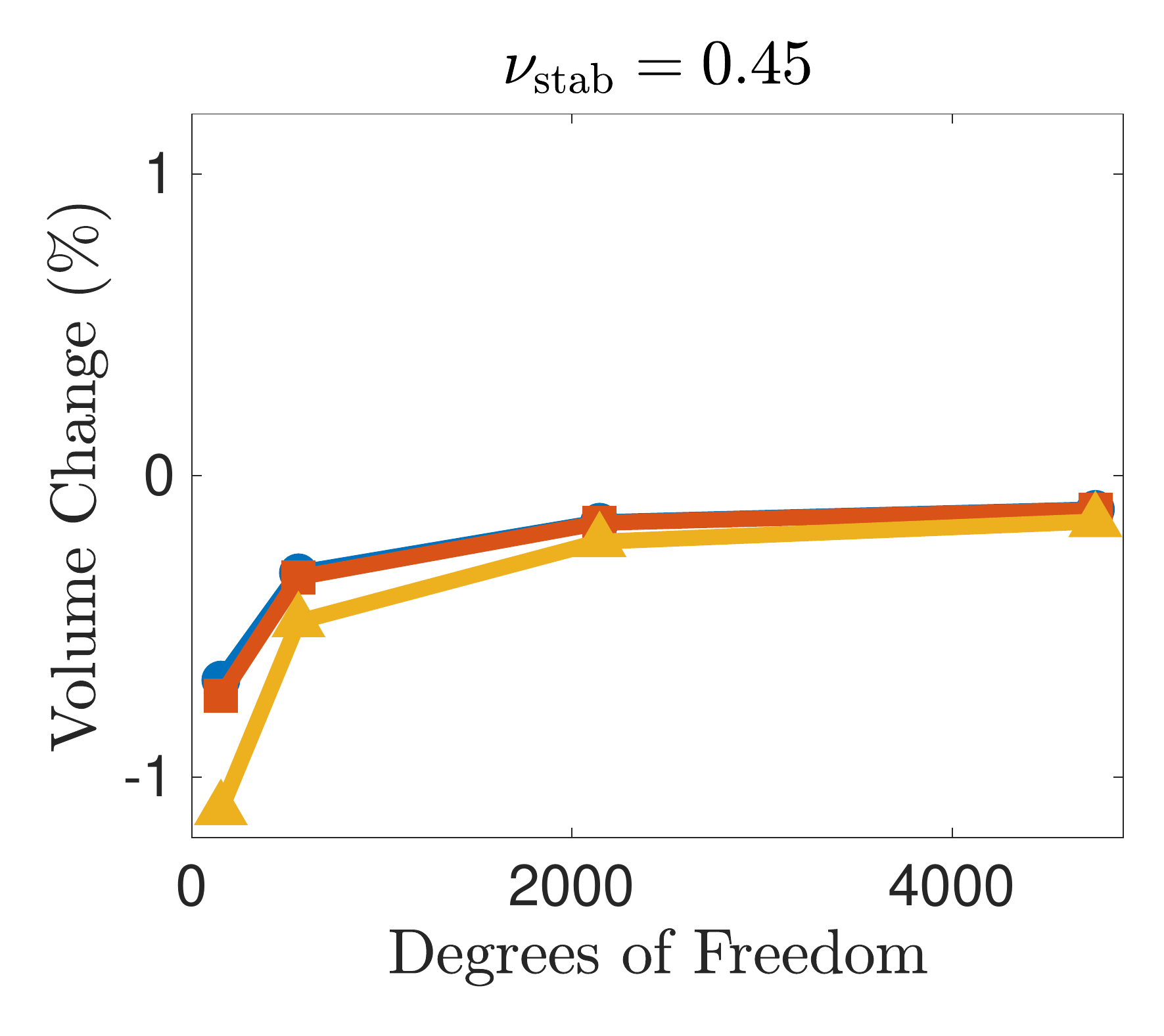}
        \end{subfigure} 
        \begin{subfigure}{.3\textwidth}
                \includegraphics[width=\textwidth]{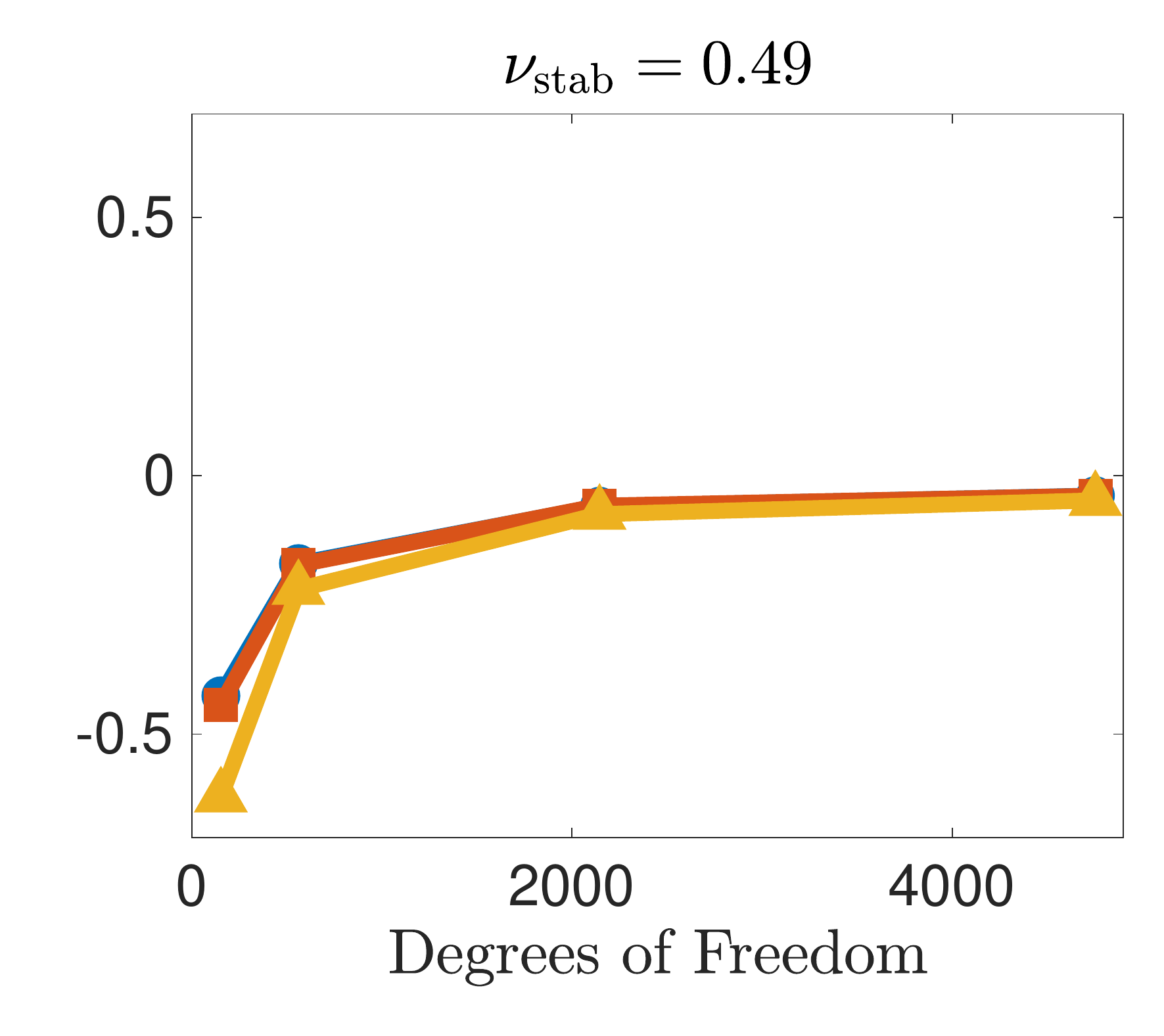}
        \end{subfigure}
        \begin{subfigure}{.3\textwidth}
                \includegraphics[width=\textwidth]{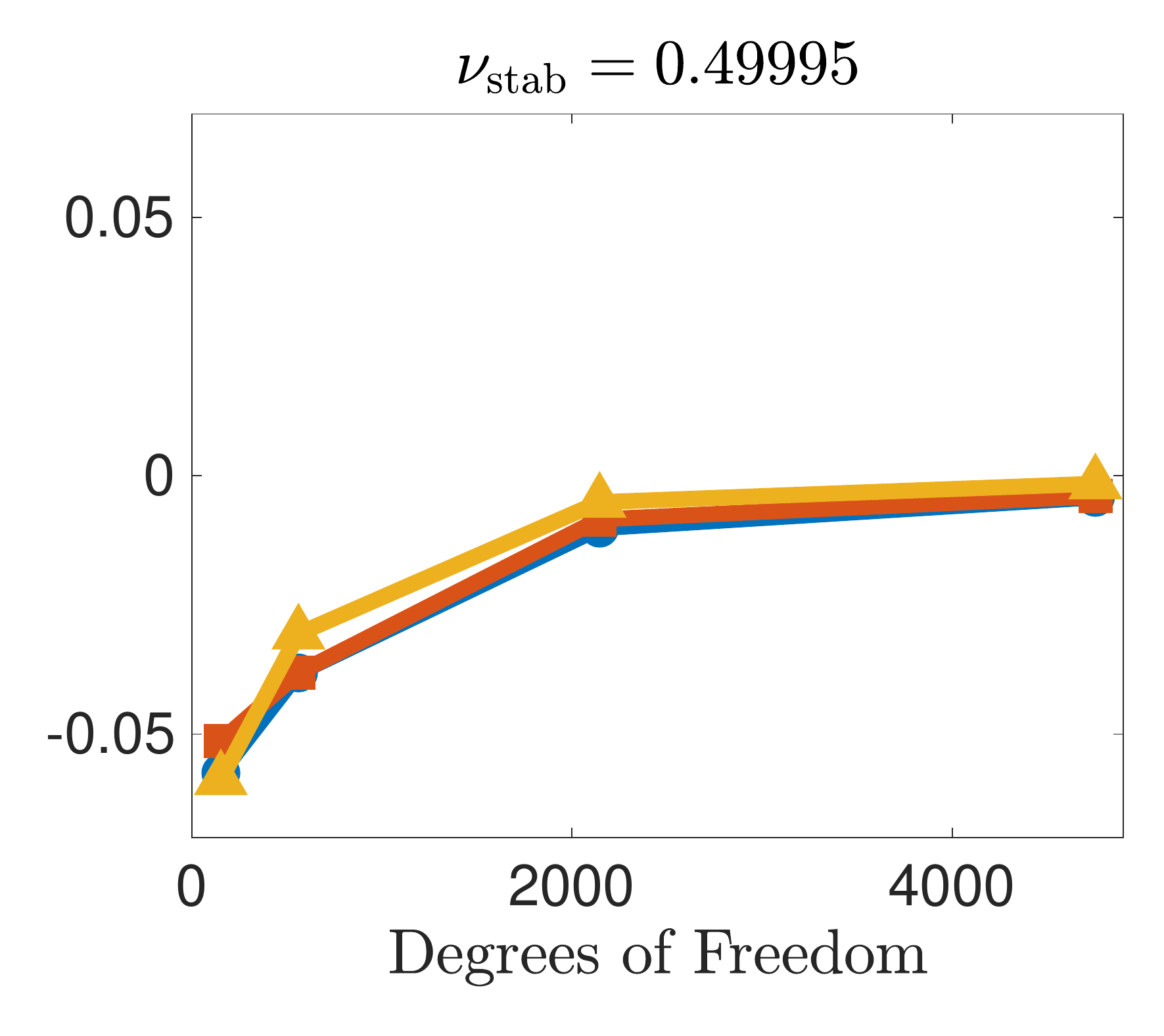}
        \end{subfigure} \\
    \end{tabular}
    \caption{Volume change of the compressed block for different choices of horizon size $\horizonsize$ and numerical Poisson's ratio $\nu_{\mathrm{stab}}$. The solid DoF range from $153$ to $4753$. The largest change is approximately $2.7 \%$.}
    \label{f:Compression_vol}
\end{figure}

Fig.~\ref{f:Compression_vol} shows the volume change observed under deformation for different grid spacings. 
If $\nu_{\mathrm{stab}}$ is small, slight volumetric changes occur (between $0.3\%$ and $2.7\%$) under loading. 
This volume change becomes negligible (up to $0.001 \%$) when larger values of $\nu_{\mathrm{stab}} \ge 0.4$ are used. 
This is also clear in Fig.~\ref{f:Compression_deformation}. 
IPD results agree with results obtained using IFED, with both methods exhibiting similar volume changes that range between $0.0004\%$ and $2.1\%$. 
Under grid refinement, negligible spurious volume changes or locking occurs in all IPD simulations. 
In addition, relatively consistent results are obtained for all choices of the PD horizon sizes considered in the tests.

\subsubsection{Cook's membrane}
\label{s:Benchmark_Cooks}
\begin{figure}[]
\centering
    \includegraphics[width=.45\textwidth]{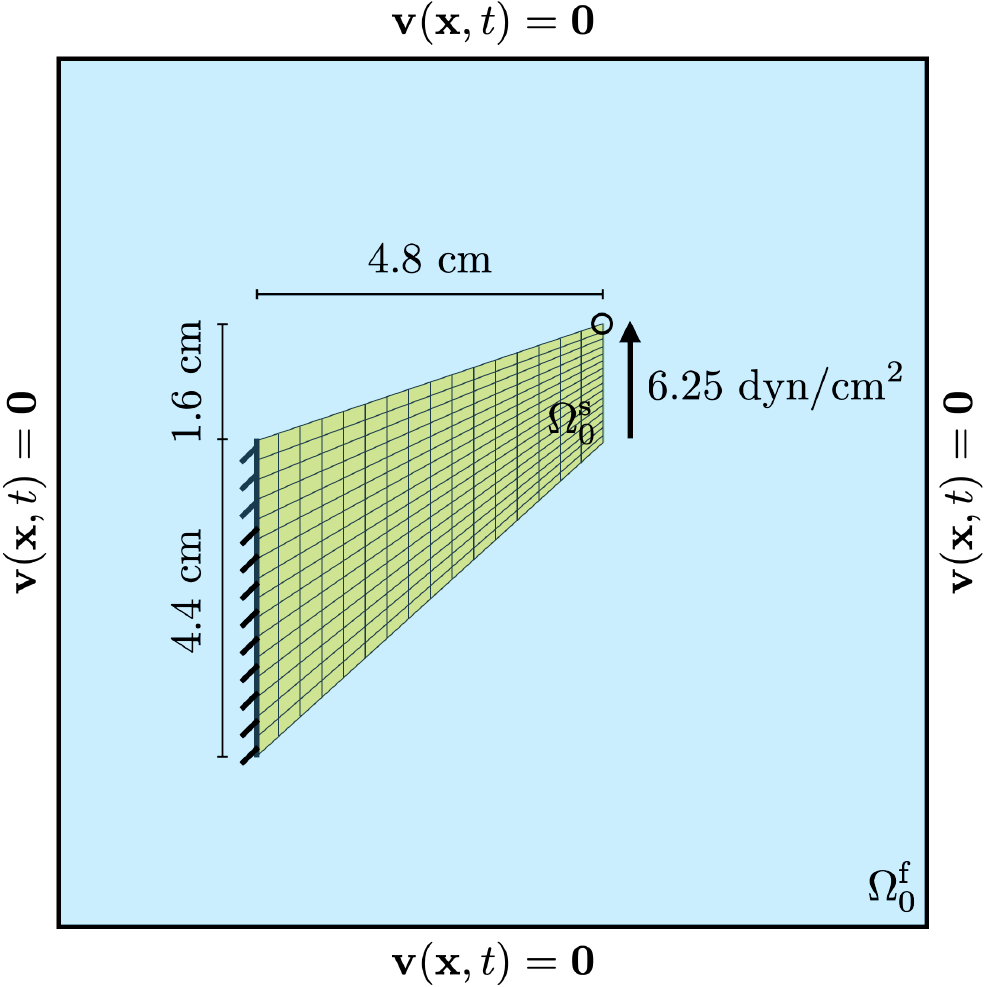}
    \caption{Schematic diagram for the Cook's membrane benchmark (Sec.~\ref{s:Benchmark_Cooks}). The initial configurations of the immersed structure and a fluid are denoted by $\Omega_0^{\text{s}}$ and $\Omega_0^{\text{f}}$, respectively. The entire computational domain is $\Omega = \Omega_0^{\text{s}} \cup \Omega_0^{\text{f}}$. Zero fluid velocity is enforced on the other boundaries of the computational domain.}
    \label{f:Cooks_schematics}
\end{figure}

Cook's membrane \cite{cook1974improved}, which is another widely used plane strain problem, is used to demonstrate the hyperelastic material response under bending and shearing.
The computational domain is $\Omega = [0, L]^2$, with $L = 40 \,  \text{cm}$. 
Zero displacement is imposed on the left boundary of the structure, and an upward traction of $6.25 \, \frac{\text{dyn}}{\text{cm}^2}$ is applied to the right boundary. 
Otherwise, stress-free boundary conditions are assumed.  
Fig.~\ref{f:Cooks_schematics} provides a schematic of this test case. 
A shear modulus of $G = 83.3333 \, \frac{\text{dyn}}{\text{cm}^2}$ is used for the incompressible hyperelastic material. 
The load time is $T_{\text{l}} = 20 \, \text{s}$, the final time is $T_{\text{f}} = 50 \, \text{s}$, and the damping parameter is set to $\eta = 4.16667 \,  \frac{\text{g}}{\text{s}}$. 
We focus on the vertical displacements of the top-right corner of the membrane to assess convergence. 
To impose the same volume fraction at each material point in the NOSB-PD formulation, we use a stair-step representation of the immersed membrane. 
Consequently, the horizon size needs to be large enough ($\horizonsize \ge \sqrt{2} \Delta X$) to ensure adequate bond connectivity throughout the material.
Bond breakage is not considered for this problem as well.

 \begin{figure}[]
\centering
	\includegraphics[width=.85\textwidth]{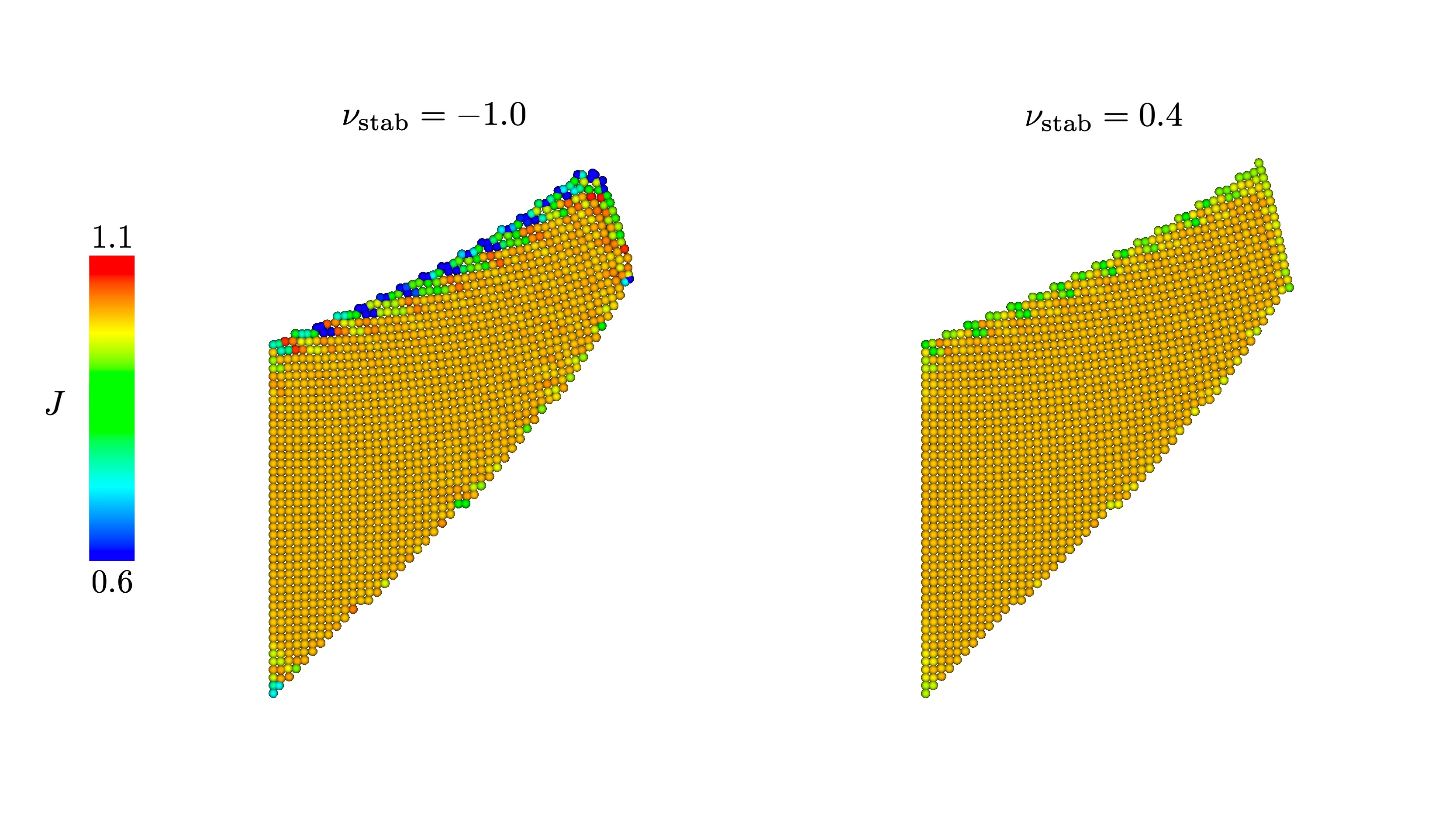} \vspace{-.3in}
    \caption{Deformations of Cook's membrane with the values of $J$ at material points using the neo-Hookean material model with $G = 83.3333 \, \frac{\text{dyn}}{\text{cm}^2}$. The deformations are represented using $1481$ solid DoF and $\horizonsize = 2.015 \Delta X$.  The left panel shows the deformation obtained using $\nu_{\mathrm{stab}}  = -1.0$, and the right panel shows the result for $\nu_{\mathrm{stab}}  = 0.4$.}
    \label{f:Cooks_deformation}
\end{figure}

\begin{figure}[]
\centering
    \begin{tabular}{cc}
        \begin{subfigure}{.3\textwidth}
          		\includegraphics[width=\textwidth]{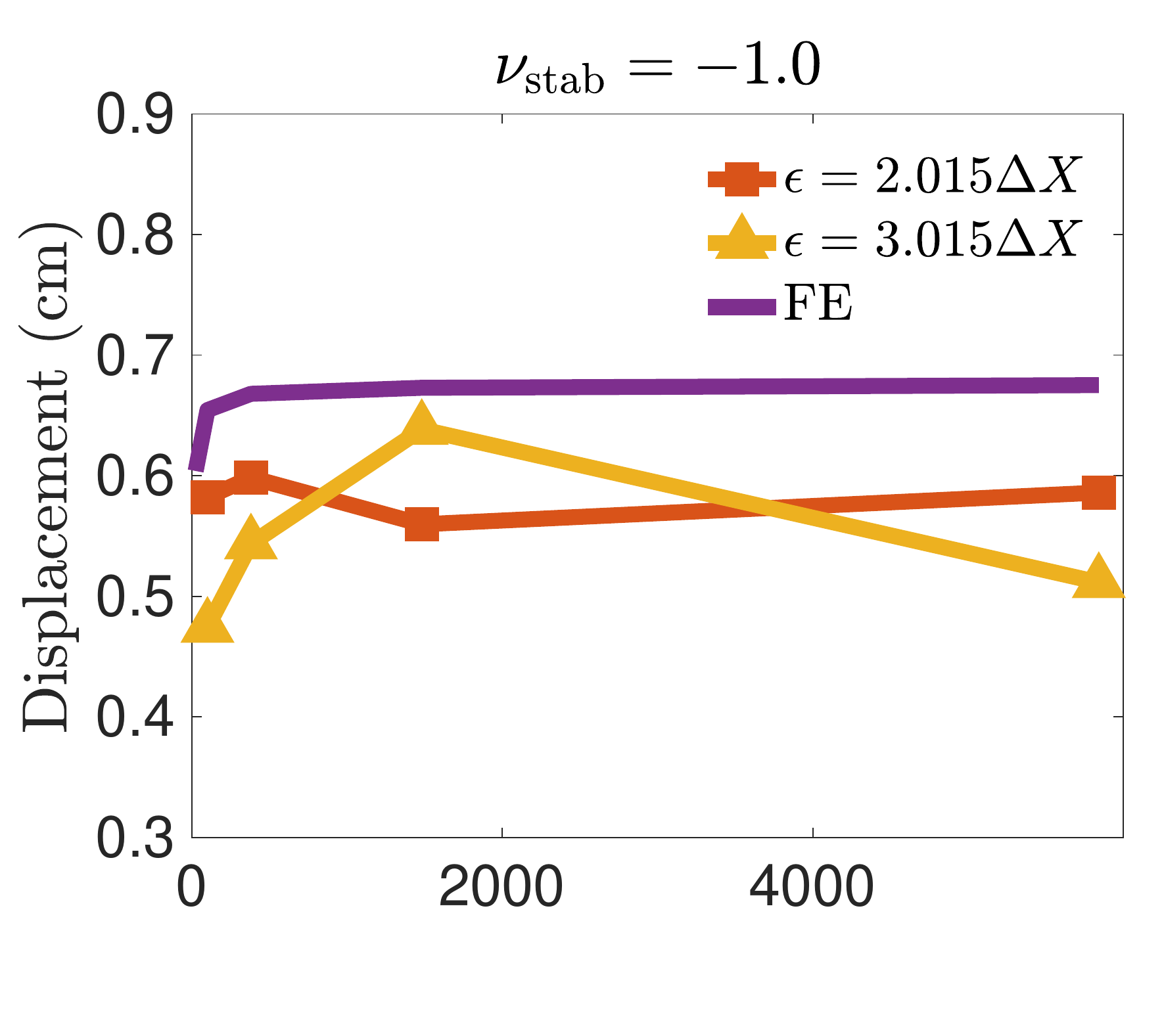}
        \end{subfigure} 
        \begin{subfigure}{.3\textwidth}
                \includegraphics[width=\textwidth]{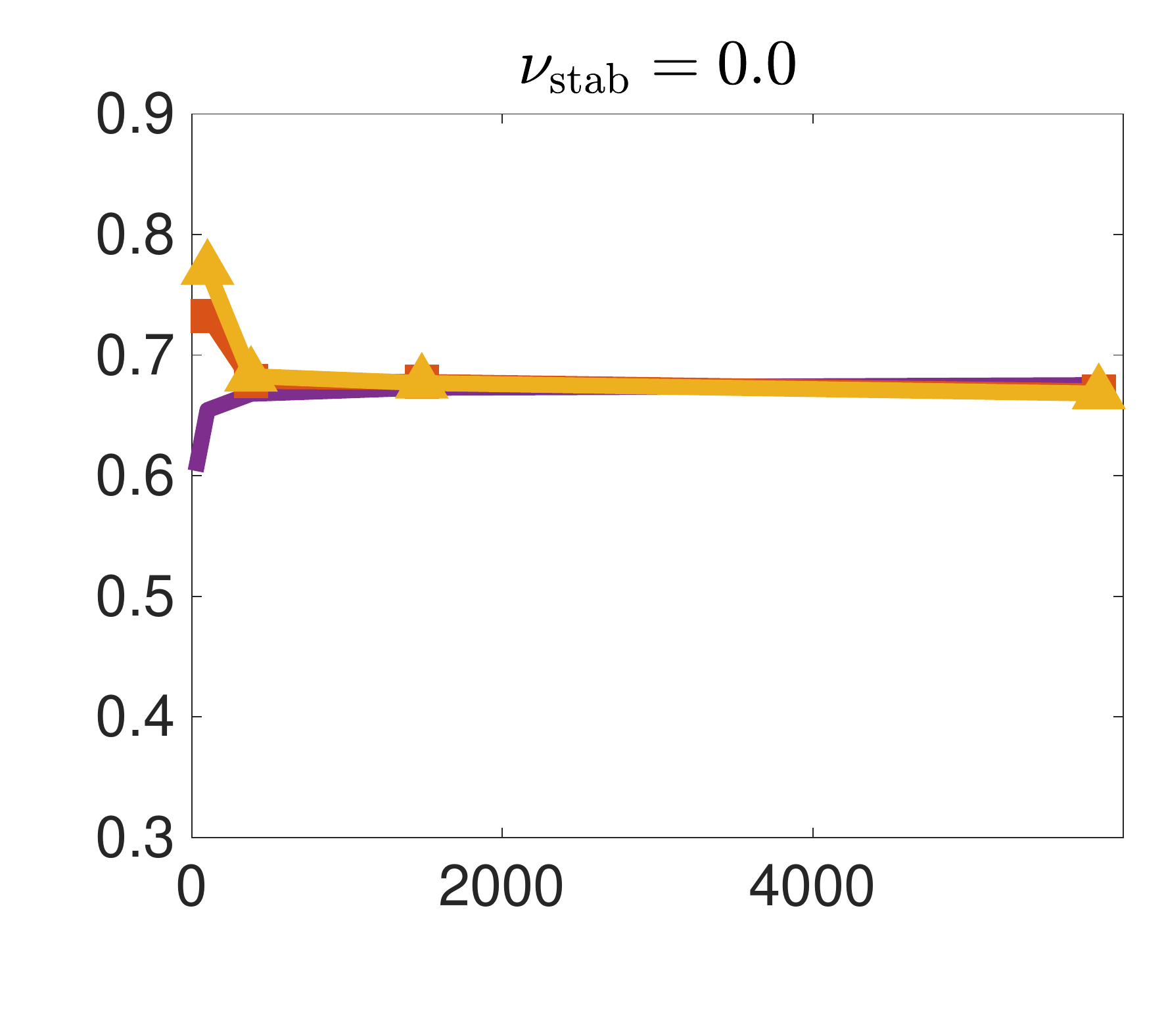}
        \end{subfigure}
        \begin{subfigure}{.3\textwidth}
                \includegraphics[width=\textwidth]{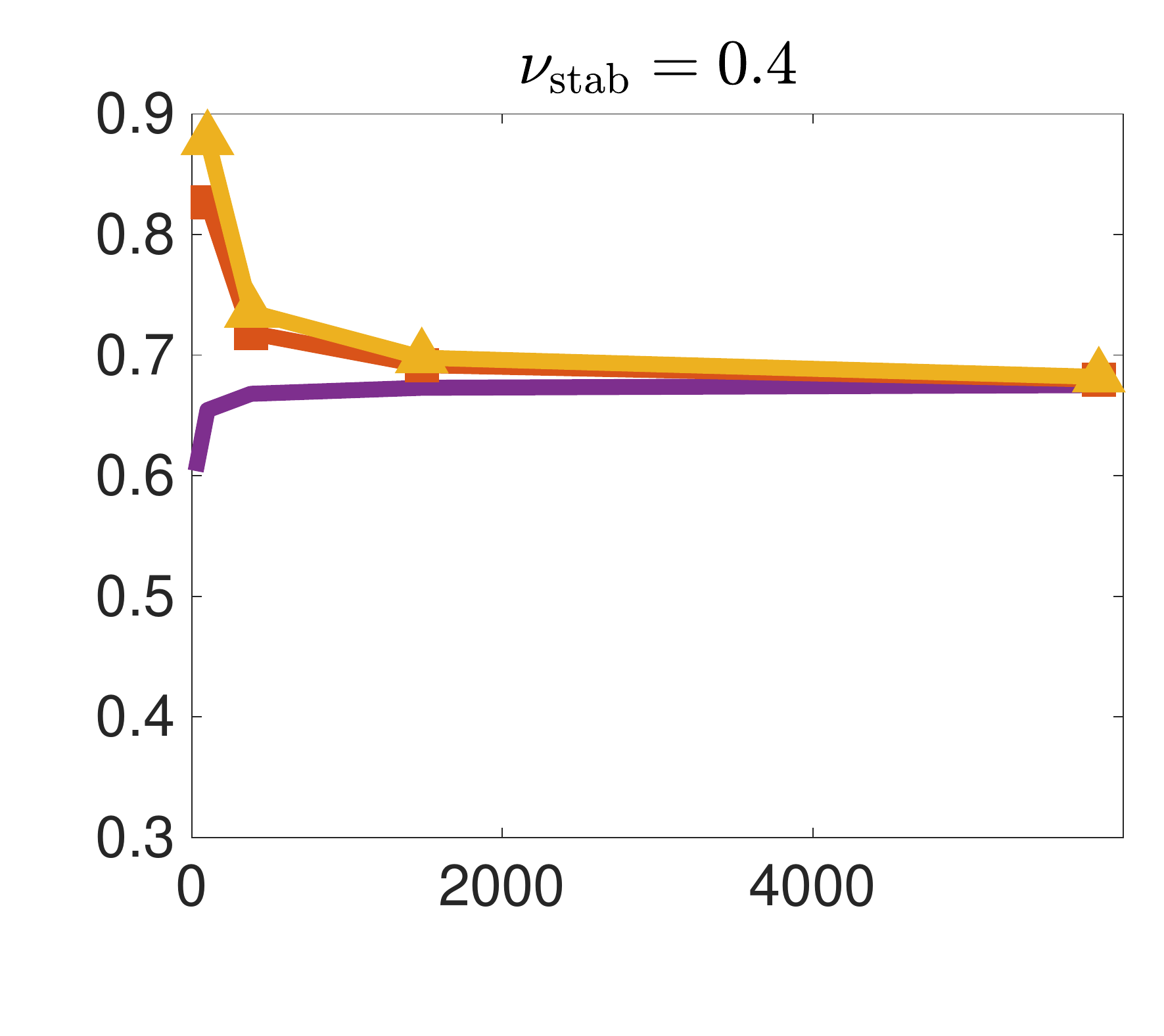}
        \end{subfigure} \\
         \begin{subfigure}{.3\textwidth}
          		\includegraphics[width=\textwidth]{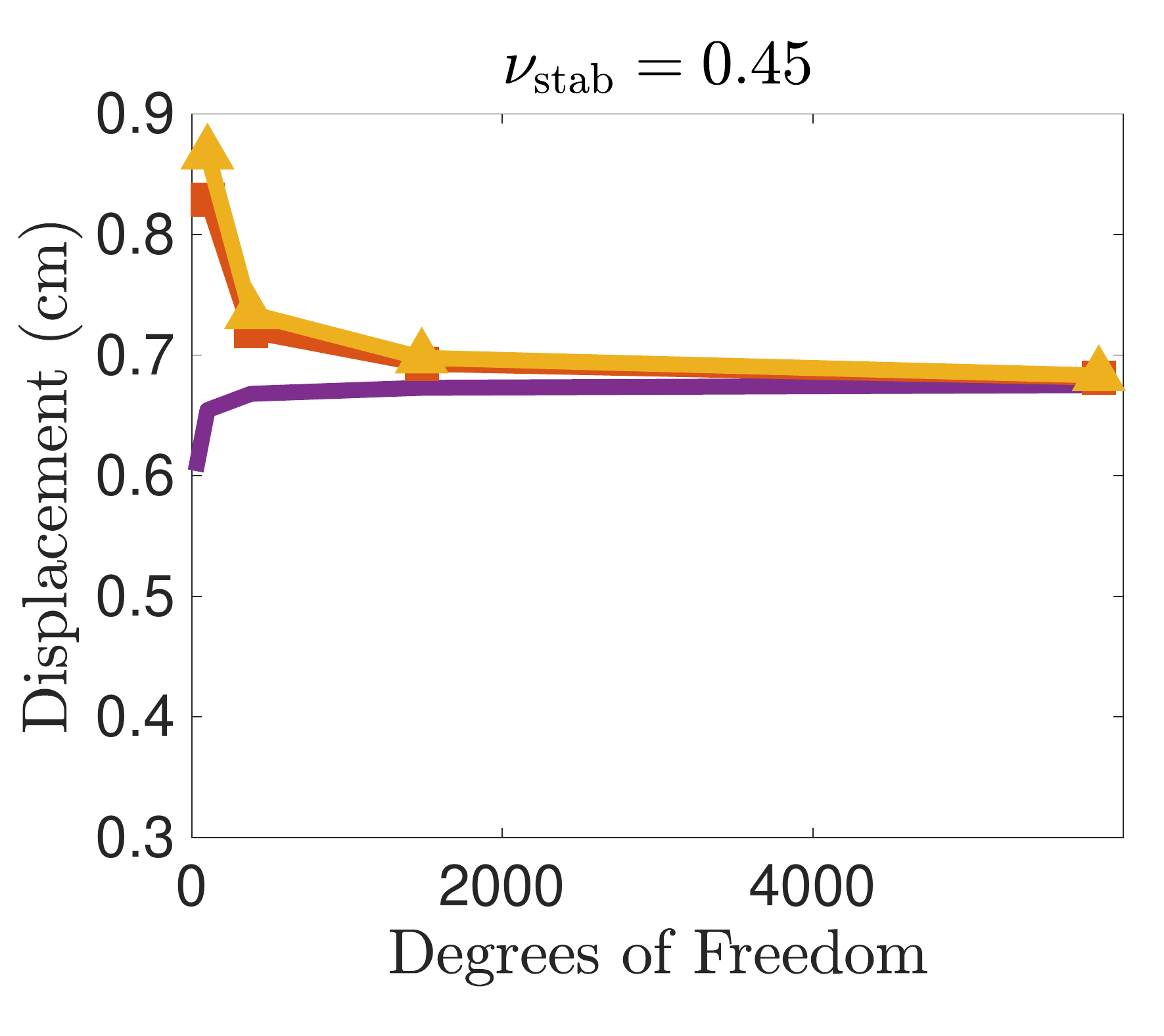}
        \end{subfigure} 
        \begin{subfigure}{.3\textwidth}
                \includegraphics[width=\textwidth]{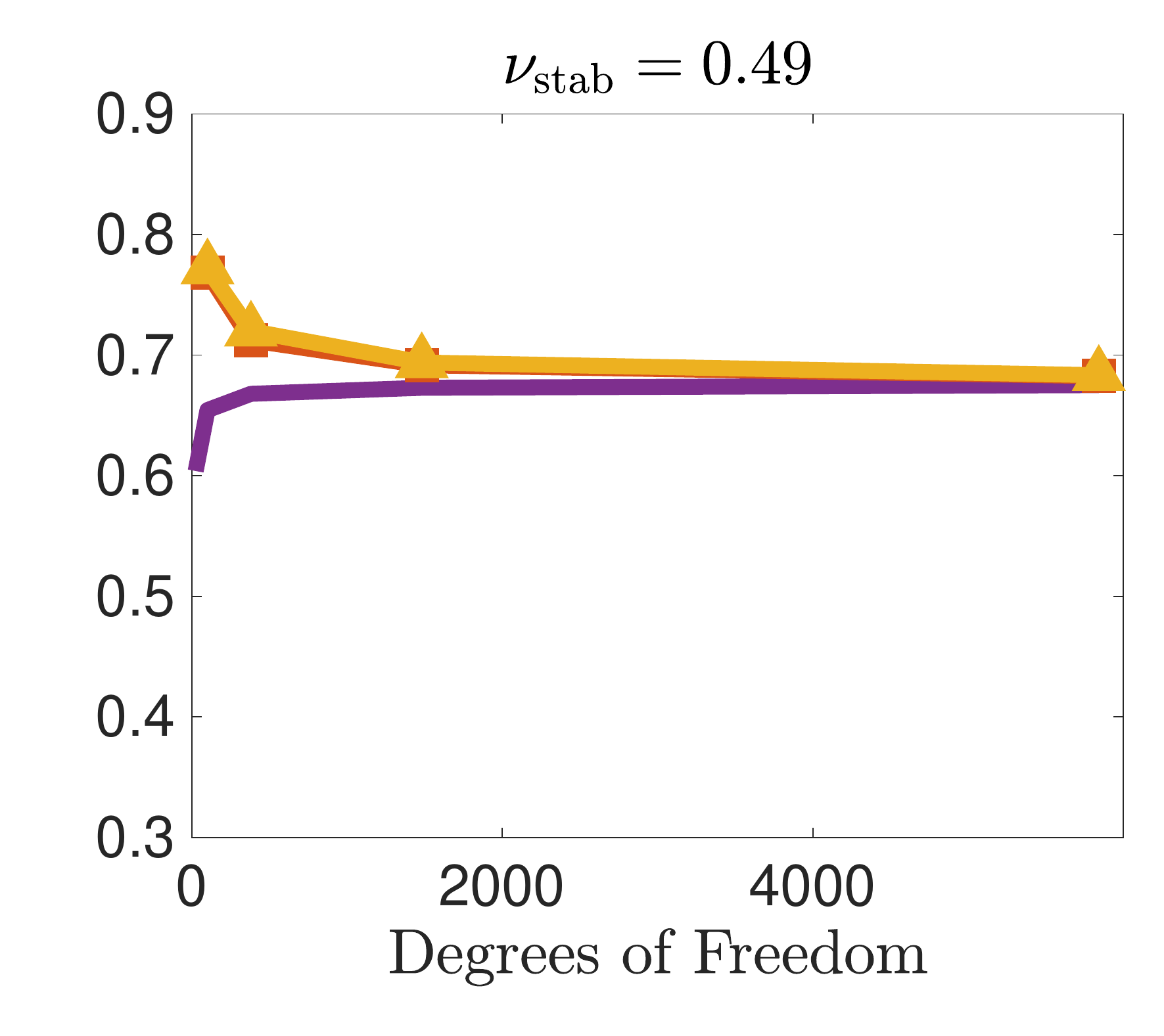}
        \end{subfigure}
        \begin{subfigure}{.3\textwidth}
                \includegraphics[width=\textwidth]{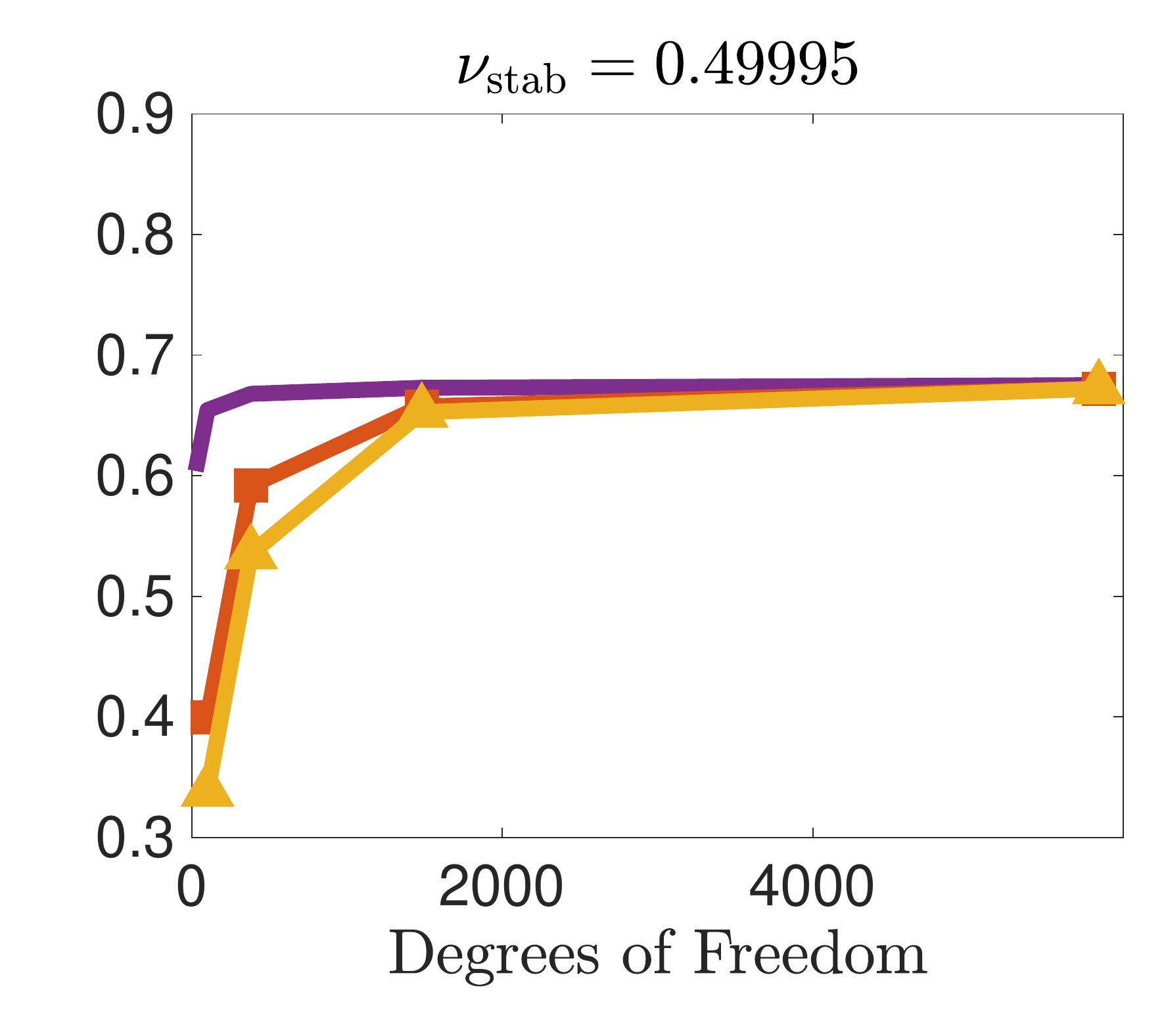}
        \end{subfigure} \\
    \end{tabular}
    \caption{Vertical displacements of the top corner point of the Cook's membrane benchmark, highlighted in Fig.~\ref{f:Cooks_schematics}, for different choices of peridynamic horizon size $\horizonsize$ and numerical Poisson's ratio $\nu_{\mathrm{stab}}$. The solid DoF range from $101$ to $5841$. Note that locking clearly occurs for $\nu_{\mathrm{stab}} = 0.49995$. As in standard computational mechanics approaches, however, the IPD formulation ultimately converges under grid refinement even with high (but fixed) levels of volumetric penalization.}
    \label{f:Cooks_disp}
\end{figure}

Fig.~\ref{f:Cooks_deformation} shows the plane sheet after the deformation along with pointwise values of the non-local Jacobian determinant of non-local deformation gradient tensor at each material point. 
Fig.~\ref{f:Cooks_disp} shows the $y$-displacement of the top-right corner in Fig.~\ref{f:Cooks_schematics} at the steady states for various numerical Poisson's ratios $\nu_{\mathrm{stab}}$ and PD horizon sizes $\horizonsize$ under grid refinement. 
Fig.~\ref{f:Cooks_disp} shows that the displacements obtained using the IPD method are comparable to those obtained using the classical FE method, and that they converge under grid refinement to approximately $0.67 \, \text{cm}$.
With $\nu_{\mathrm{stab}} = 0.49995$, a larger volumetric penalty causes volumetric locking, which results in smaller displacements when low mesh resolutions are used in the simulations.
However, under grid refinement, we ultimately recover accurate deformations for fixed finite values of numerical bulk modulus, as in classical methods for nearly incompressible elasticity.

\begin{figure}[]
\centering
   \begin{tabular}{cc}
        \begin{subfigure}{.3\textwidth}
          		\includegraphics[width=\textwidth]{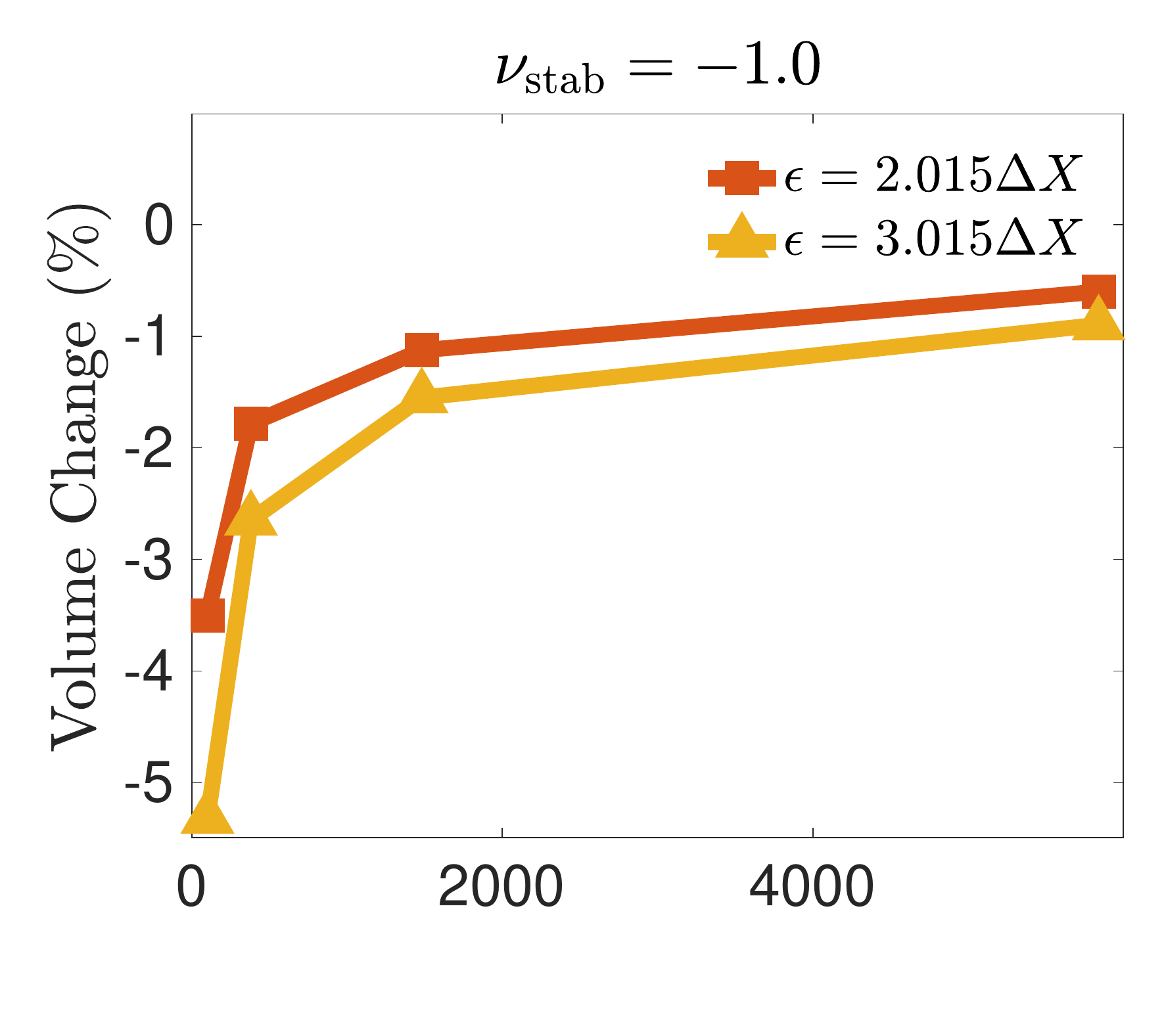}
        \end{subfigure} 
        \begin{subfigure}{.3\textwidth}
                \includegraphics[width=\textwidth]{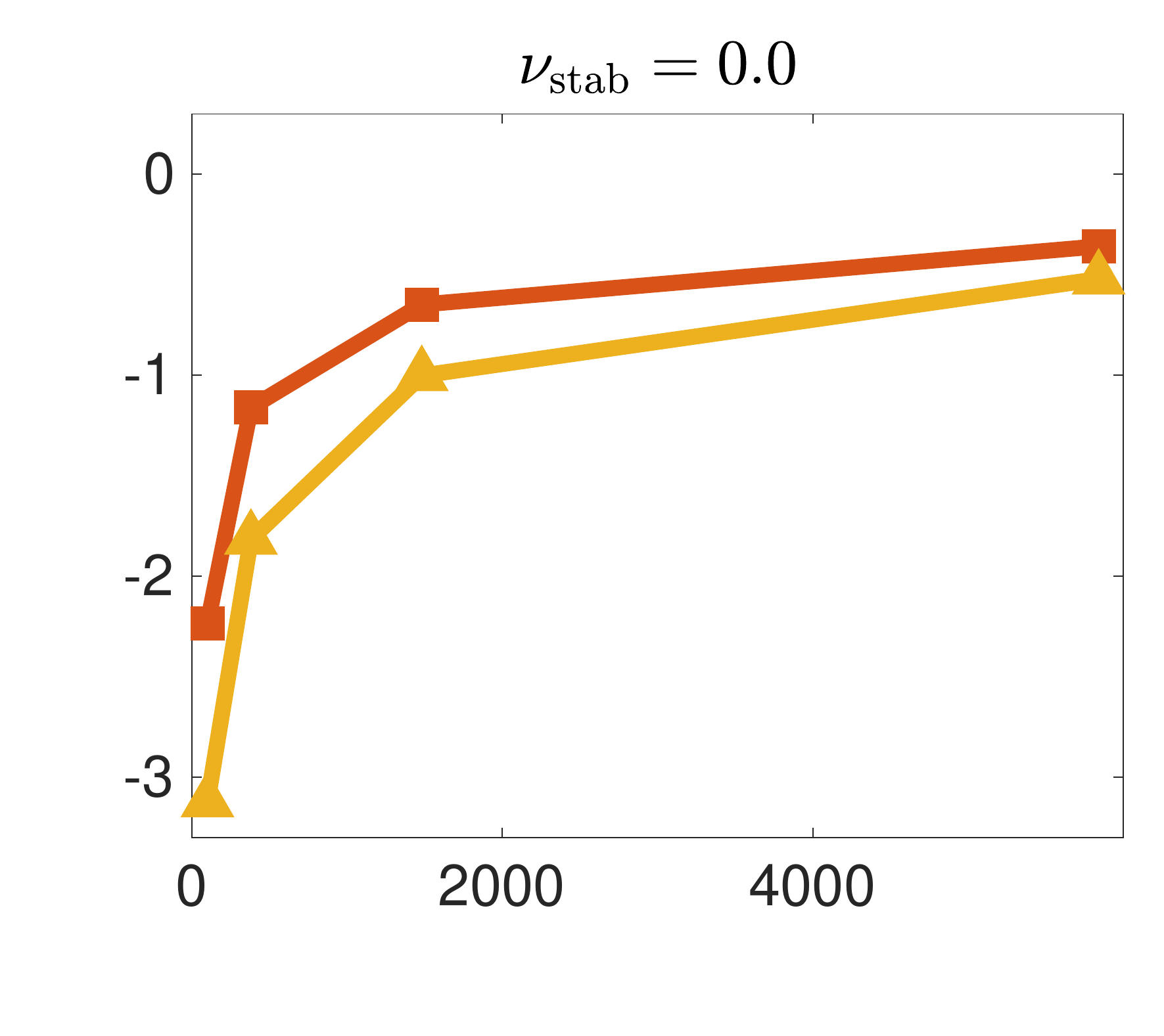}
        \end{subfigure}
        \begin{subfigure}{.3\textwidth}
                \includegraphics[width=\textwidth]{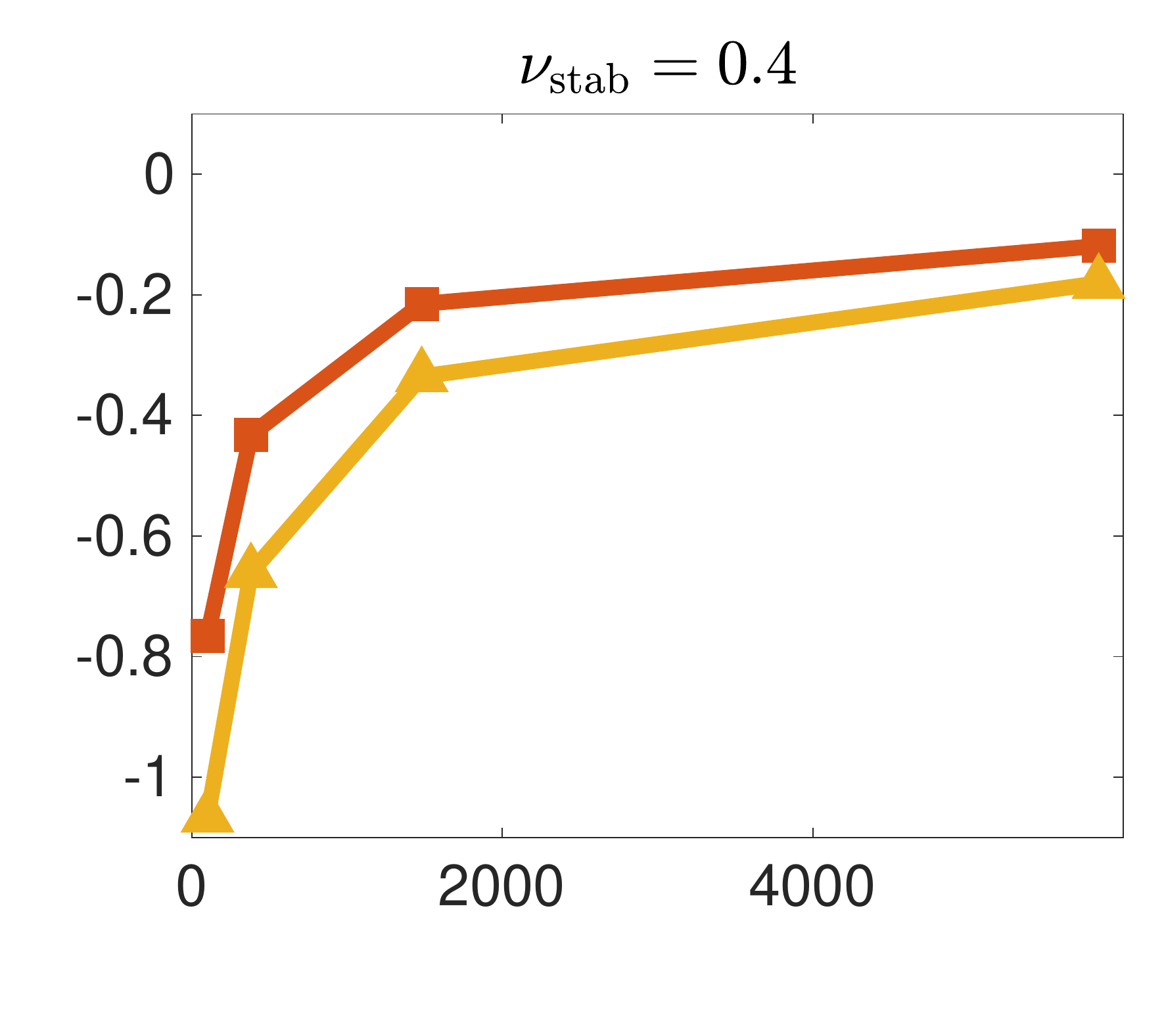}
        \end{subfigure} \\
         \begin{subfigure}{.3\textwidth}
          		\includegraphics[width=\textwidth]{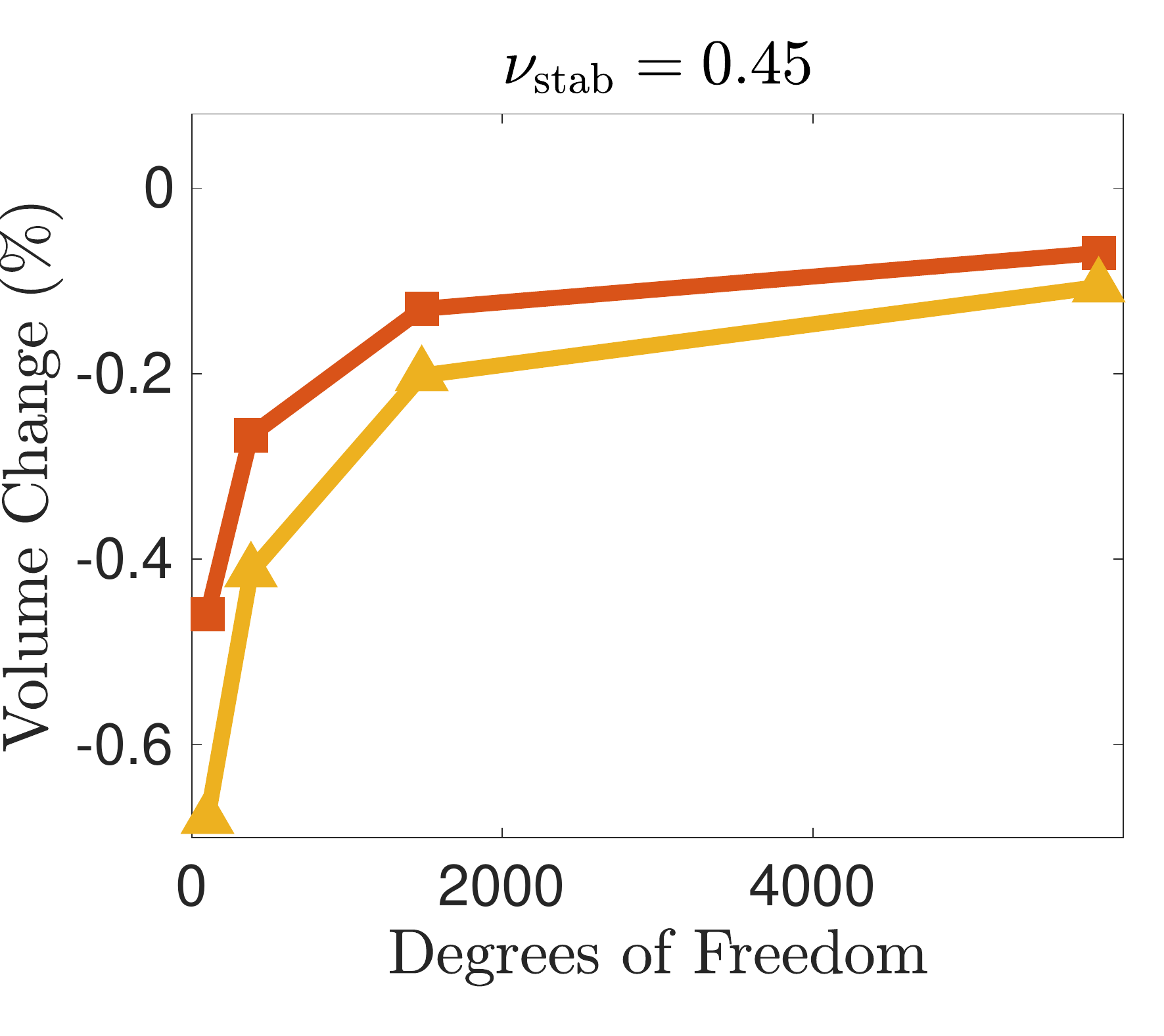}
        \end{subfigure} 
        \begin{subfigure}{.3\textwidth}
                \includegraphics[width=\textwidth]{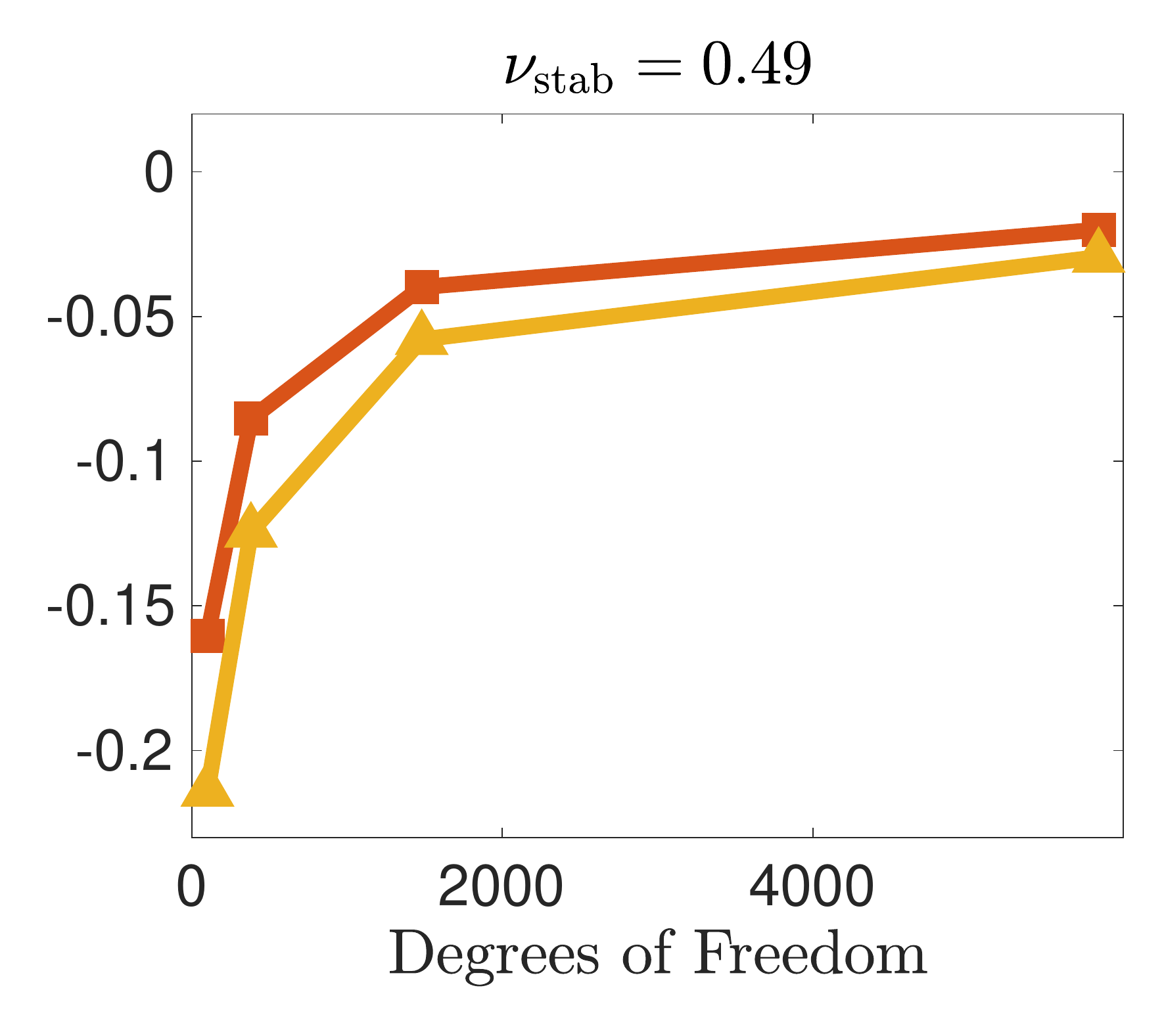}
        \end{subfigure}
        \begin{subfigure}{.3\textwidth}
                \includegraphics[width=\textwidth]{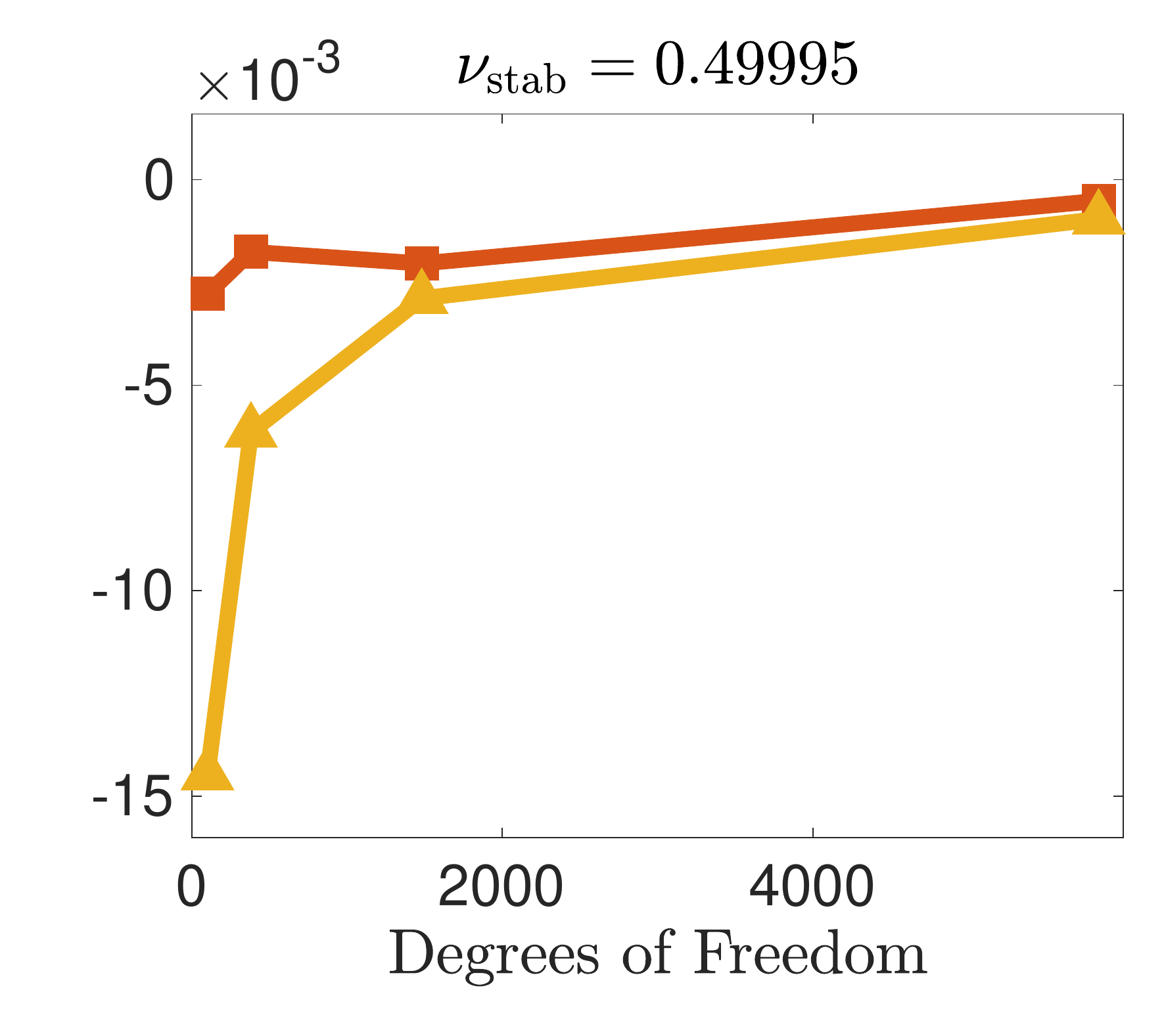}
        \end{subfigure} \\
    \end{tabular}
    \caption{Volume change of the Cook's membrane benchmark for different choice of horizon size $\horizonsize$ and numerical Poisson's ratio $\nu_{\mathrm{stab}}$.The solid DoF range from $101$ to $5841$. The largest change is approximately $5.3 \%$.}
    \label{f:Cooks_vol}
\end{figure}

Fig.~\ref{f:Cooks_vol} shows volume changes observed under deformation for different grid spacings. 
With smaller values of $\nu_{\mathrm{stab}}$, slight volume changes (between $0.6\%$ and $5.3\%$) are observed under loading. 
The volume change becomes negligible (up to $0.0005\%$) when larger values of $\nu_{\mathrm{stab}} \ge 0.4$ are used. 
The volume change using IPD is comparable to the results obtained using IFED, which is between $0.000021\%$ and $0.1\%$. 
We can expect negligible volume leaking or locking under grid refinement in all IPD simulations. 
In addition, relatively consistent results are obtained for all considered choices of the PD horizon size.

\subsubsection{Torsion}
\label{s:Benchmark_Torsion}
\begin{figure}[]
\centering
    \includegraphics[width=.65\textwidth]{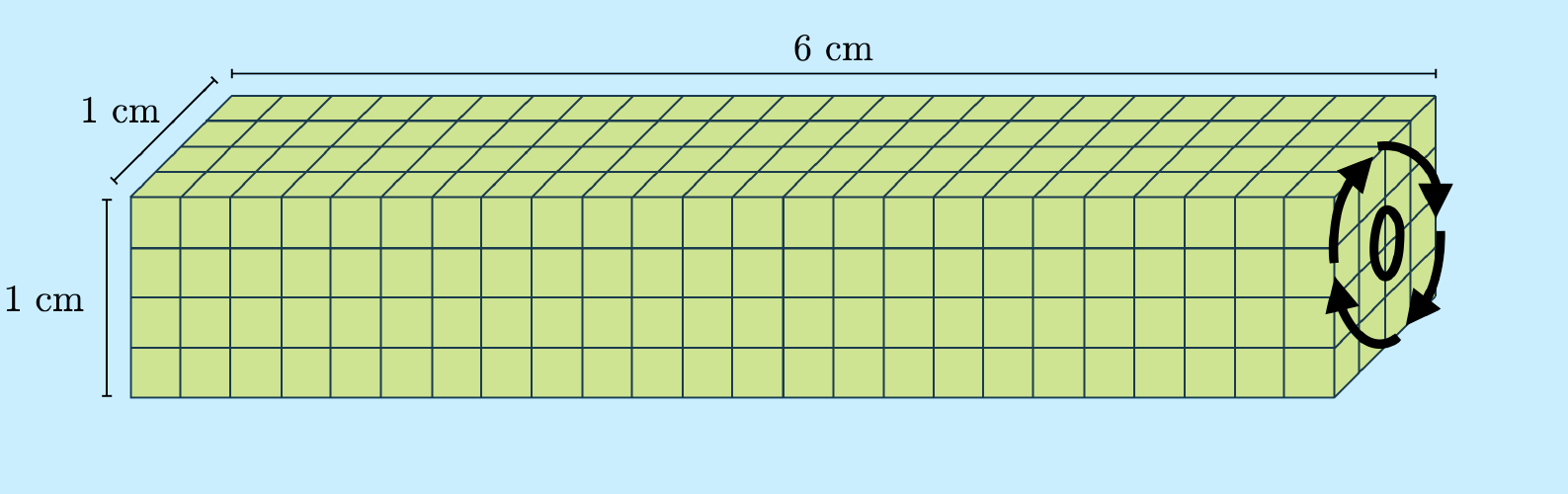}
    \caption{Schematic diagram for the three-dimensional torsion benchmark (Sec.~\ref{s:Benchmark_Torsion}). The computational domain is $\Omega = [0,L]^3$, with $L = 9  \, \text{cm}$, and the three-dimensional beam is placed at the center of the domain. Zero fluid velocity is enforced on the outer boundaries of the computational domain.}
    \label{f:torsion_schematics}
\end{figure}

We use a three-dimensional beam under torsion to investigate three dimensional hyperelastic material responses. 
This benchmark is based on a test suggested by Bonet et al.~\cite{bonet2015computational} and later modified by Vadala-Roth et al.~\cite{vadala2020stabilization} to use with the IFED method.
The computational domain is the cube $\Omega = [0,L]^3$, with $L = 9 \,  \text{cm}$. 
One side of the beam is fixed in place, and a torsion is applied to the opposite end through displacement boundary conditions. 
Fig.~\ref{f:torsion_schematics} provides a schematic of the test case.
The right surface is rotated by the linear function $\theta(t)$ from $0$ to $\theta_{T_{\text{f}}} = 2.5 \pi$ in time. 
The maximum angle of rotation is achieved at $T_{\text{l}} = 0.4 T_{\text{f}}$, with $T_{\text{f}} = 5 \,  \text{s}$. 
All other surfaces have zero traction boundary conditions. 
Material damage and failure are not considered.

We use the modified Mooney-Rivlin material model \cite{vadala2020stabilization} 
\begin{align}\label{mooney_rivlin}
\Psi &= c_1\left(J^{-2/3} I_1 - 3 \right) + c_2 \left( \frac{J^{-4/3}}{2} I_2 - 3 \right)  + \frac{\kappa_{\text{stab}}}{2} \left( \ln J \right)^2,\\
\PP &= 2 c_1 J^{-2/3} \left( \FF - \frac{I_1}{3} \FF^{-T} \right) + 2c_2 J^{-4/3}\left(I_1 \FF - \FF \mathbb{C} - \frac{I_2}{3} \FF^{-T}  \right) + \kappa_\mathrm{stab} \ln \left(J\right) \FF^{-T},
\end{align}
in which $c_1$ and $c_2$ are material parameters, $I_1 = \text{tr}\left(\mathbb{C}\right)$, and $I_2 =\text{tr}\left(\mathbb{C}\right)^2 - \text{tr}\left(\mathbb{C}^2\right)$.  
The material parameters are set to $c_1 = c_2 = 9000 \, \frac{\text{dyn}}{\text{cm}^2}$, and $G = \left(c_1 + c_2 \right)$ is used to determine the numerical bulk modulus. 
The density and viscosity of the surrounding fluid are set to $\rho = 1.0 \,  \frac{\text{g}}{\text{cm}^3}$ and  $\mu = 0.04 \, \frac{\text{dyn$\cdot$s}}{\text{cm}^2}$, respectively. 
This larger value of viscosity compared to the previous benchmarks is used to accelerate reaching the steady state.

Overall deformations and volume conservation using the numerical Poisson's ratio of $\nu_{\text{stab}} = 0.4$ in Sec.~\ref{s:Benchmark_Compression} and Sec.~\ref{s:Benchmark_Cooks} are consistent with the classical elasticity results under grid refinement. 
In general, a larger value of numerical bulk modulus requires a smaller time step size. 
For the remainder of the IPD simulations presented herein, we only consider two fixed numerical Poisson's ratios; $\nu_{\text{stab}} = 0.4$ for nearly incompressible hyperelastic material models and $\nu_{\text{stab}} = -1.0$ for comparison tests.

 \begin{figure}[]
\centering
\hspace{-.1in}
	\begin{tabular}{cc}
	\includegraphics[width=.95\textwidth]{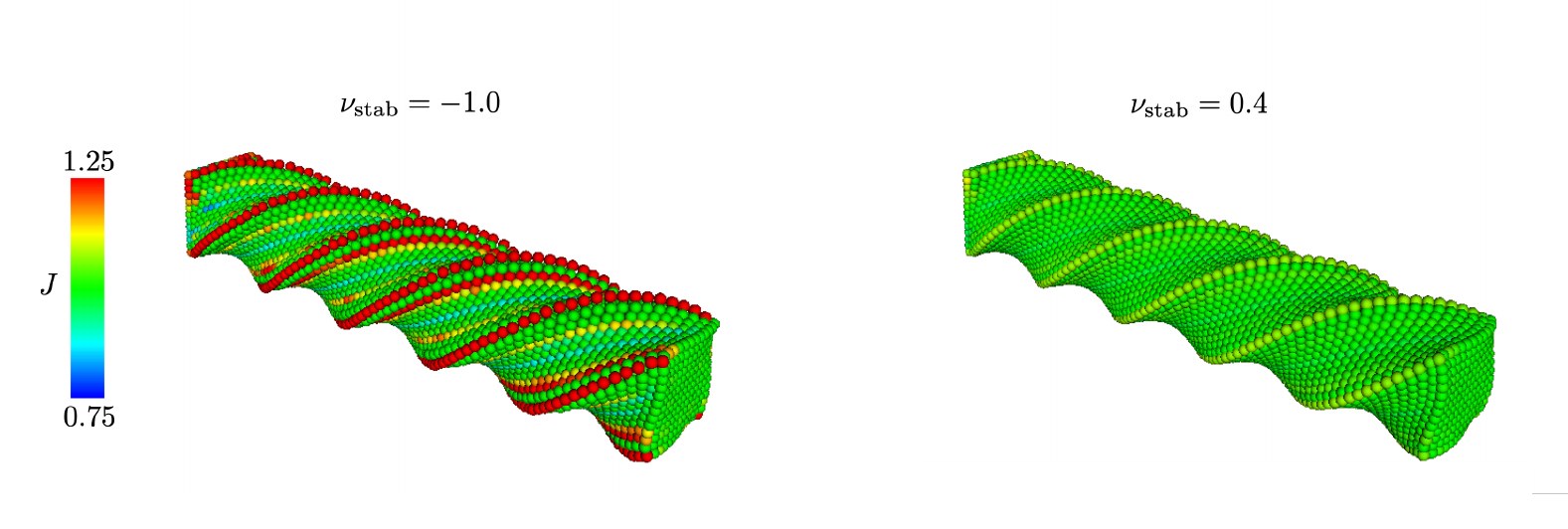}
	\end{tabular}
    \caption{Deformation of the three-dimensional beam with the values of $J$ at material points using the Mooney-Rivlin material model with $c_1 = c_2 = 9000 \, \frac{\text{dyn}}{\text{cm}^2}$. The deformations are computed using 12337 solid DoF and $\horizonsize = 2.015 \Delta X$. The left panel shows the deformation obtained using $\nu_{\text{stab}}  = -1.0$, and the right panel shows the result for$\nu_{\text{stab}}  = 0.4$.}
    \label{f:torsion_deformation}
\end{figure}

\begin{figure}[]
\centering
   \begin{tabular}{cc}
        \begin{subfigure}{.4\textwidth}
          		\includegraphics[width=\textwidth]{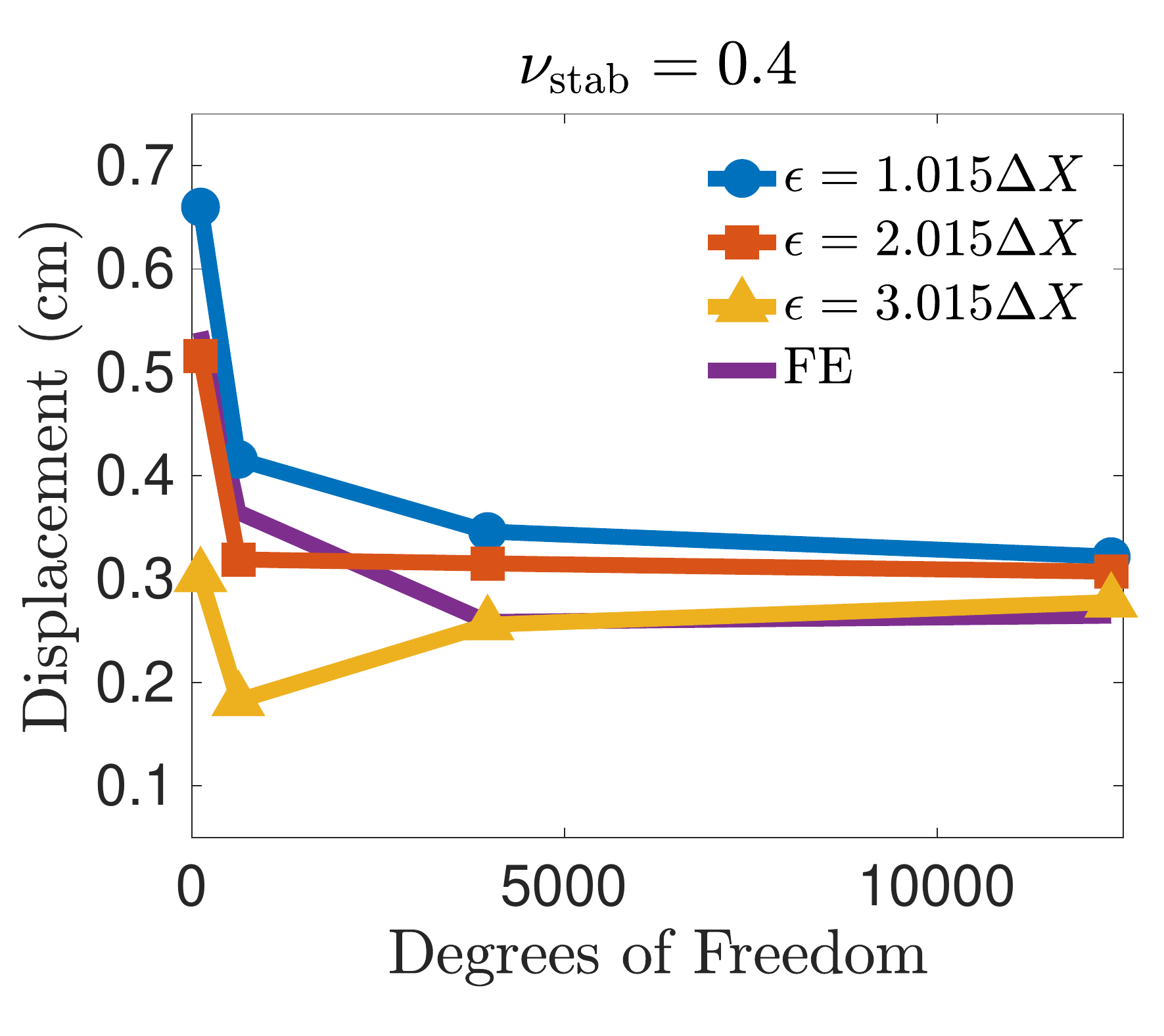}
          		\caption{}
          		 \label{f:Torsion_disp}
        \end{subfigure} 
         \begin{subfigure}{.4\textwidth} \hspace{.05\textwidth}
         		\includegraphics[width=\textwidth]{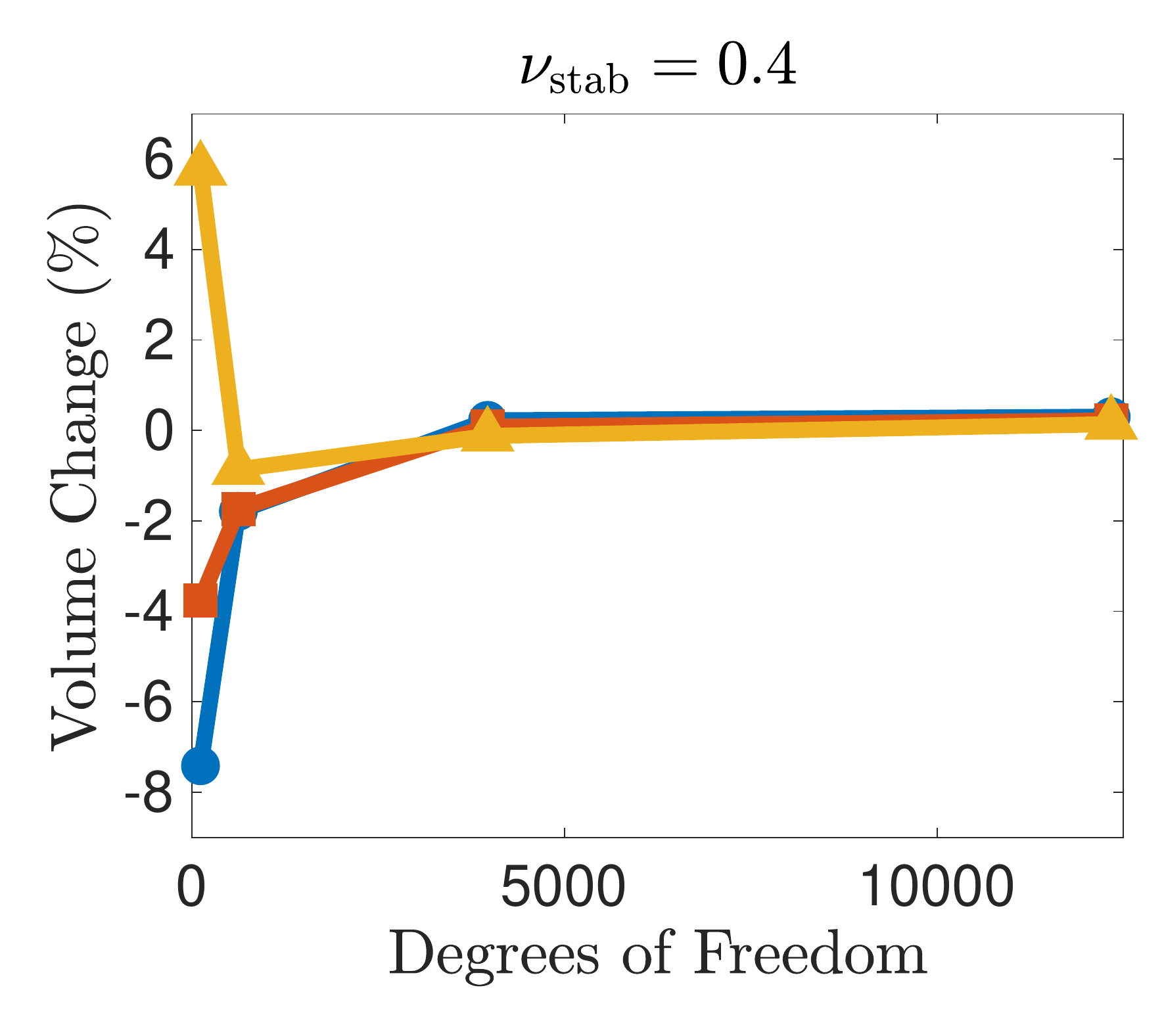}
          		\caption{}
          		 \label{f:Torsion_vol}
        \end{subfigure} 
    \end{tabular}
    \caption{(a) Displacements of the point of interest, highlighted in Fig.~\ref{f:torsion_schematics}, for different choices of horizon size $\horizonsize$. The solid DoF range from $117$ to $12337$. (b) Volume change of the beam for different choices of horizon size $\horizonsize$. }
\end{figure}

Fig.~\ref{f:torsion_deformation} illustrates the deformations of the beam under torsion along with the Jacobian determinant of the non-local deformation gradient tensor at each material point. 
Fig.~\ref{f:Torsion_disp} shows the maximum displacement of the center of the top surface in Fig.~\ref{f:torsion_schematics} at the steady states for various sizes of $\horizonsize$ under grid refinement. 
The displacements using IPD are comparable to the classical FE results and converge under grid refinement to a value of approximately $0.27 \, \text{cm}$. 
Fig.~\ref{f:Torsion_vol} shows the volume change of the beam for different numbers of the solid DoF, ranging from $0.13 \%$ to $7.41\%$. 
Volume changes obtained using the IFED method are between $0.16 \%$ and $11 \%$, which are comparable to the IPD simulations.

\subsubsection{Elastic band}
\label{s:Benchmark_ELS}
\begin{figure}[]
\centering
    \includegraphics[width=.8\textwidth]{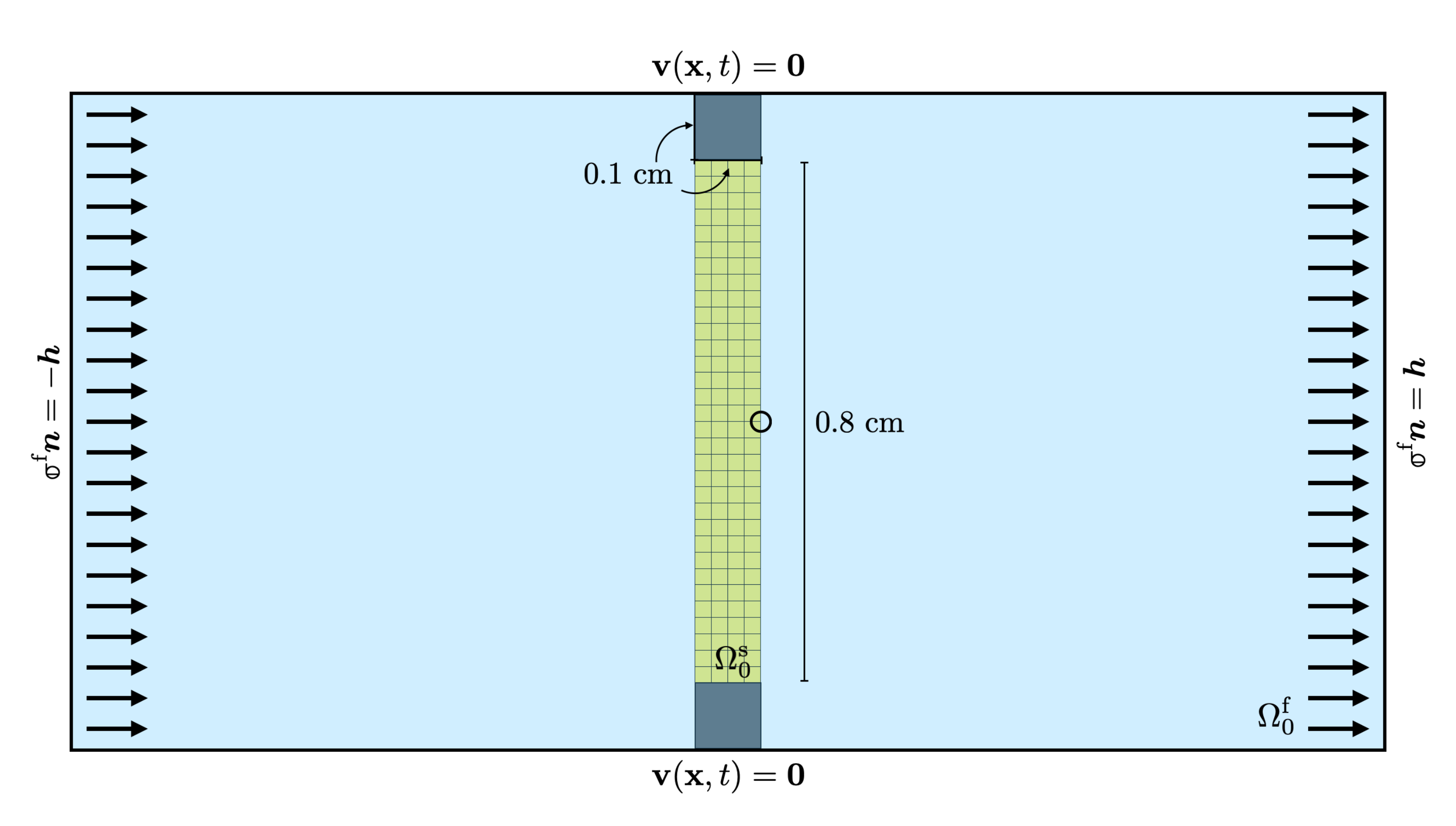}
    \caption{Schematic diagram for the elastic band benchmark (Sec.~\ref{s:Benchmark_ELS}).  The initial configurations of the immersed structure and a fluid are denoted by $\Omega_0^{\text{s}}$ and $\Omega_0^{\text{f}}$, respectively. The entire computational domain is $\Omega = \Omega_0^{\text{s}} \cup \Omega_0^{\text{f}}$. Zero fluid velocity is enforced on the top and bottom boundaries of the computational domain, and fluid traction boundary conditions are applied to the left and right boundaries. Fluid traction is set to $\bm{h}(t) = \left(10 \sin\left(\frac{\pi t}{2 T_{\text{l}}}\right) , 0\right) \, \frac{\mathrm{dyn}}{\mathrm{cm}^2}$ when $t< T_{\text{l}}$ and  $\bm{h}(t) = \left(10 , 0\right)\, \frac{\mathrm{dyn}}{\mathrm{cm}^2}$ otherwise.}
    \label{f:ELS_BND_schematics}
\end{figure}

This benchmark examines  deformations of an elastic band that are driven by fluid forces themselves. 
The deformations of the elastic band are simulated under fluid pressure loading. 
The computational domain is $\Omega = [0,2L] \times [0,L]$, with $L = 1 \, \mathrm{cm}$. 
Fluid traction boundary conditions are imposed on the boundaries of the computational domain as $\vec{\bbsigma}^{\text{f}} (\x,t) \bm{n} (\x) = \bm{h}(t)$ and $\vec{\bbsigma}^{\text{f}} (\x,t) \bm{n} (\x) = - \bm{h}(t)$ on the left and right, respectively, in which $\vec{\bbsigma}^{\text{f}} $ is the fluid stress tensor and $\bm{h}(t) = \left(10 \sin\left(\frac{\pi t}{2 T_{\text{l}}}\right) , 0\right) \, \frac{\mathrm{dyn}}{\mathrm{cm}^2}$ when $t< T_{\text{l}}$ and $\bm{h}(t) = \left(10 , 0\right) \, \frac{\mathrm{dyn}}{\mathrm{cm}^2}$ otherwise.
The load time is $T_{\text{l}} = 5 \, \mathrm{s}$. 
Zero fluid velocity conditions are applied to the top and bottom boundaries of the computational domain. 
Different from the previous benchmarks, the difference in pressure across the computational domain causes the deformations of the elastic band. 
Both top and bottom surfaces of the elastic band are attached to fixed blocks.
The stationary blocks serve to block the fluid flow between the wall and the flexible band. 
Fig.~\ref{f:ELS_BND_schematics} provides a schematic of this test case. 
We measure the maximum displacement of the point of interest, the encircled point in Fig.~\ref{f:ELS_BND_schematics}, at $T_{\text{f}} = 15\,  \mathrm{s}$. 
A shear modulus of $G = 200\, \frac{\mathrm{dyn}}{\mathrm{cm}^2}$ is used for the nearly incompressible neo-Hookean material model. 
Damage and failure of the elastic band are not allowed. 
An additional damping force in the structure is used with $\eta = 10 \,  \frac{\text{g}}{\text{s}}$, to accelerate reaching steady state.

\begin{figure}[]
\centering
   \begin{tabular}{ccc}
    \includegraphics[width=.85\textwidth]{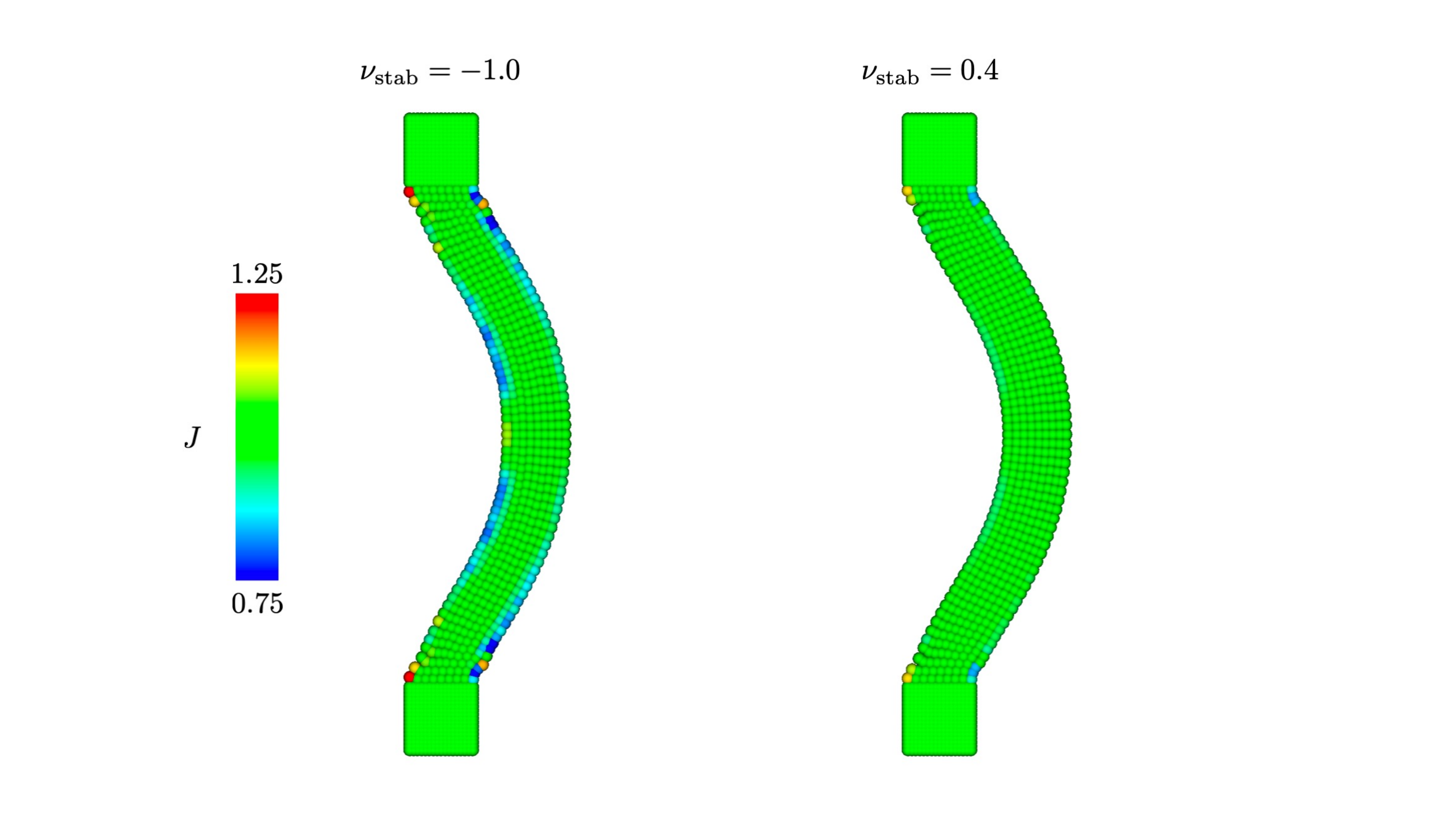}
	\end{tabular}
    \caption{Deformations of an elastic band with the values of $J$ at material points using the neo-Hookean material model with $G = 200 \, \frac{\mathrm{dyn}}{\mathrm{cm}^2}$. The deformations are computed using 1261 solid degrees of freedom and $\horizonsize = 2.015 \Delta X$. Note that the number of solid DoF only considers the number of Lagrangian points of the band. The left panel shows the deformation obtained using $\nu_{\text{stab}} = -1.0$, and the right panel shows the result for $\nu_{\text{stab}} = 0.4$.}
    \label{f:ELS_BND_deformation}
\end{figure}

\begin{figure}[]
\centering
   \begin{tabular}{cc}
        \begin{subfigure}{.4\textwidth}
          		\includegraphics[width=\textwidth]{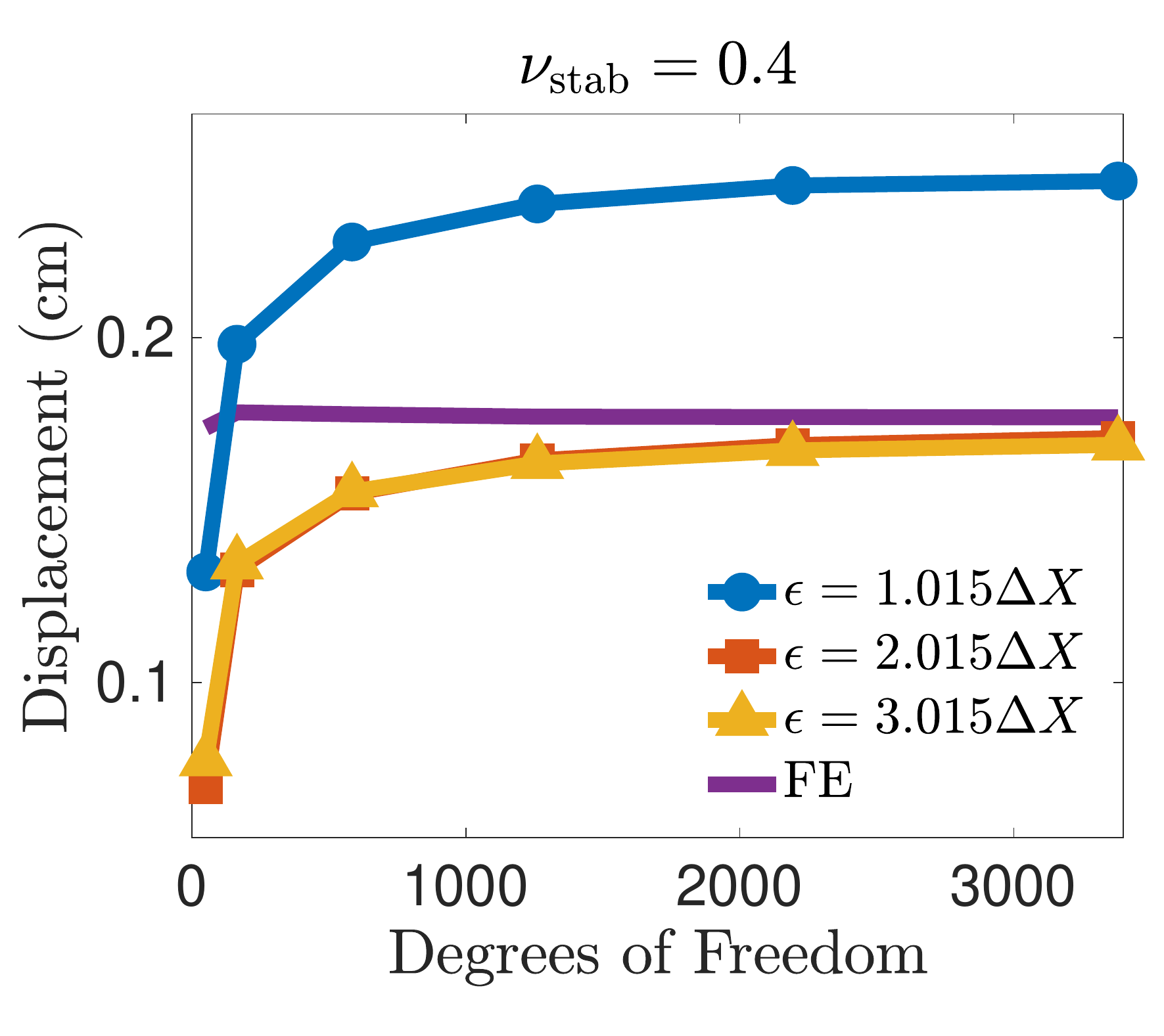}
          		\caption{}
          		 \label{f:ELS_BND_disp}
        \end{subfigure} \hspace{.05\textwidth}
         \begin{subfigure}{.4\textwidth}
        			\includegraphics[width=\textwidth]{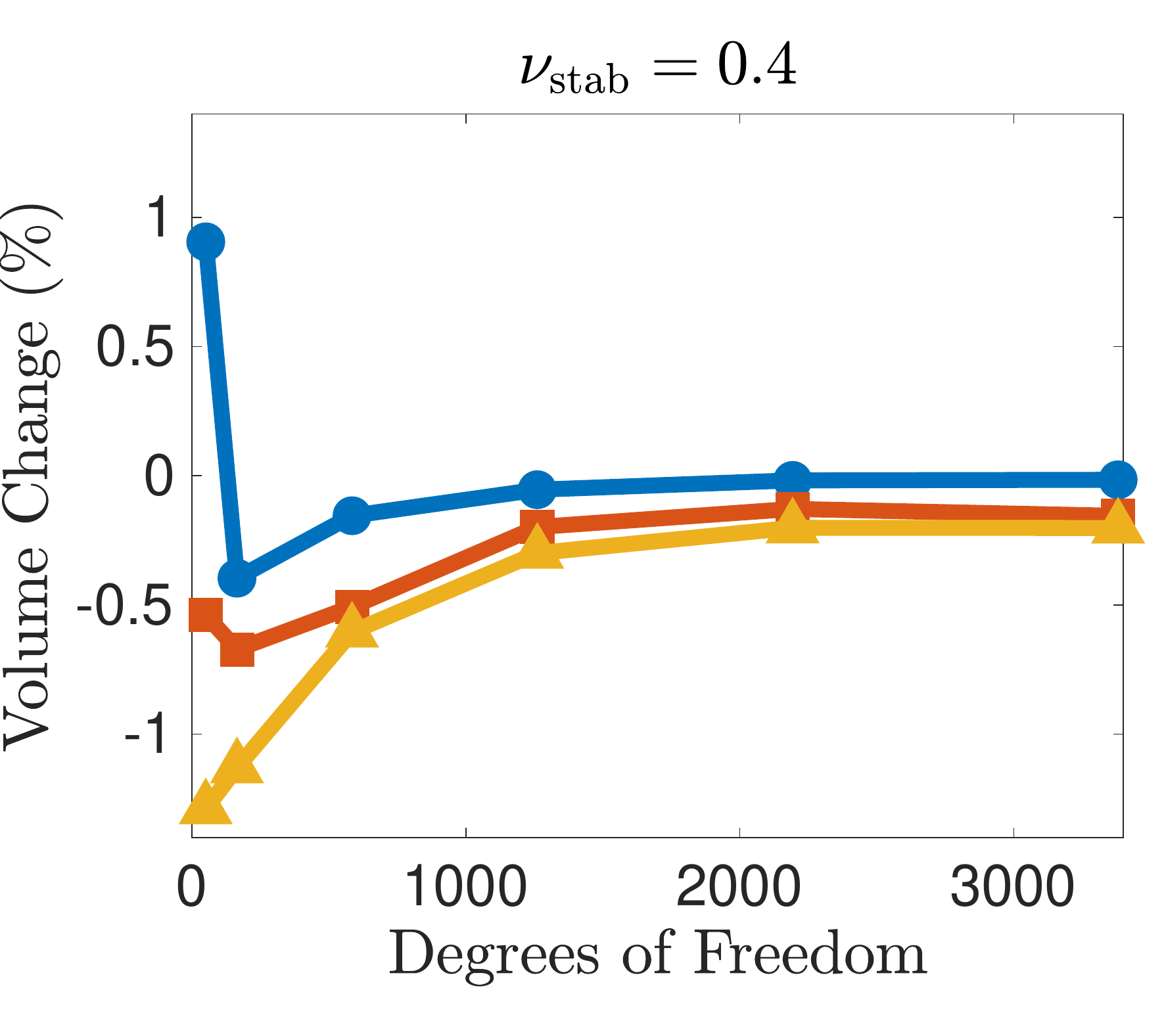}
          		\caption{}
          		 \label{f:ELS_BND_vol}
        \end{subfigure} 
    \end{tabular}
    \caption{(a) Horizontal displacements of the point of interest, highlighted in Fig.~\ref{f:ELS_BND_schematics}, for different choices of horizon size $\horizonsize$. The solid DoF range from $51$ to $3381$. (b) Volume change of the band for different choices of horizon size $\horizonsize$. }
\end{figure}

Fig.~\ref{f:ELS_BND_deformation} shows the deformations of the elastic band under pressure loading at steady state. 
Fig.~\ref{f:ELS_BND_disp} shows the maximum horizontal displacement of the point of interest in Fig.~\ref{f:ELS_BND_schematics} at the steady states for different choices of the PD horizon sizes under grid refinement, which is approximately $0.17\, \text{cm}$. 
Except for $\horizonsize = 1.015 \Delta X$, the results obtained using the IPD method are comparable to FE results. 
The absence of a diagonal connectivity in the immersed structure with $\horizonsize = 1.015 \Delta X$ causes a lack of resistance to this type of bending. 
Fig.~\ref{f:ELS_BND_vol} shows the volume change of the band for different numbers of the solid DoF.
We observe slight volume changes with smaller values of solid DoF, however, these are resolved under grid refinement. 
The volume conservation achieved by the IPD method (between $0.015 \%$ and $1.2\%$)  is comparable to the results generated by the IFED method (between $0.0015 \%$ and $2.1 \%$).

To investigate nontrivial fluid dynamics in this benchmark, we test the transient behavior of a dynamic version of the elastic band. 
Instead of gradually applying  the fluid traction as the static problem, fluid traction conditions on the boundaries of the computational domain are set to $\vec{\bbsigma}^{\text{f}} (\x,t) \bm{n} (\x) = \bm{h}(t)$, in which $\bm{h}(t) = \left(-10 , 0\right) \, \frac{\mathrm{dyn}}{\mathrm{cm}^2}$ and $\bm{h}(t) = \left(10,0\right) \, \frac{\mathrm{dyn}}{\mathrm{cm}^2}$ on the left and right, respectively. 
The final simulation time is set to $T_{\text{f}} = 10 \ \mathrm{s}$. 
Otherwise, we fix the rest of the test parameters as in the static problem, with no damping. 

\begin{figure}[]
\centering
\includegraphics[width=\textwidth]{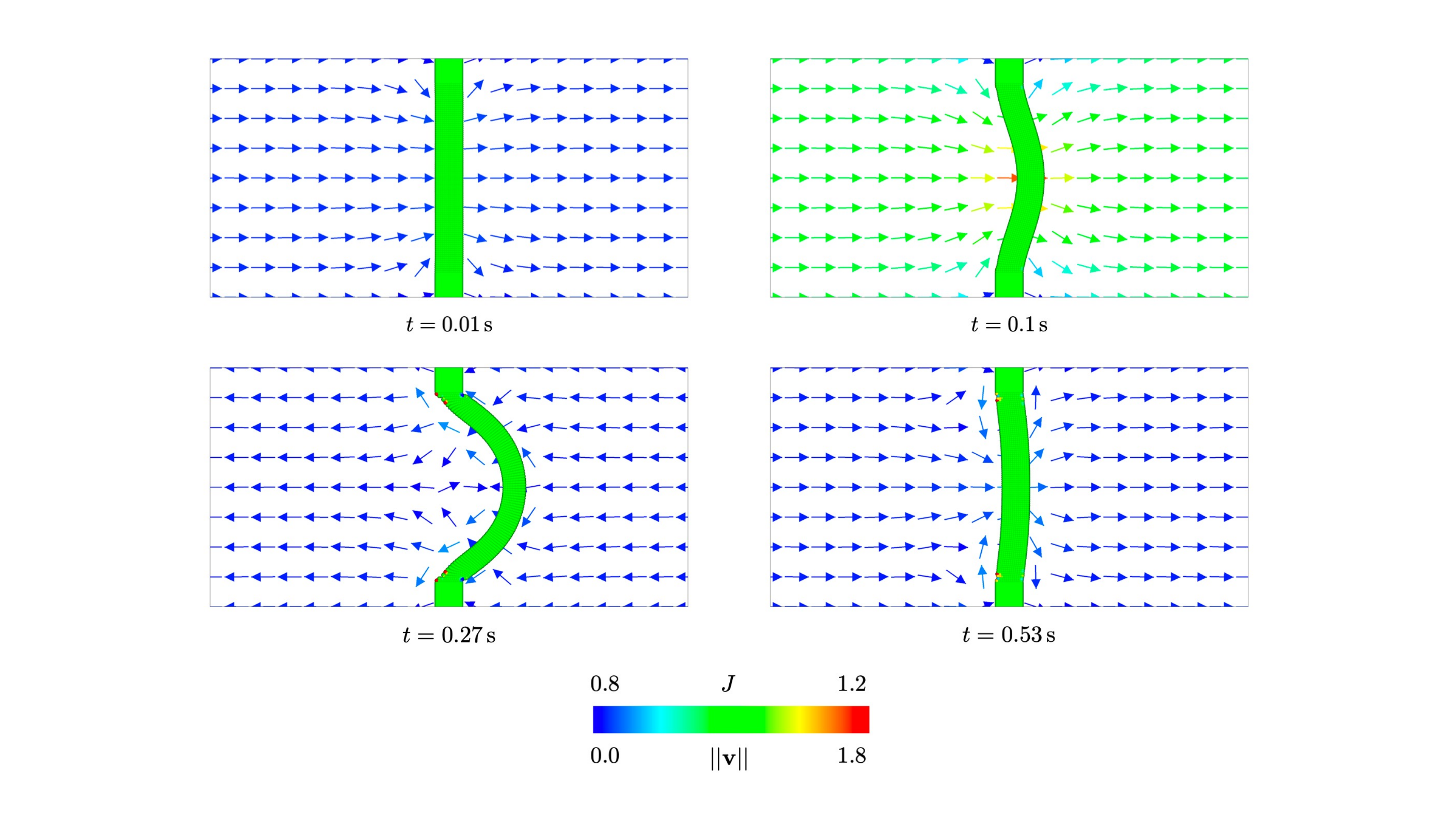}
\vspace{-.25in}
     \caption{Deformations of a dynamic version of the elastic band under the fluid traction force. The color represents the magnitude of the Eulerian velocity at each spatial point and the values of $J$ at material points. The deformations are computed using $1261$ solid DoF, $\horizonsize = 2.015 \Delta X$, and $\nu_{\text{stab}} = 0.4$. The band undergoes its largest deformation at $t = 0.27 \, \text{s}$ and enters another period of oscillation at $t = 0.53 \, \text{s}$.}
    \label{f:ELS_BND_dynamics}    
\end{figure}

\begin{figure}[]
\centering
\hspace{-.3in}
   \begin{tabular}{cc}
        \begin{subfigure}{.3\textwidth}
          		\includegraphics[width=\textwidth]{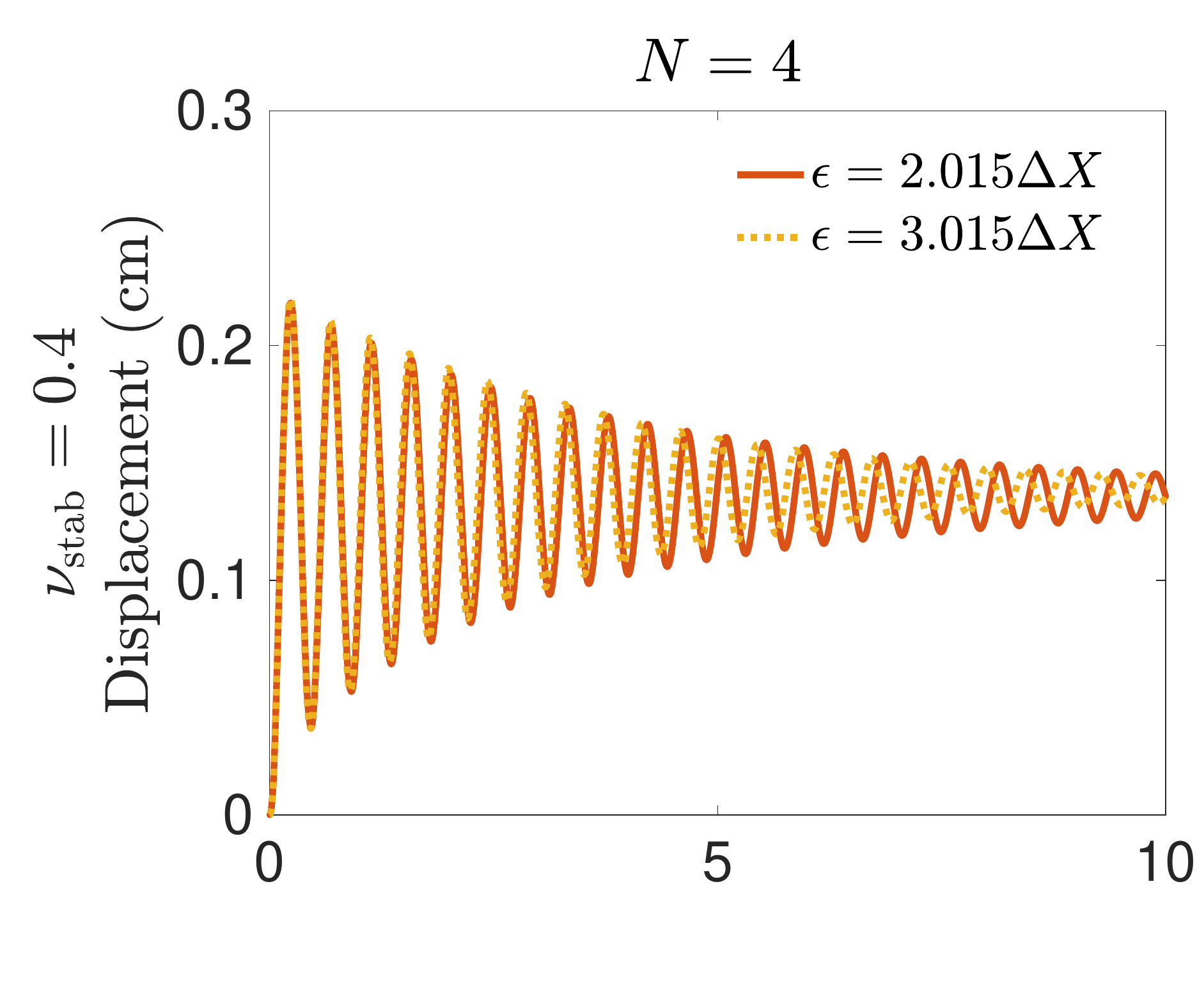}
        \end{subfigure} 
        \begin{subfigure}{.3\textwidth}
                \includegraphics[width=\textwidth]{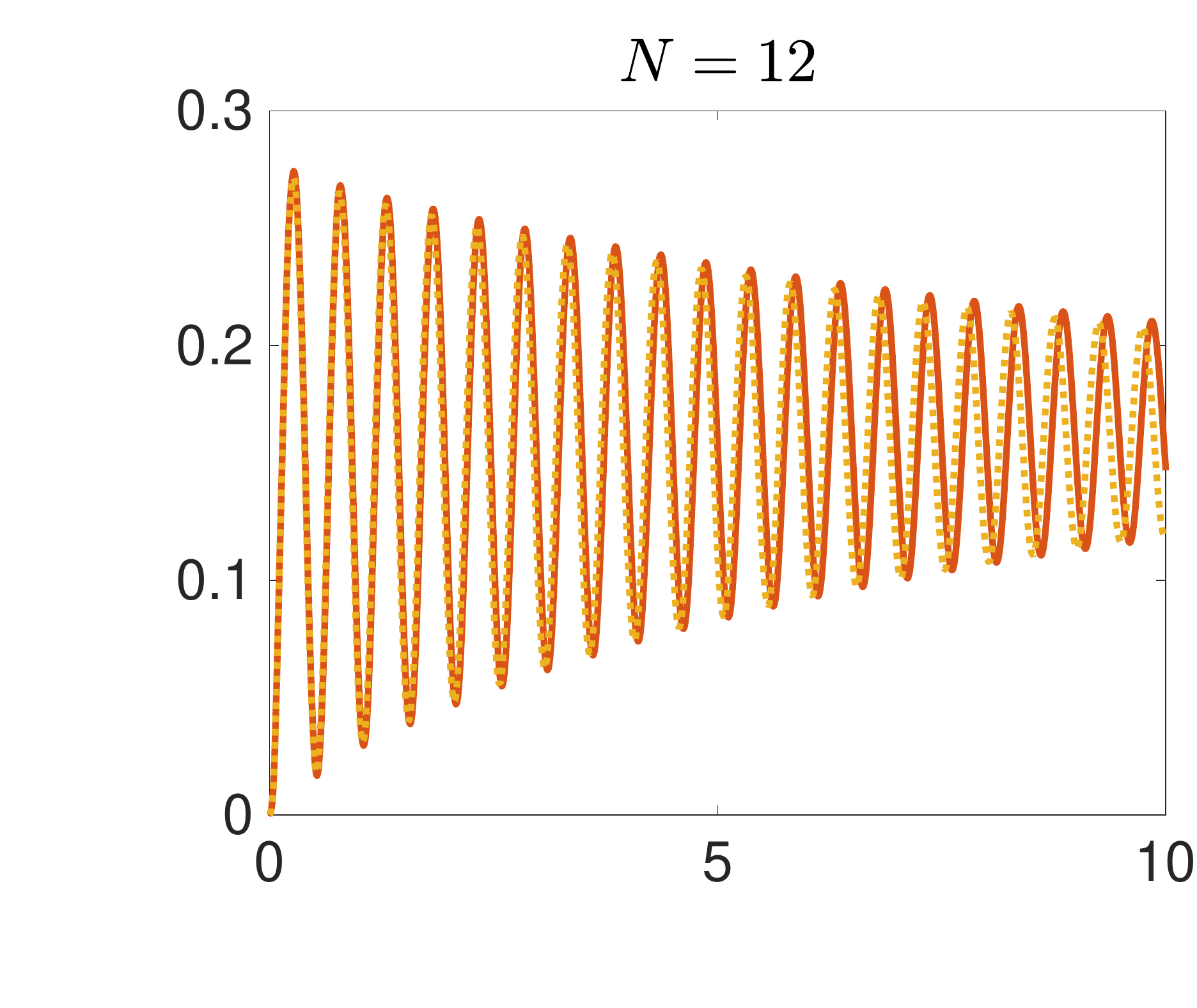}
        \end{subfigure}
        \begin{subfigure}{.3\textwidth}
                \includegraphics[width=\textwidth]{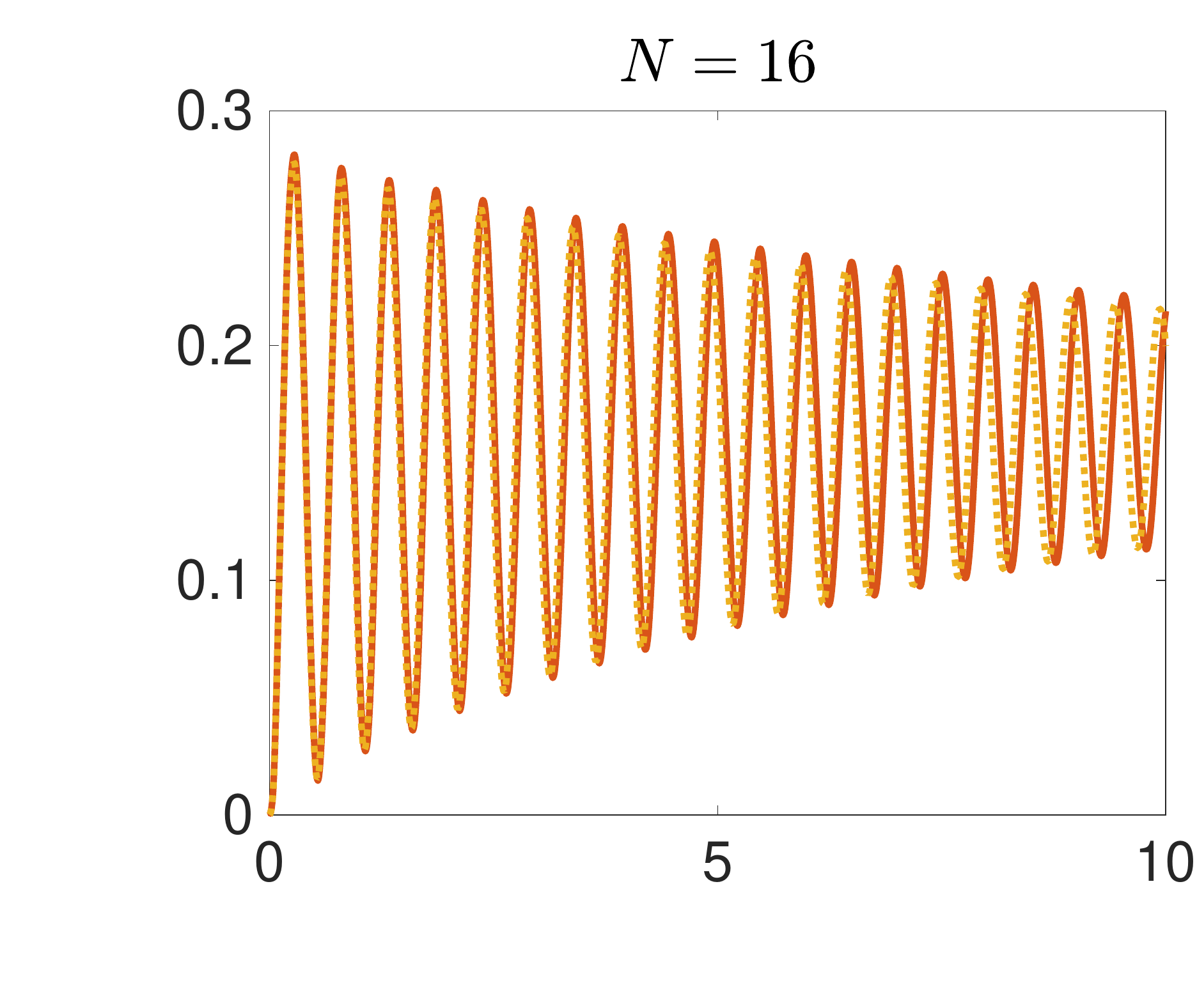}
        \end{subfigure} \\
         \begin{subfigure}{.3\textwidth}
          		\includegraphics[width=\textwidth]{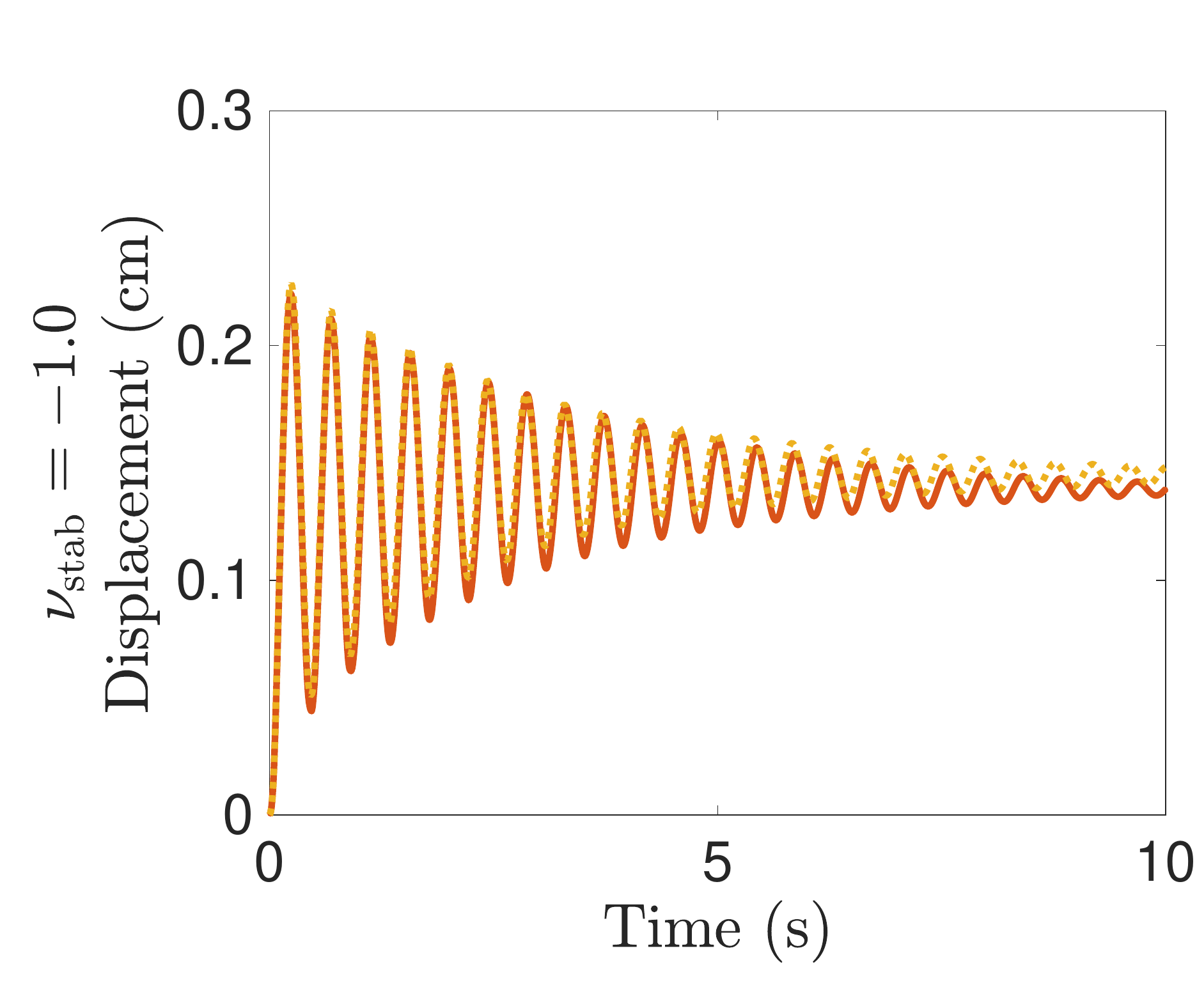}
        \end{subfigure} 
        \begin{subfigure}{.3\textwidth}
                \includegraphics[width=\textwidth]{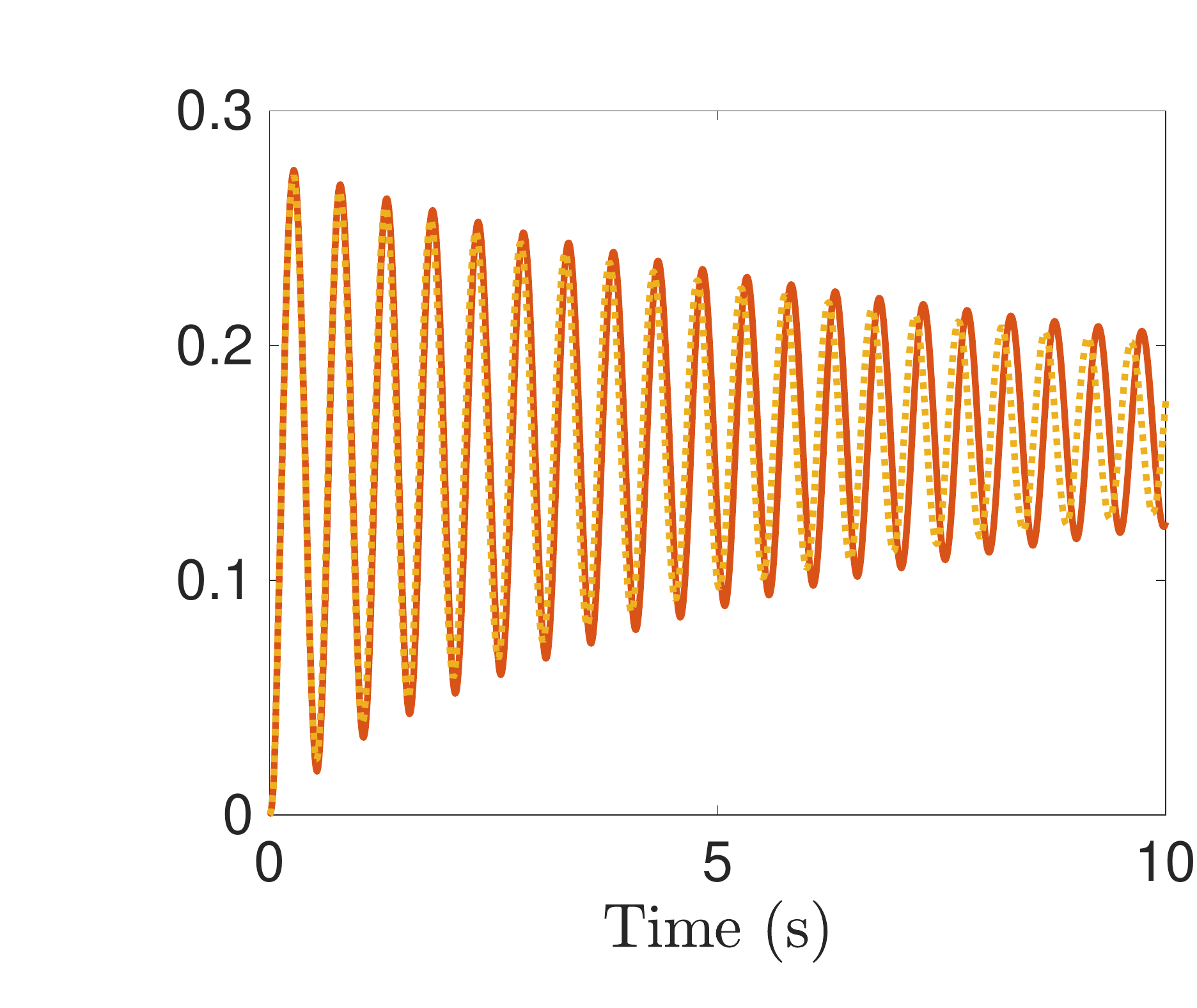}
        \end{subfigure}
        \begin{subfigure}{.3\textwidth}
                \includegraphics[width=\textwidth]{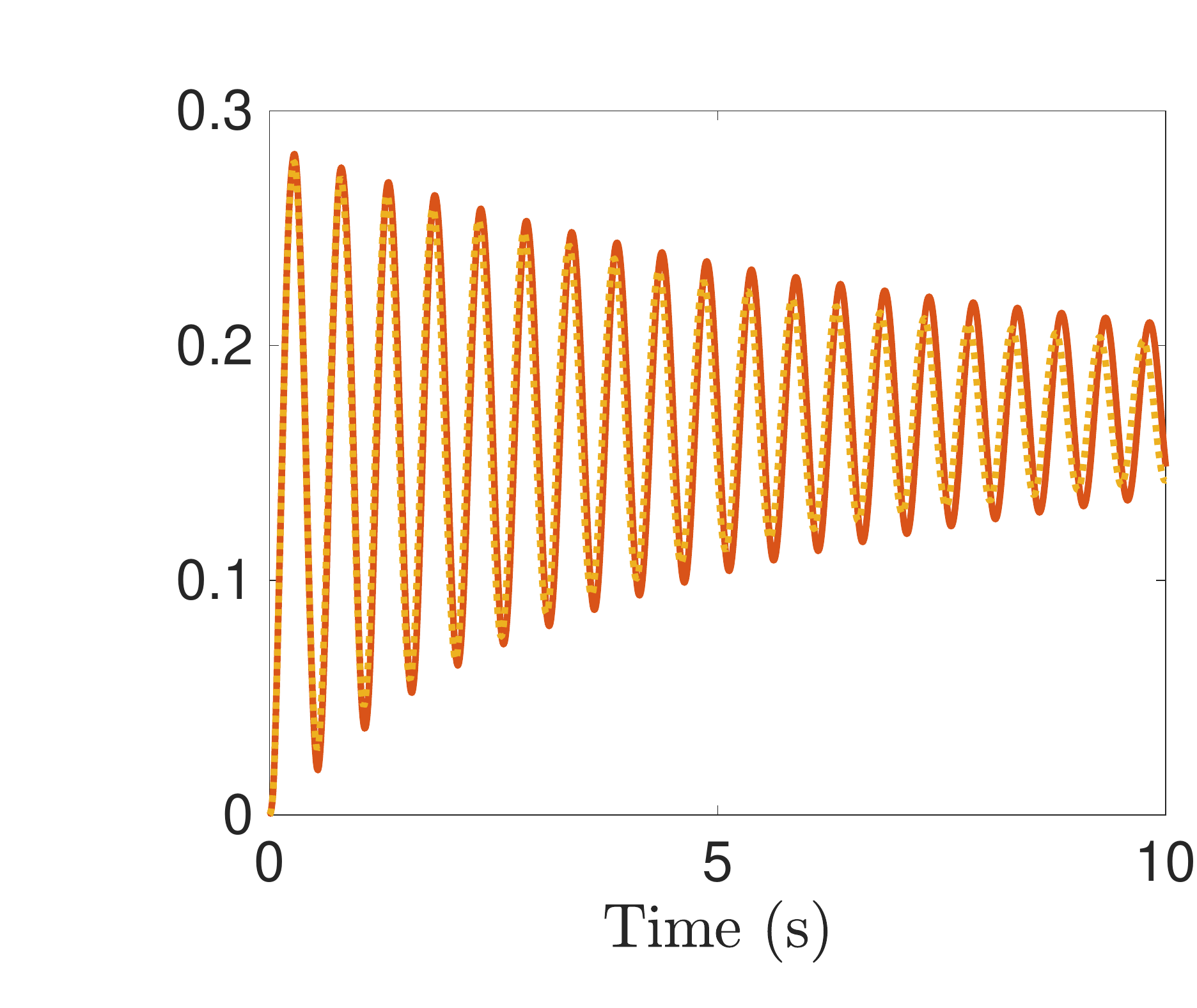}
        \end{subfigure} \\
    \end{tabular}
    \caption{Horizontal displacements of the point of interest, highlighted in Fig.~\ref{f:ELS_BND_schematics}, for different choices of horizon size $\horizonsize$ and numerical Poisson's ratio $\nu_{\text{stab}}$ under grid refinement.  $N=4$ corresponds to $165$ solid DoF, $N=12$ corresponds to $1261$ solid DoF, and $N = 16$ corresponds to $2193$ solid DoF.}
\label{f:ELS_BND_dynamics_displacement}	
\end{figure}

\begin{figure}[t!]
\centering
\hspace{-.2in}
   \begin{tabular}{cc}
        \begin{subfigure}{.3\textwidth}
          		\includegraphics[width=\textwidth]{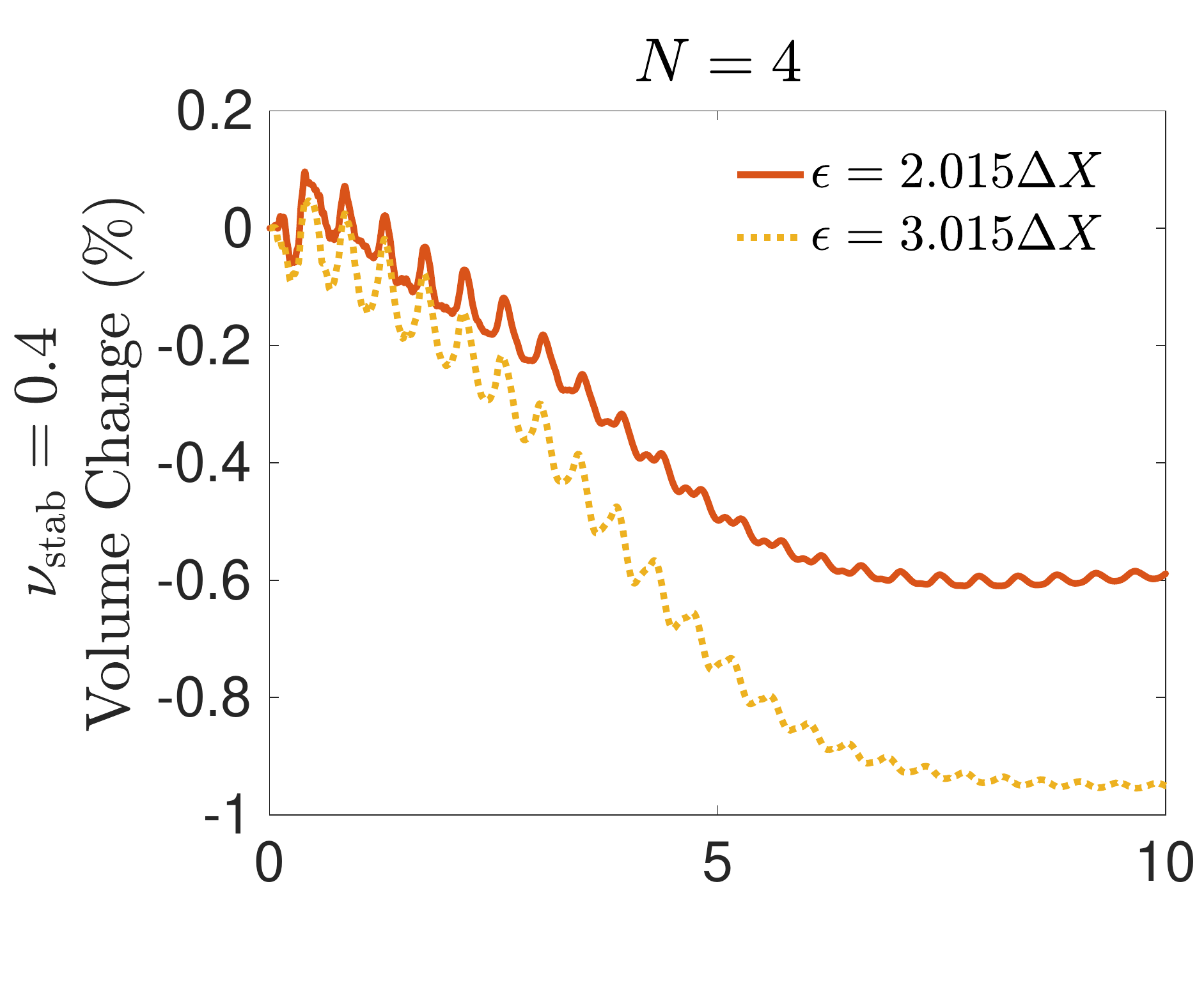}
        \end{subfigure} 
        \begin{subfigure}{.3\textwidth}
                \includegraphics[width=\textwidth]{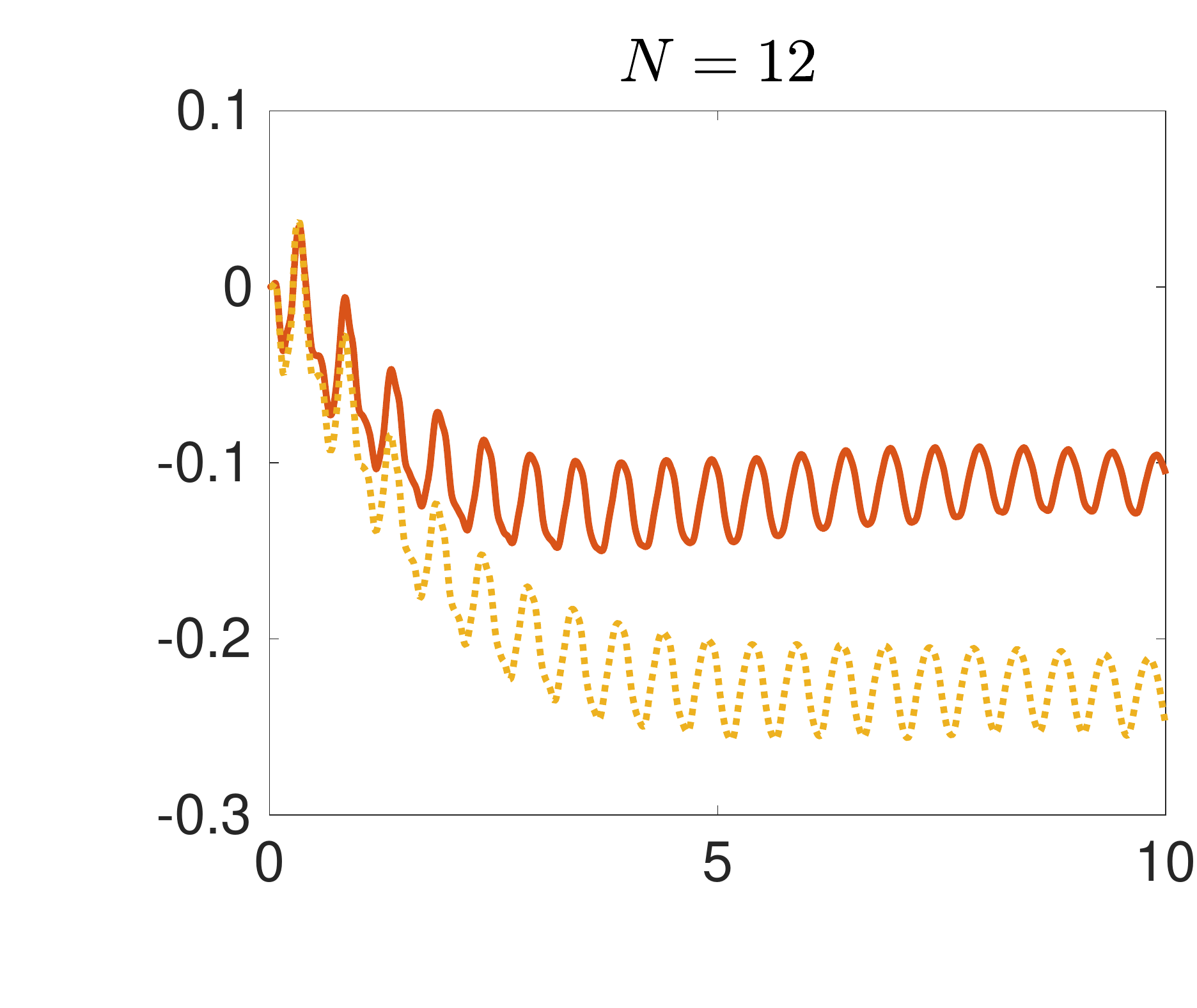}
        \end{subfigure}
        \begin{subfigure}{.3\textwidth}
                \includegraphics[width=\textwidth]{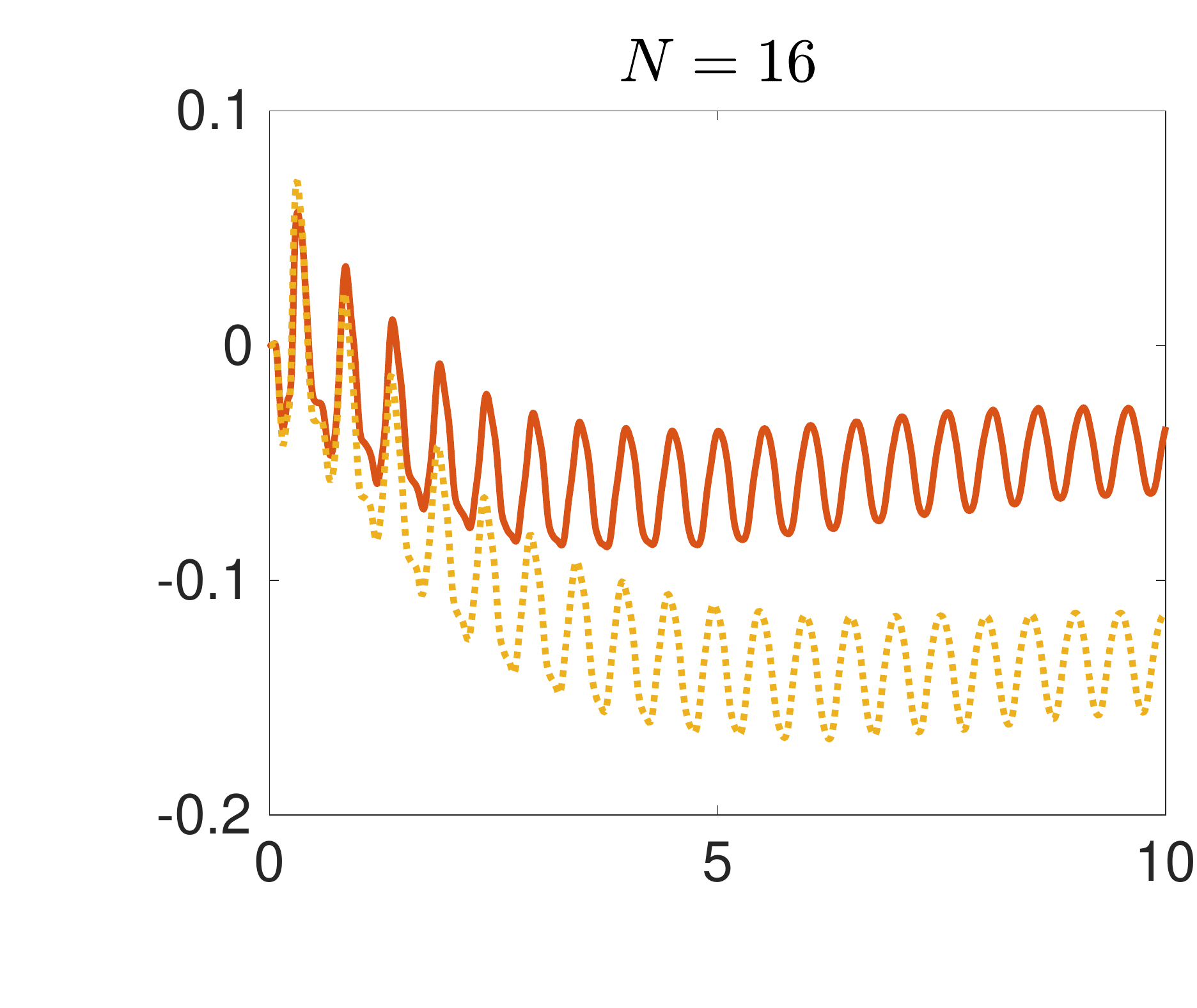}
        \end{subfigure} \\
         \begin{subfigure}{.3\textwidth}
          		\includegraphics[width=\textwidth]{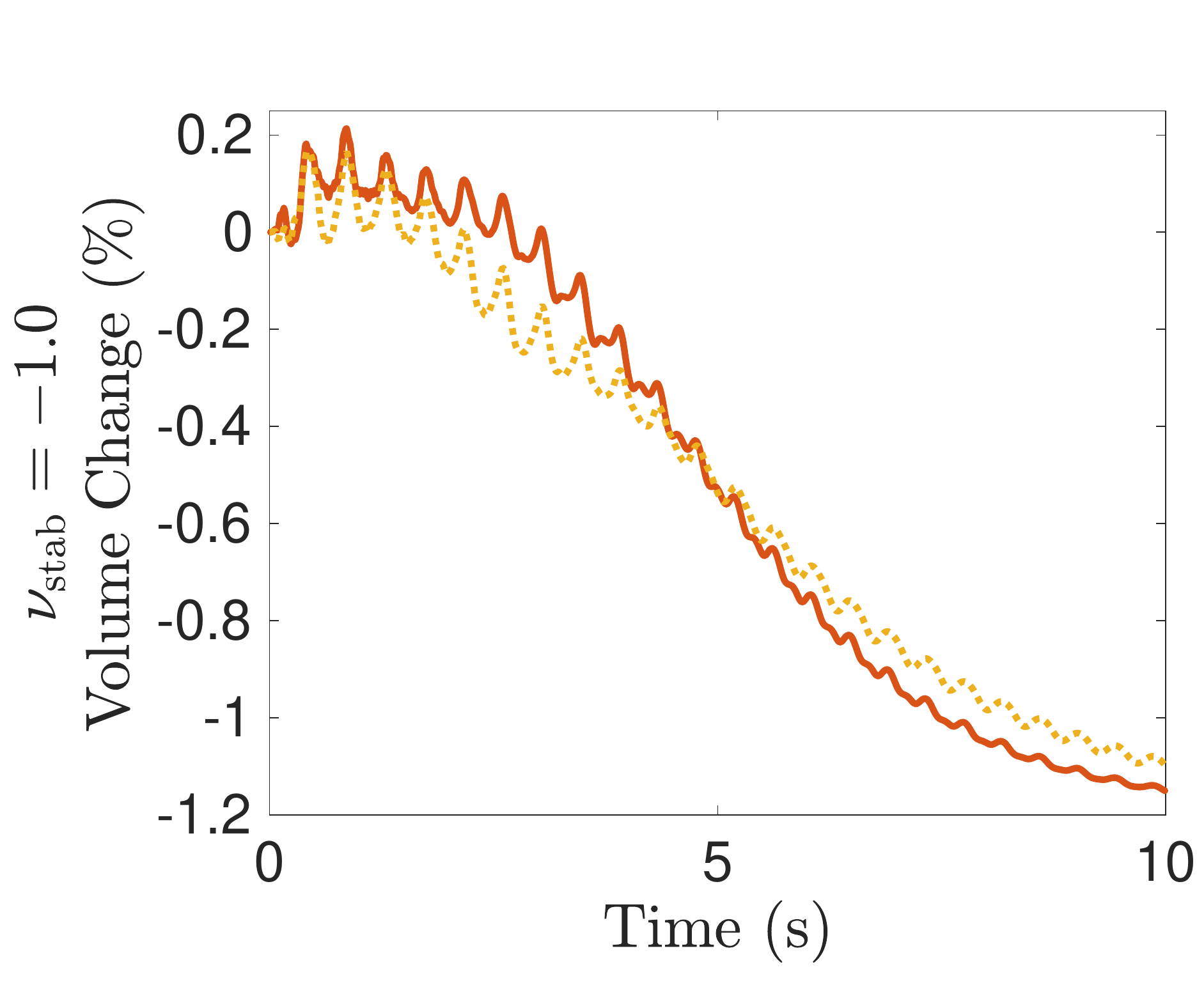}
        \end{subfigure} 
        \begin{subfigure}{.3\textwidth}
                \includegraphics[width=\textwidth]{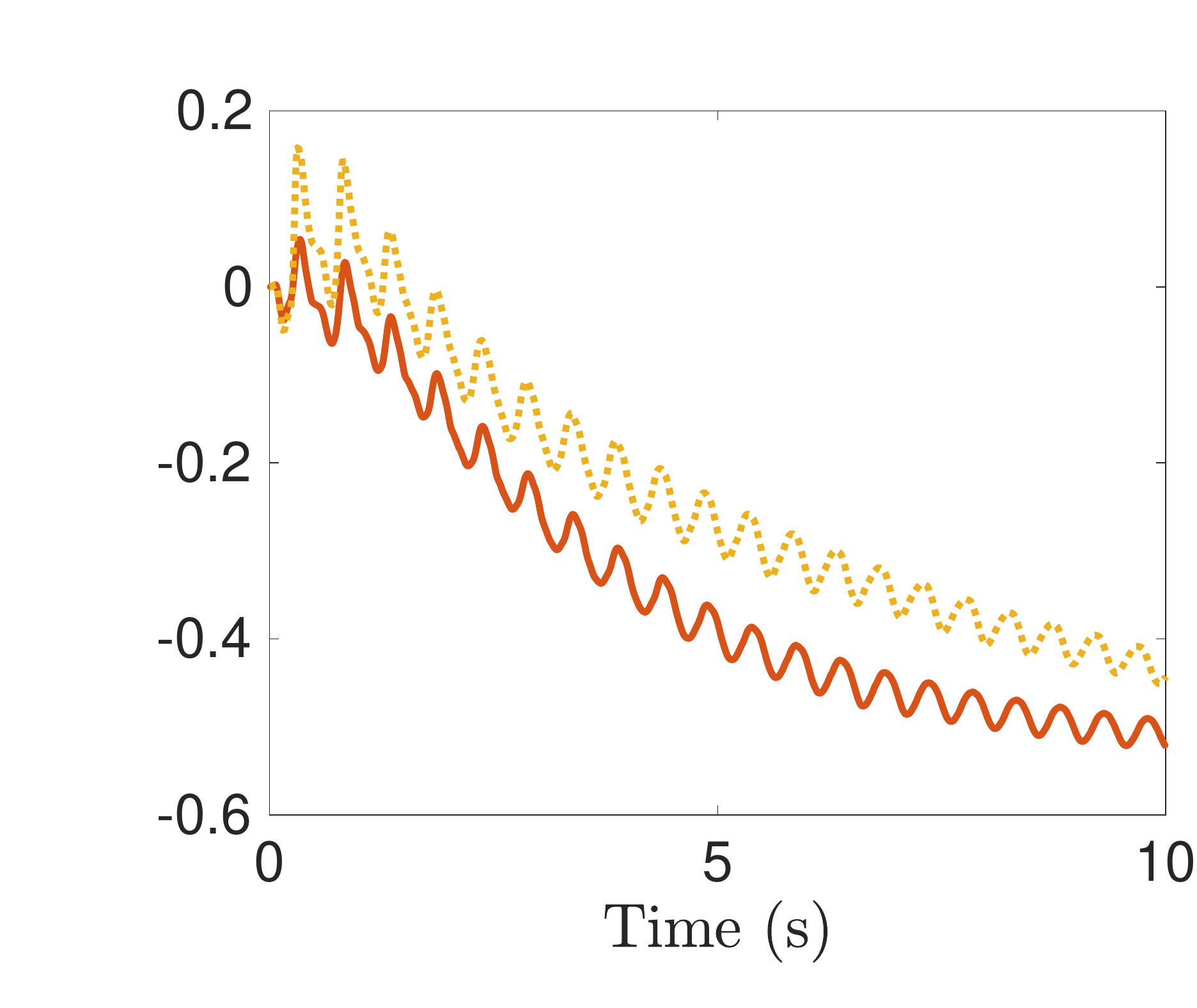}
        \end{subfigure}
        \begin{subfigure}{.3\textwidth}
                \includegraphics[width=\textwidth]{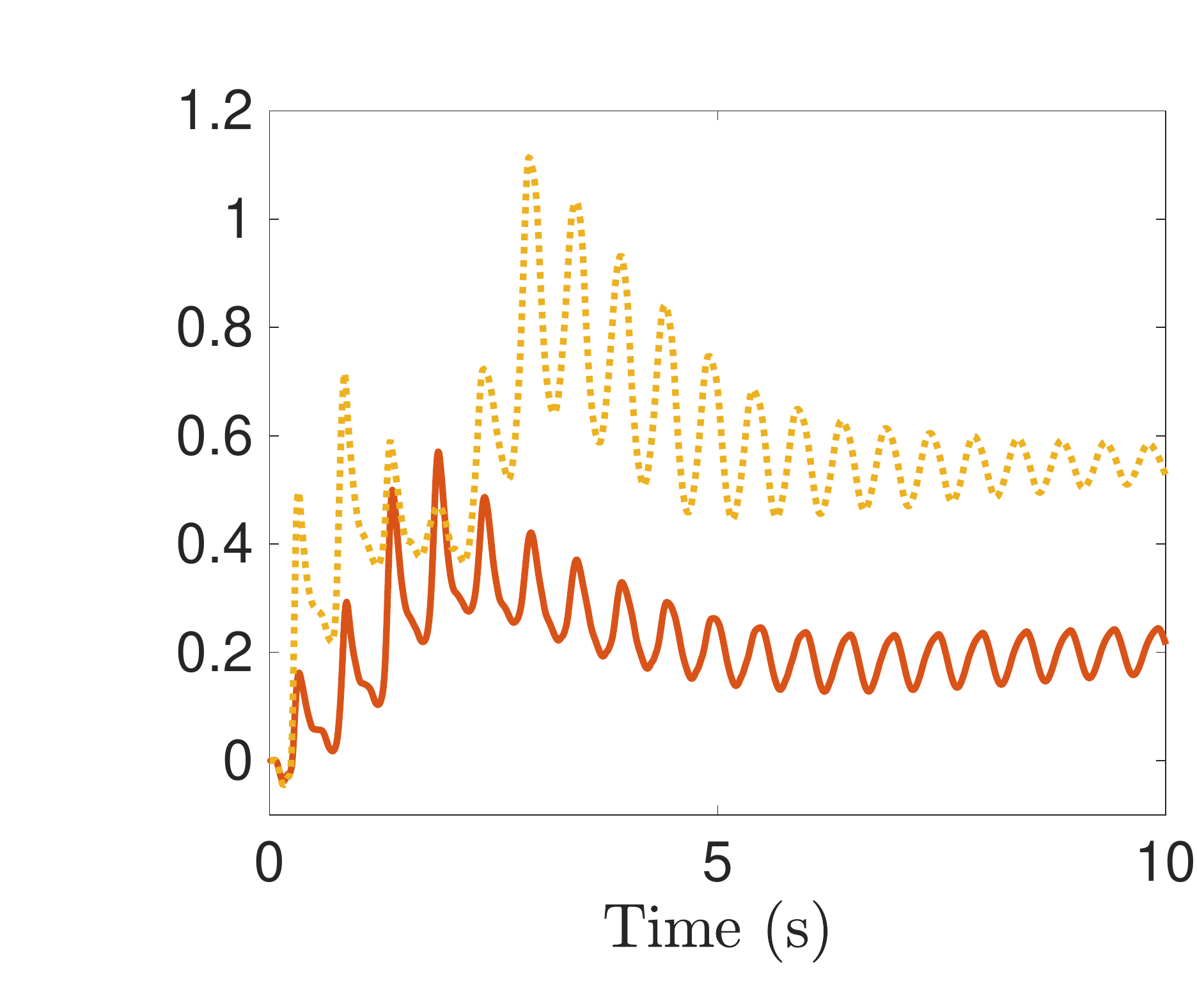}
        \end{subfigure} \\
    \end{tabular}
    \caption{Volume change of the band in Fig.~\ref{f:ELS_BND_schematics} for different choices of horizon size $\horizonsize$ and numerical Poisson's ratio $\nu_{\text{stab}}$ under grid refinement. $N=4$ corresponds to $165$ solid DoF, $N=12$ corresponds to $1261$ solid DoF, and $N = 16$ corresponds to $2193$ solid DoF.}
    \label{f:ELS_BND_dynamics_volume_change}
\end{figure}

\begin{figure}[t!]
\centering
   \begin{tabular}{cc}
        \begin{subfigure}{.4\textwidth}
          		\includegraphics[width=\textwidth]{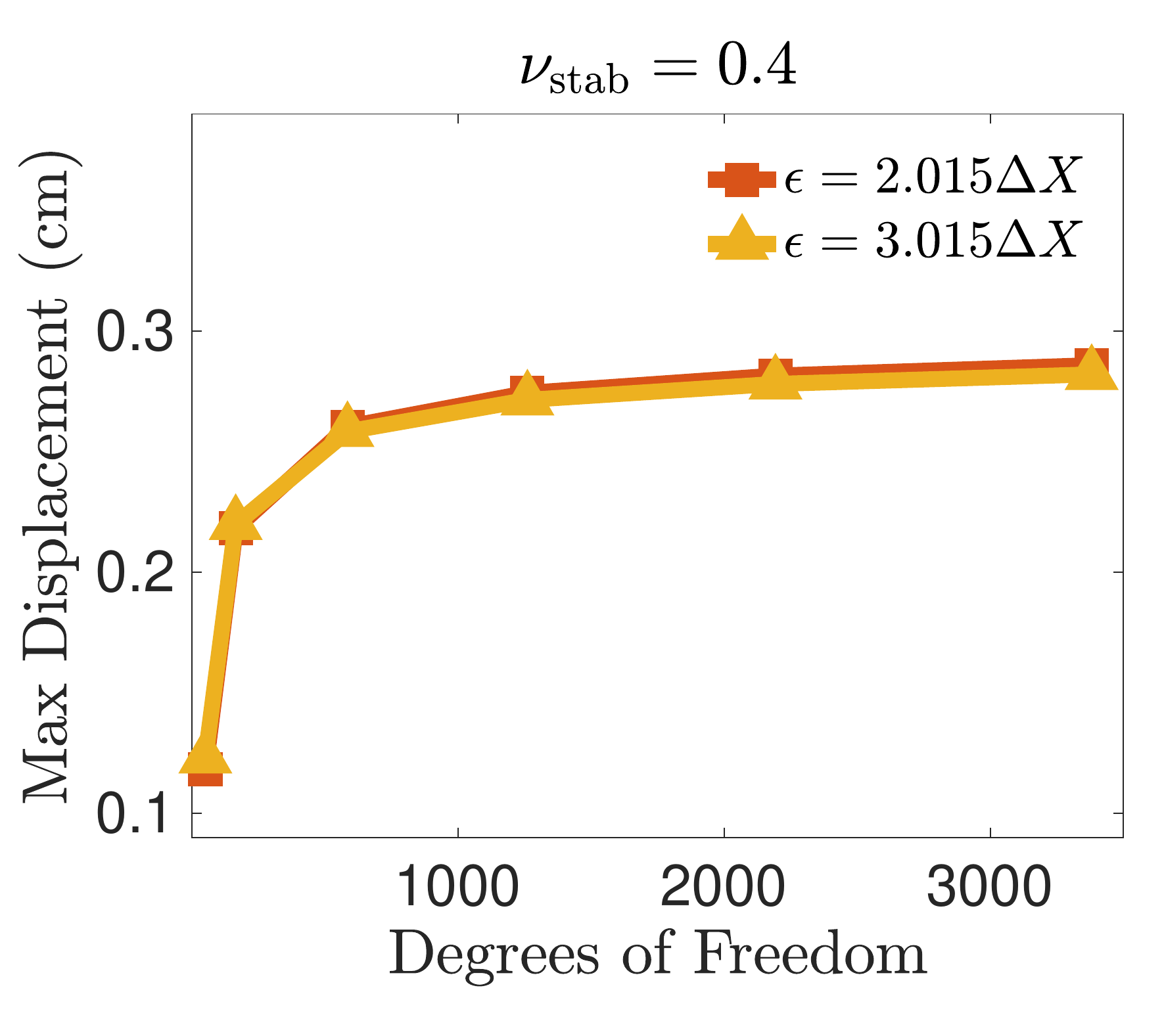}
          		\caption{}
          		 \label{f:ELS_BND_dynamics_max_disp}
        \end{subfigure} \hspace{.05\textwidth}
         \begin{subfigure}{.4\textwidth}
        			\includegraphics[width=\textwidth]{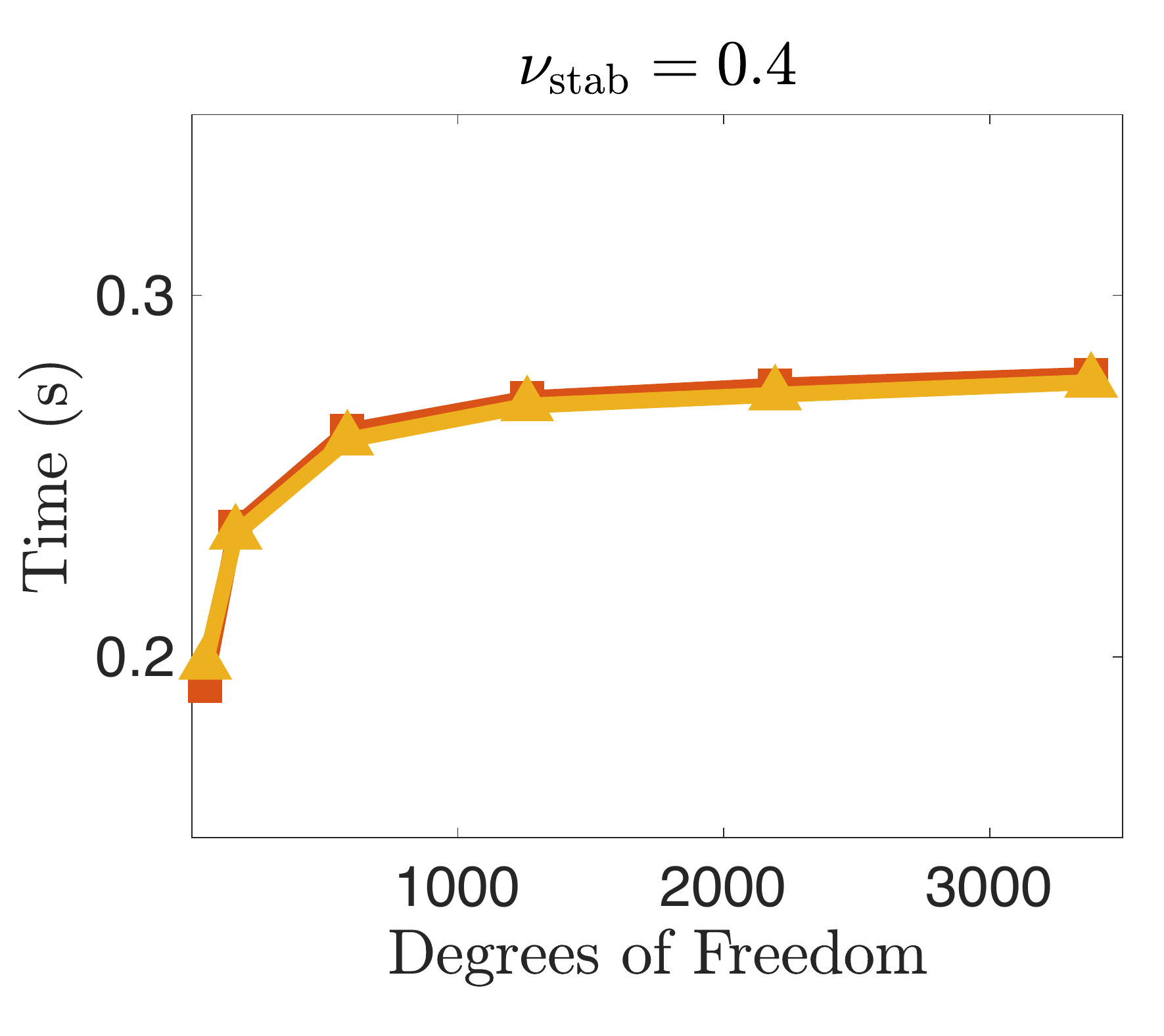}
          		\caption{}
          		 \label{f:ELS_BND_dynamics_max_time}
        \end{subfigure} 
    \end{tabular}
    \caption{(a) Maximum displacements of the dynamic version of the elastic band for different choices of horizon size $\horizonsize$ under grid refinement. The solid DoF range from $51$ to $3381$. (b) Time to reach the maximum displacement of the dynamic version of the elastic band for different choices of horizon size $\horizonsize$ under grid refinement. }
\end{figure}

\begin{figure}[t!]
\centering
   \begin{tabular}{cc}
        \begin{subfigure}{.4\textwidth}
          		\includegraphics[width=\textwidth]{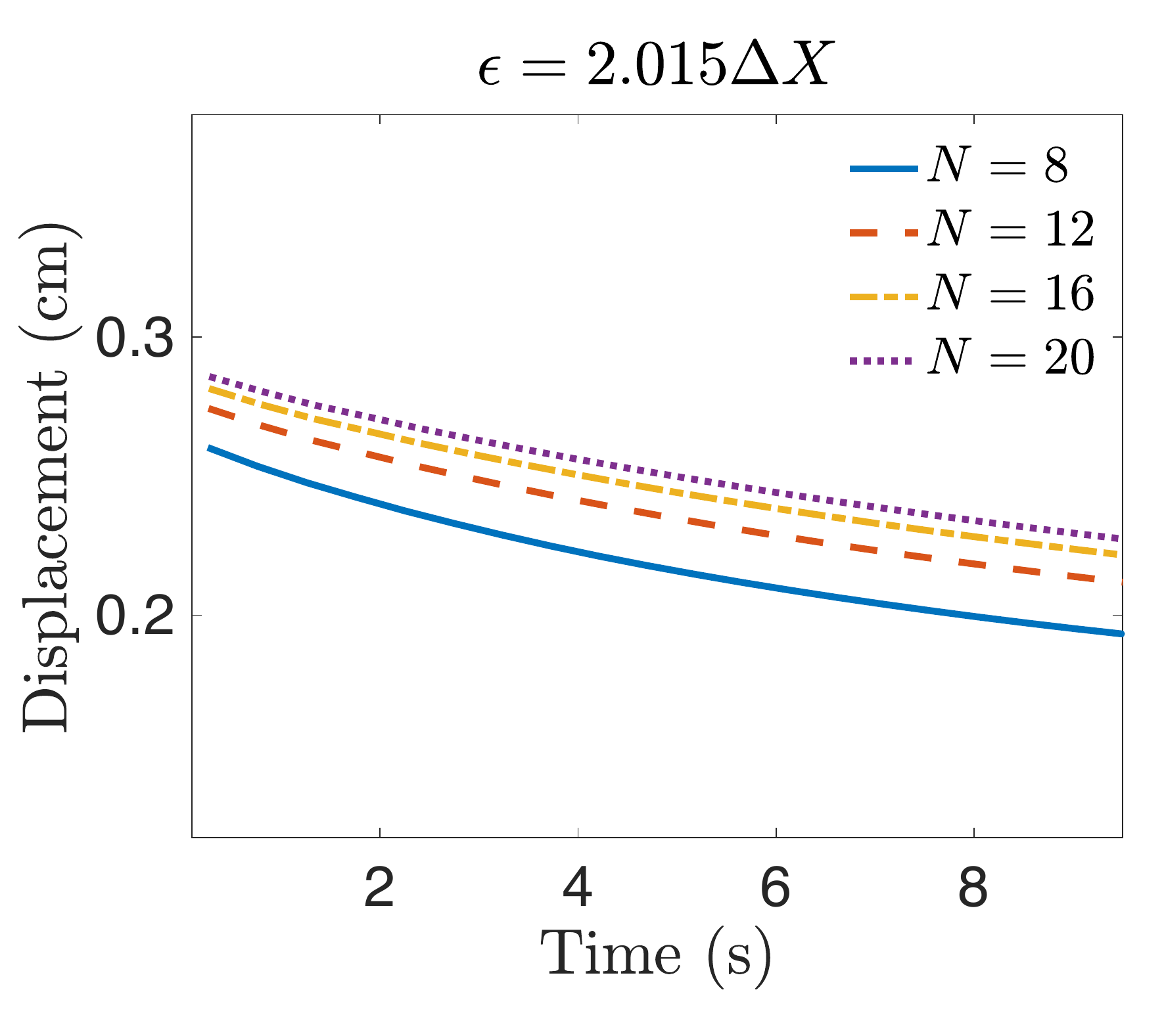}
          		\caption{}
        \end{subfigure} \hspace{.05\textwidth}
         \begin{subfigure}{.4\textwidth}
        			\includegraphics[width=\textwidth]{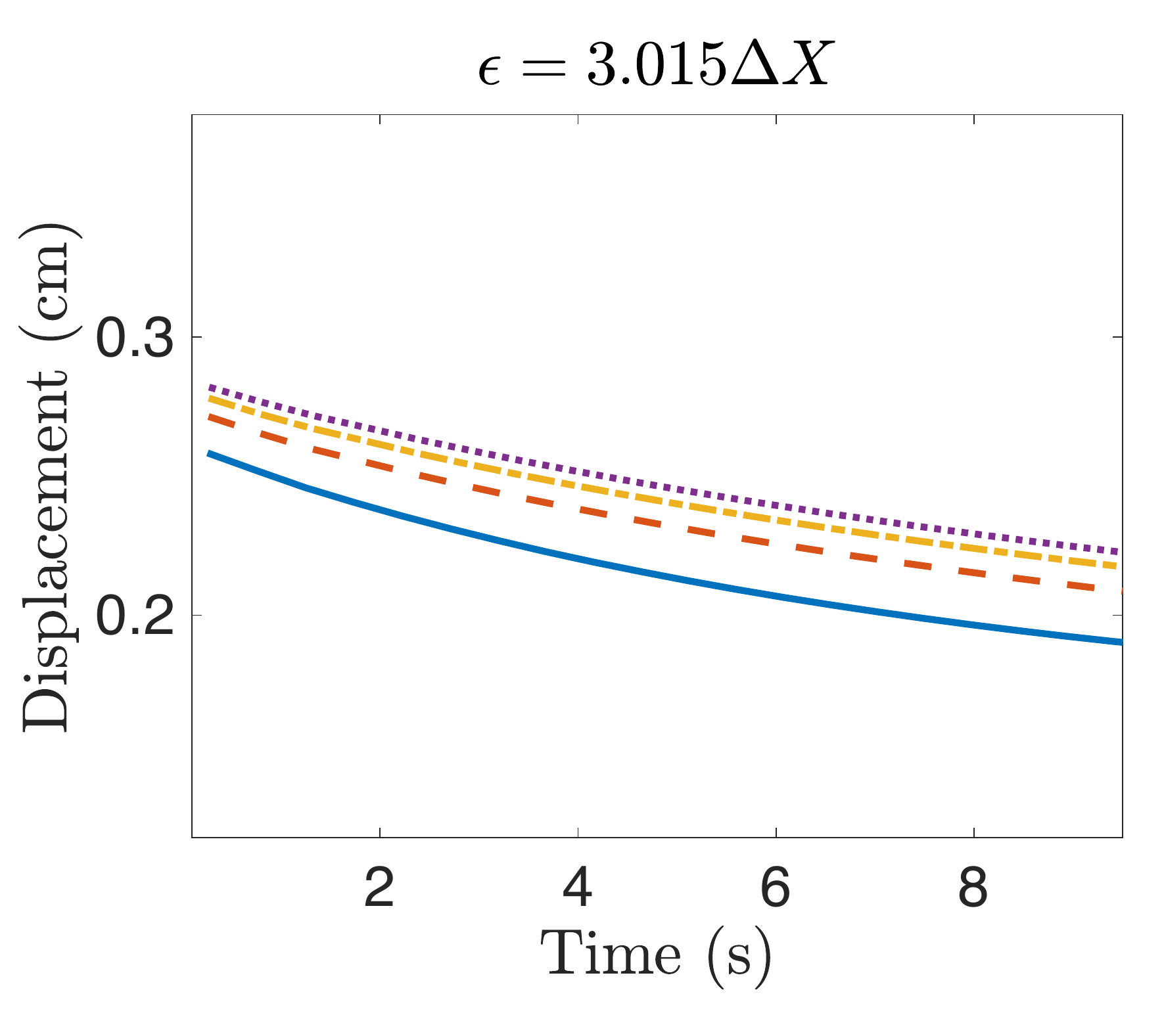}
          		\caption{}
        \end{subfigure} 
    \end{tabular}
    \caption{Upper envelopes of the oscillations in Fig.~\ref{f:ELS_BND_dynamics_displacement} for different choices of horizon size $\horizonsize$ under grid refinement with $\nu_{\text{stab}} = 0.4$. The solid DoF range from $585$ to $3381$. $N = 8$ corresponds to $585$ solid DoF, $N=12$ corresponds to $1261$ solid DoF, $N=16$ corresponds to $2193$ solid DoF, and $N=20$ corresponds to $3381$ solid DoF. }
 \label{f:ELS_BND_dynamics_upper_envelops}
\end{figure}

Fig.~\ref{f:ELS_BND_dynamics} shows the structural deformations of the band along with the Eulerian velocity field and the values of Jacobian of the non-local deformation tensor. 
Fig.~\ref{f:ELS_BND_dynamics_displacement} shows the transient behavior of the elastic band against time for various numerical Poisson's ratios $\nu_{\text{stab}}$ and peridynamic horizon sizes $\horizonsize$ under grid refinement. 
Fig.~\ref{f:ELS_BND_dynamics_volume_change} shows the volume change of the band under deformation for different choices of $\nu_{\text{stab}}$ and $\horizonsize$. 
The total volume change decreases under grid refinement as in the static case, and the range is comparable to the static case as well. 
With a larger value of the numerical bulk modulus, the volume change noticeably decreases. 
Fig.~\ref{f:ELS_BND_dynamics_max_disp} and Fig.~\ref{f:ELS_BND_dynamics_max_time} show that the maximum displacements of the oscillations and time to reach the maximum displacements converge under grid refinement. 
Fig.~\ref{f:ELS_BND_dynamics_upper_envelops} shows the upper envelops of the oscillations for different choices of $\horizonsize$ under grid refinement.

\subsection{Failure benchmarks}
\label{s:failure}
This section presents modified elastic band benchmark problems that allow bond breakage to simulate the fluid-driven deformations of a material that can experience damage and, ultimately, failure.

\subsubsection{Rupture of an elastic band}
This benchmark considers dynamic material deformations and fracture of the elastic band under fluid-driven forces. 
In this benchmark, the critical bond stretch is set to $\sc = 4.5$ to demonstrate the effectiveness of simulating crack initiation and propagation using the IPD method. 
In general, the critical bond stretch of a material must be experimentally determined.
The critical bond stretch $\sc$ used here is determined based on preliminary simulations.
The pressure loading is three times larger than the value used in Sec.~\ref{s:Benchmark_ELS} on each side; $\vec{\bbsigma}^{\text{f}} (\x,t) \bm{n} (\x) = \bm{h}(t)$, in which $\bm{h}(t) = \left(- 30 , 0\right) \, \frac{\mathrm{dyn}}{\mathrm{cm}^2}$ and $\bm{h}(t) = \left(30,0\right) \, \frac{\mathrm{dyn}}{\mathrm{cm}^2}$ on the left and right, respectively. 
The final simulation time is set to $T_{\text{f}} = 0.25 \, \mathrm{s}$.
Otherwise, we use the same parameters as the non-failure case of the dynamic elastic band benchmark.  
The horizon size is set to $\horizonsize = 3.015 \Delta X$, as suggested in the PD literature \cite{behera2020peridynamic}.

\begin{figure}[]
\centering
\includegraphics[width=\textwidth]{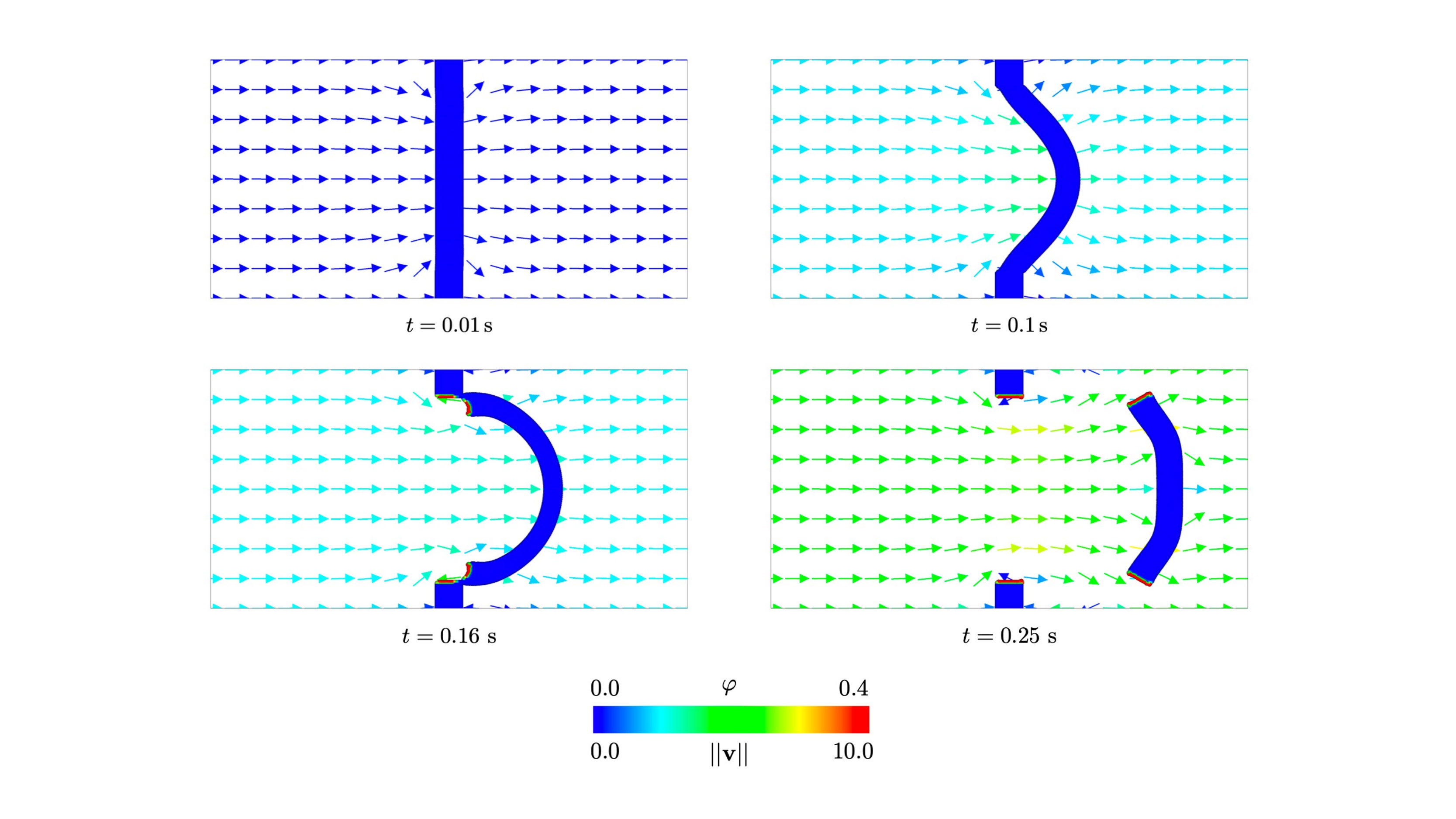}
\vspace{-.25in}
    \caption{Dynamic failure process of the elastic band with the local damage in NOSB-PD along with the corresponding Eulerian velocity field. Note that $\varphi = 0$ implies that all initial bonds are connected, and $\varphi = 1$ implies that all initial bonds are disconnected. The deformations are computed using $3381$ solid DoF, $\horizonsize = 3.015 \Delta X$, and $\nu_{\text{stab}} = 0.4$.}
    \label{f:ELS_BND_failure_dynamics}
\end{figure}

\begin{figure}[]
\centering
   \begin{tabular}{cc}
        \begin{subfigure}{.4\textwidth}
          		\includegraphics[width=\textwidth]{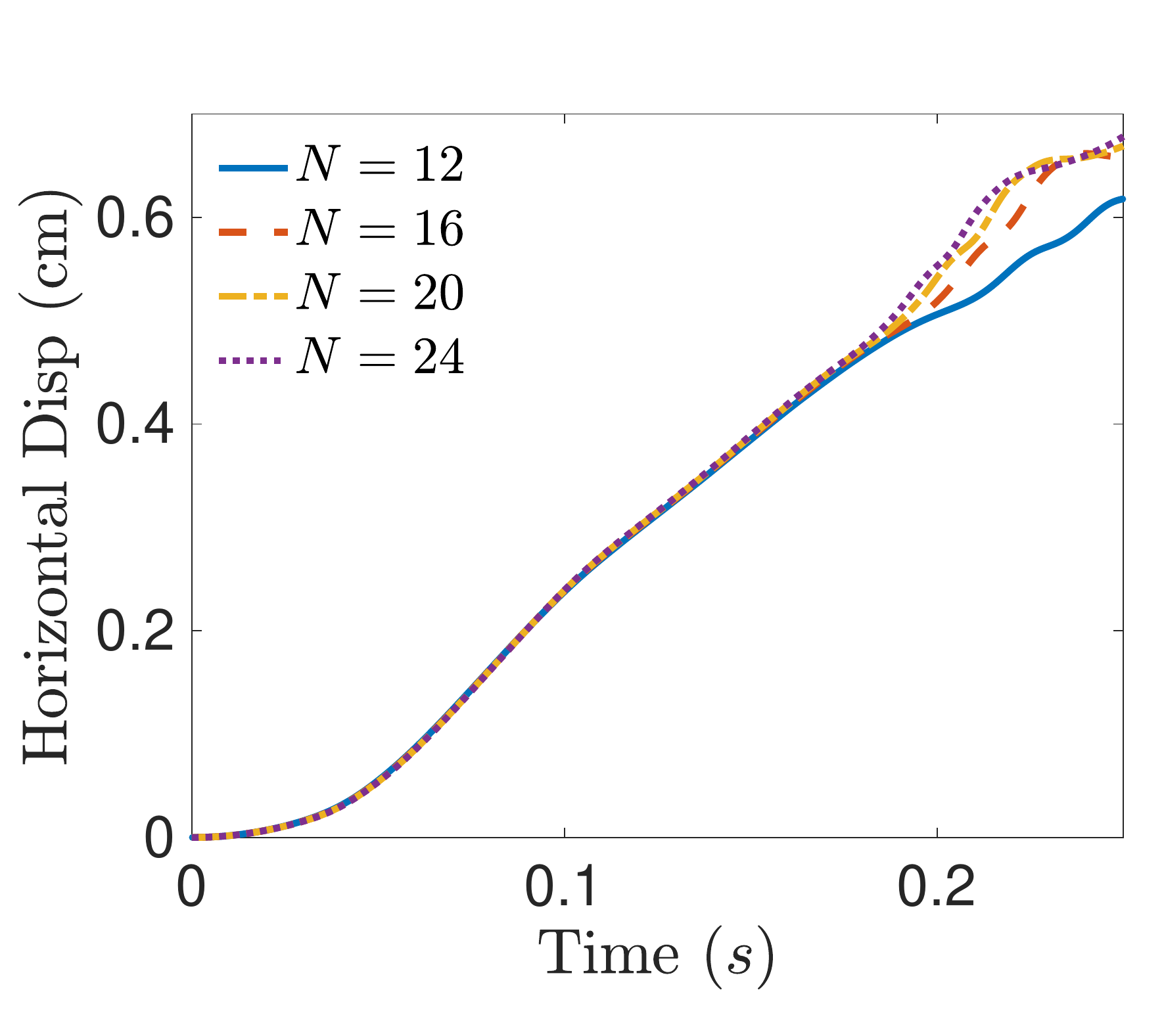}
          		\caption{}
          		 \label{f:ELS_BND_rupture_disp_1}
        \end{subfigure} \hspace{.05\textwidth}
         \begin{subfigure}{.4\textwidth}
        			\includegraphics[width=\textwidth]{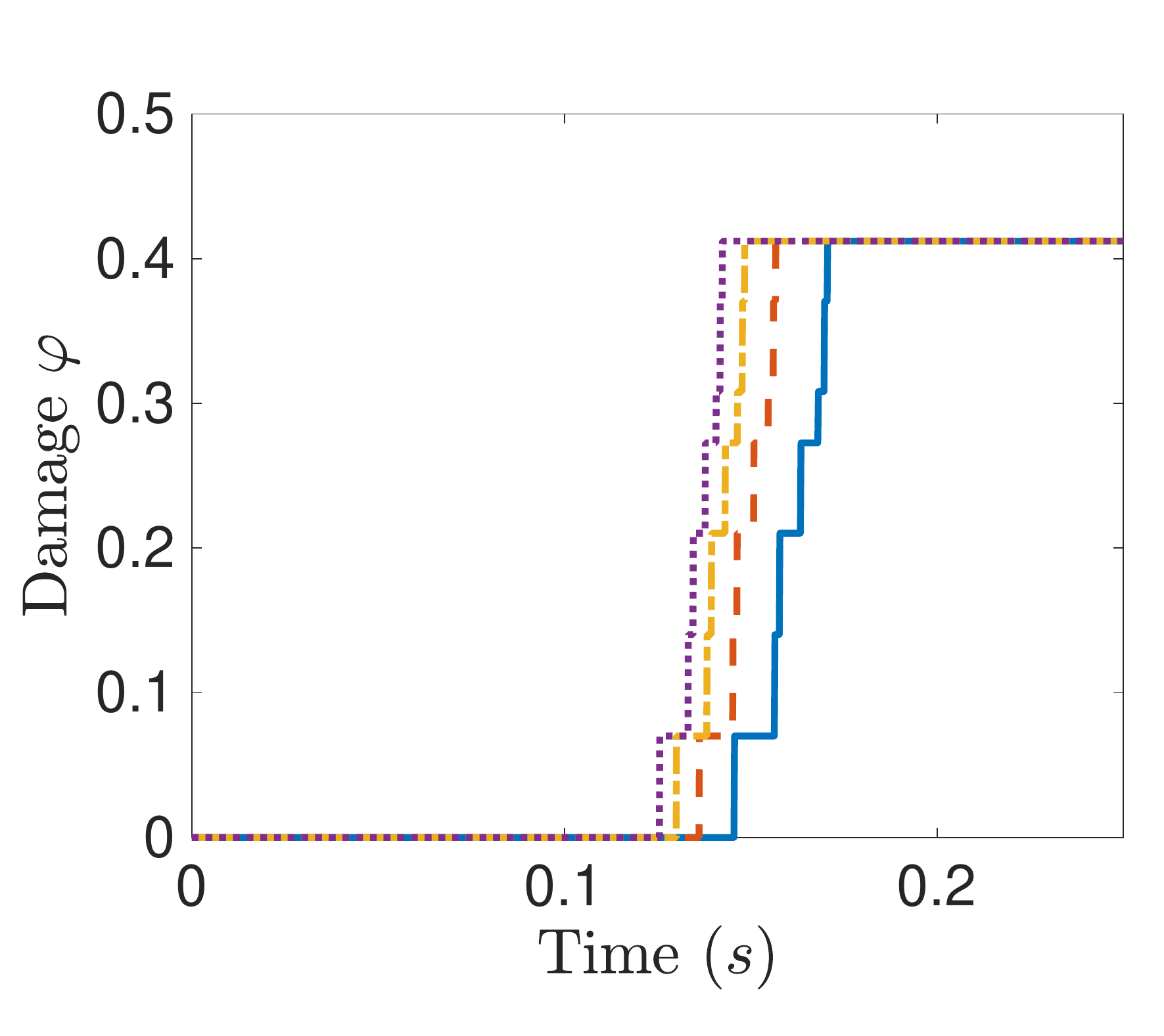}
          		\caption{}
          		 \label{f:ELS_BND_rupture_damage}
        \end{subfigure} 
    \end{tabular}
    \caption{(a): Horizontal displacements of the point of interest, highlighted in Fig.~\ref{f:ELS_BND_schematics}, under grid refinement. (b): Local damage growth at the top left corner of the detached band during the failure process under grid refinement. $N=12$ corresponds to $1261$ solid DoF, $N=16$ corresponds to $2193$ solid DoF, $N=20$ corresponds to $3381$ solid DoF, and $N=24$ corresponds to $4825$ solid DoF.}
\end{figure}

Fig.~\ref{f:ELS_BND_failure_dynamics} shows the crack nucleation and propagation of the dynamic version of the elastic band benchmark with an Eulerian velocity field. 
The crack formulation is initiated near the junctions between the fixed blocks and the band, and the band gets entirely detached from the block when the bonds exceed the critical bond stretch $\sc$.
Fig.~\ref{f:ELS_BND_rupture_disp_1}  shows the horizontal displacements of the point of interest, highlighted in Fig.~\ref{f:ELS_BND_schematics}, for different grid spacings. 
Fig.~\ref{f:ELS_BND_rupture_damage} shows the local damage growth at the top left corner of the detached band during the failure process under grid refinement. 

\subsubsection{Elastic band with a notch}
\label{s:Benchmark_ELS_notch}
\begin{figure}[]
\centering
    \includegraphics[width=.8\textwidth]{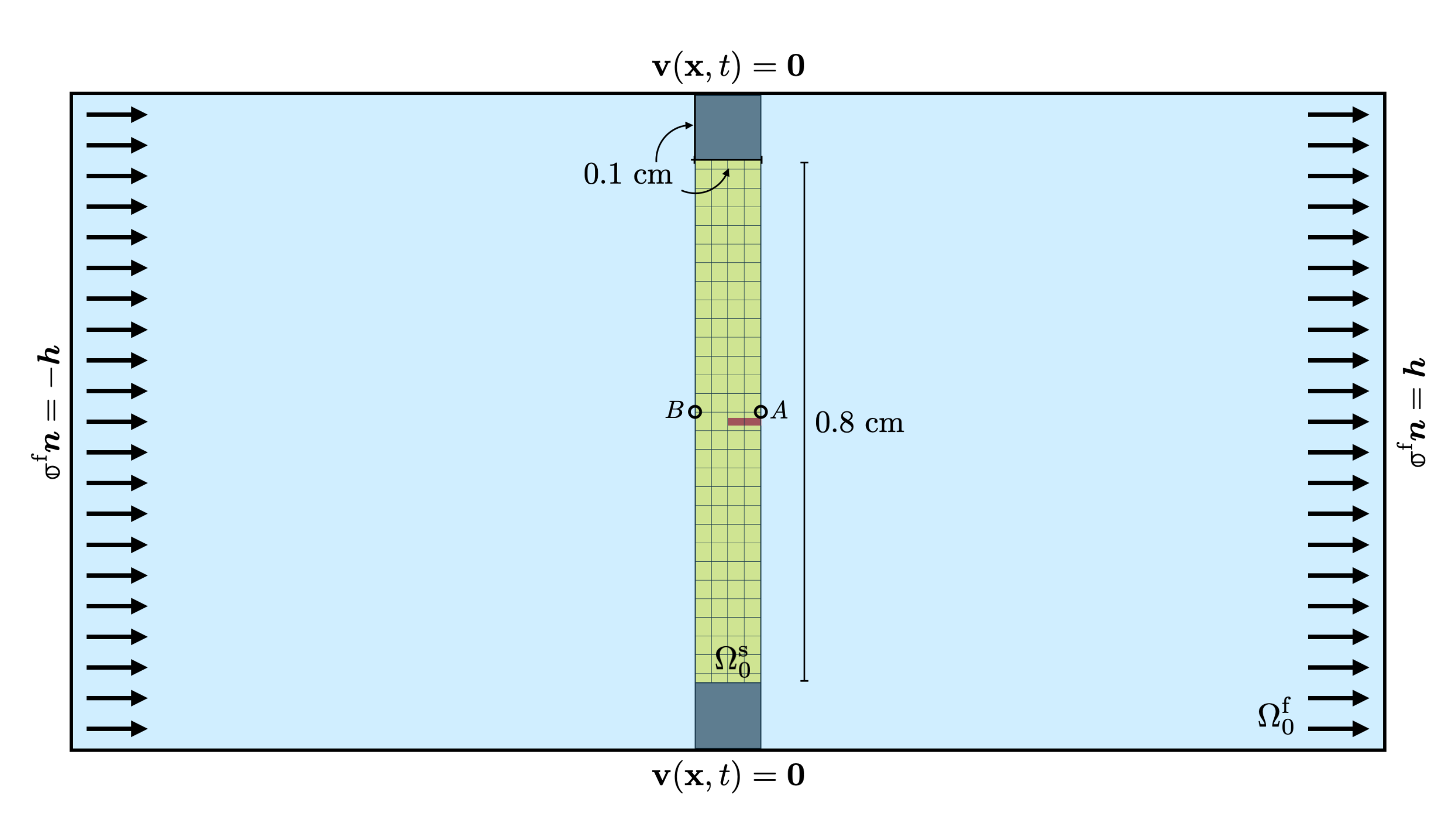}
    \caption{Schematic diagram for the failed elastic band benchmark (Sec.~\ref{s:Benchmark_ELS_notch}). The initial configurations of the immersed structure and a fluid are denoted by $\Omega_0^{\text{s}}$ and $\Omega_0^{\text{f}}$, respectively. The entire computational domain is $\Omega = \Omega_0^{\text{s}} \cup \Omega_0^{\text{f}}$. Zero fluid velocity is enforced on the top and bottom boundaries of the computational domain, and fluid traction boundary conditions are applied to the left and right boundaries. Fluid tranction is set to $\bm{h}(t) = \left(20, 0\right) \, \frac{\mathrm{dyn}}{\mathrm{cm}^2}$.}
    \label{f:ELS_BND_notch_schematics}
\end{figure}

We next consider the dynamics of an elastic band with a pre-existing crack. 
A notch is placed on the center-right of the band with the length of $0.05 \, \mathrm{cm}$. 
The pressure loading is set to $\vec{\bbsigma}^{\text{f}} (\x,t) \bm{n} (\x) = \bm{h}(t)$, in which $\bm{h}(t) = \left(- 20 , 0\right) \, \frac{\mathrm{dyn}}{\mathrm{cm}^2}$ and $\bm{h}(t) = \left(20,0\right) \, \frac{\mathrm{dyn}}{\mathrm{cm}^2}$ on the left and right, respectively, with a zero loading time. 
In this benchmark, the critical bond stretch is set to $\sc = 4.5$ based on preliminary tests. 
Otherwise, we use the same parameters as the non-failure case of the dynamic elastic band benchmark. 
To obtain symmetric fracture in the middle of the band, we use even numbers of Lagrangian points in the vertical direction of the elastic band. 
Fig.~\ref{f:ELS_BND_notch_schematics} provides a schematic of this test case. 
The horizon size is set to $\horizonsize = 3.015 \Delta X$, as suggested in the PD literature \cite{behera2020peridynamic}. 
The final simulation time is set to $T_{\text{f}} = 0.3 \, \mathrm{s}$.

\begin{figure}[]
\centering
\includegraphics[width=\textwidth]{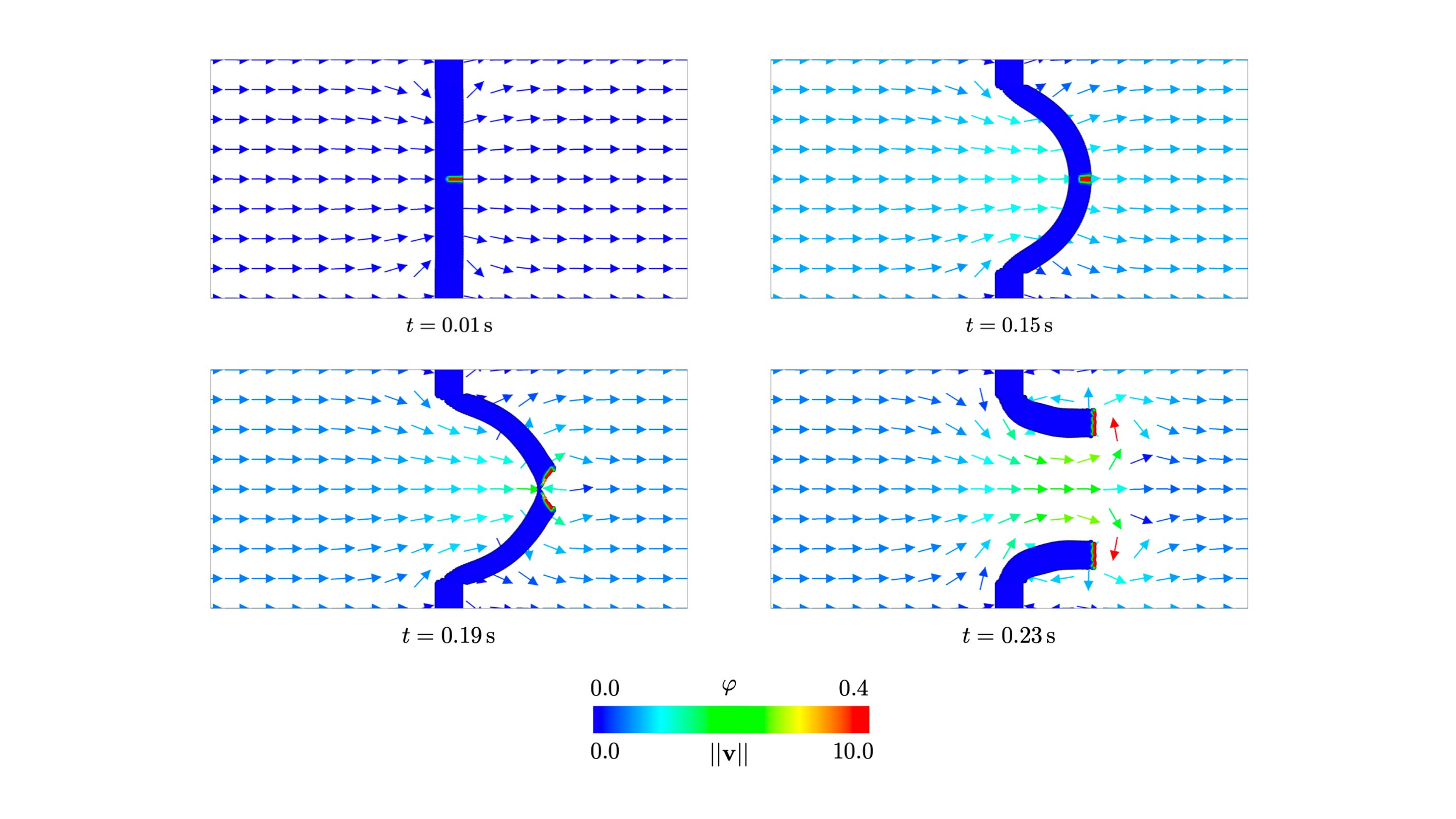}
\vspace{-.25in}
\caption{Dynamic failure process of the elastic band with the local damage in NOSB-PD along with the corresponding Eulerian velocity field. Note that $\varphi = 1$ indicates all initially connected bonds are broken. The deformations are computed using $2210$ solid DoF, $\horizonsize = 3.015 \Delta X$, and $\nu_{\text{stab}} = 0.4$.}
\label{f:ELS_BND_notch_dynamics}
\end{figure}

\begin{figure}[]
\centering
   \begin{tabular}{cc}
        \begin{subfigure}{.4\textwidth}
          		\includegraphics[width=\textwidth]{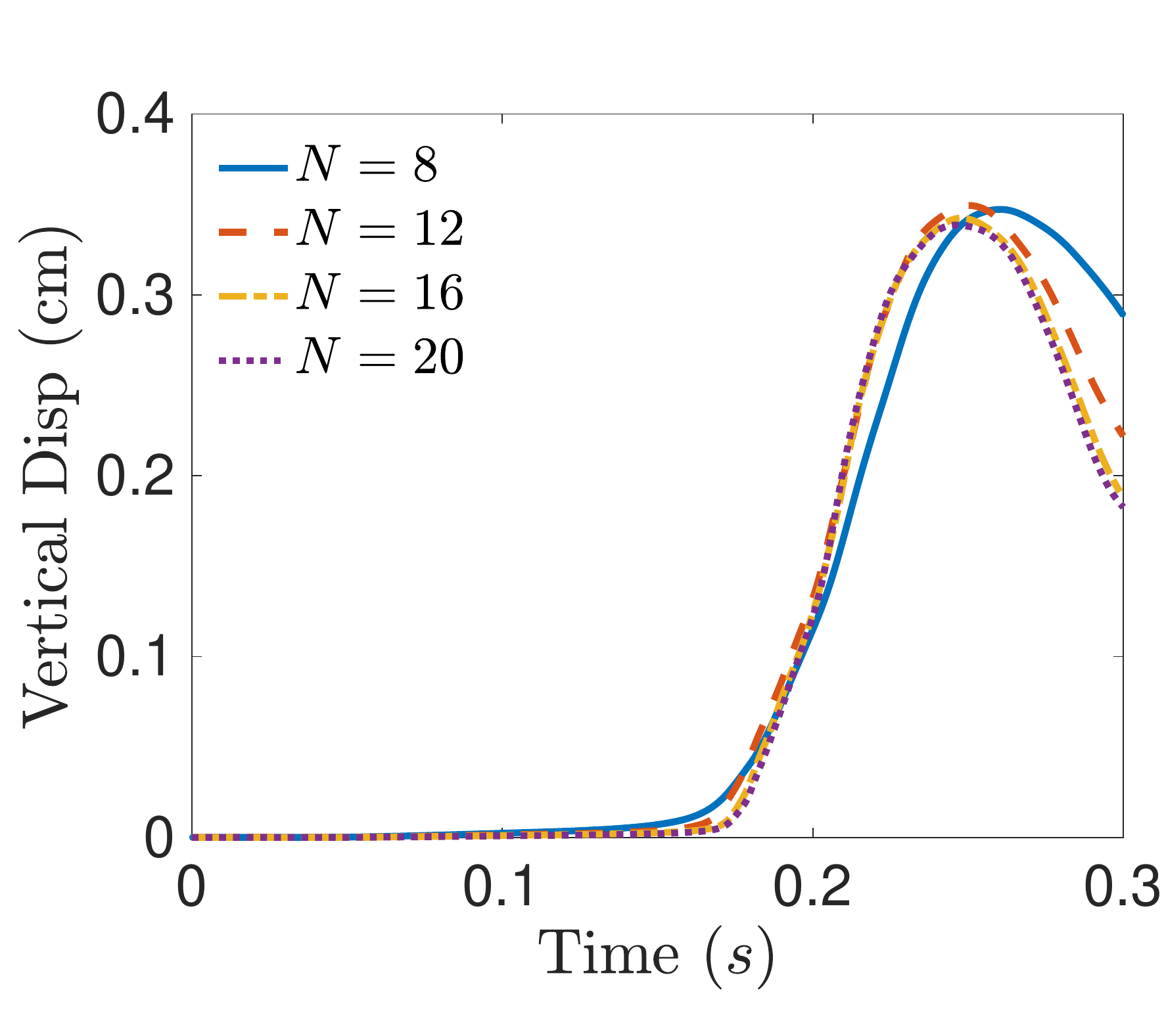}
          		\caption{}
          		 \label{f:ELS_BND_notch_disp_1}
        \end{subfigure} \hspace{.05\textwidth}
         \begin{subfigure}{.4\textwidth}
        			\includegraphics[width=\textwidth]{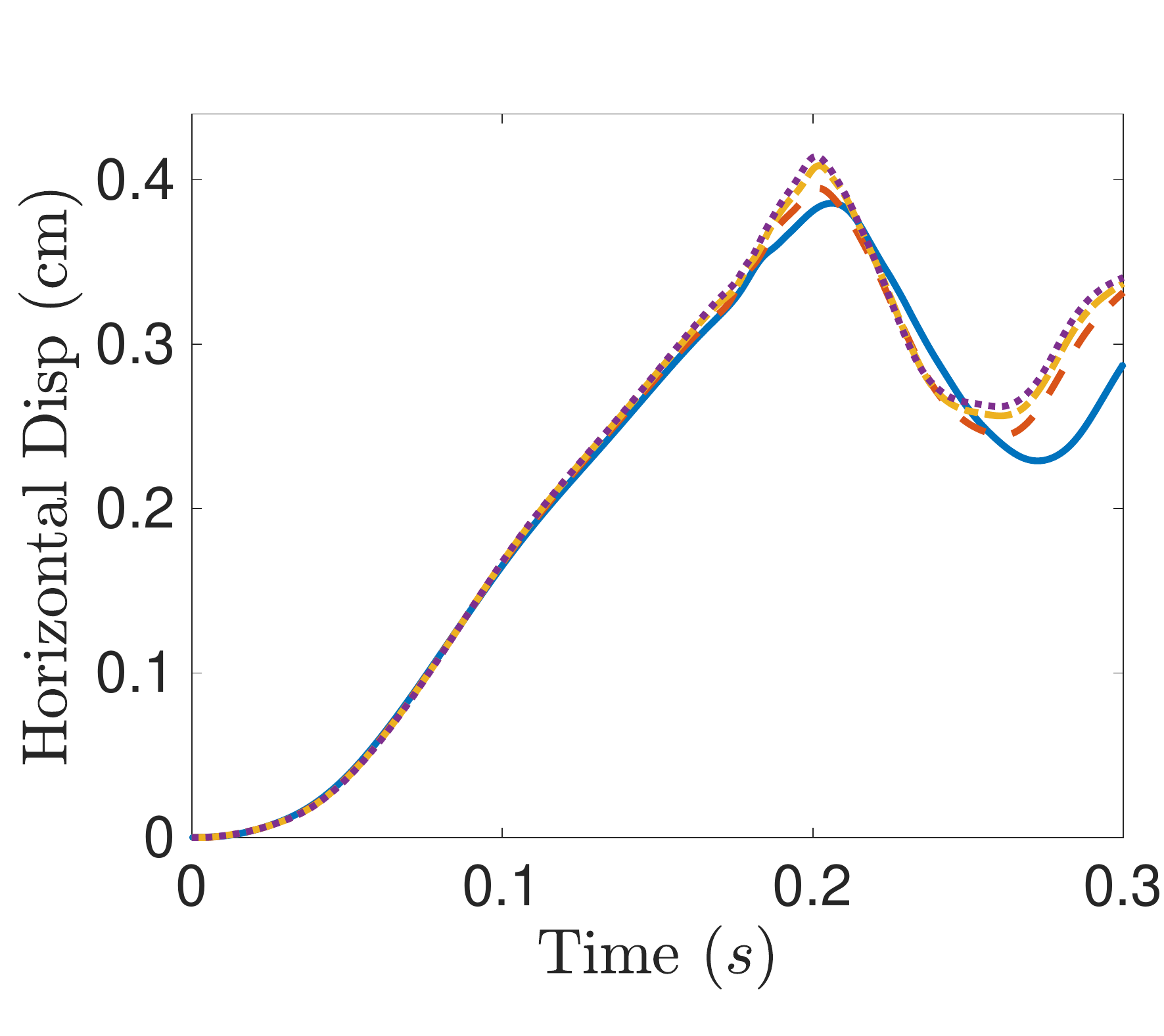}
          		\caption{}
          		 \label{f:ELS_BND_notch_disp_2}
        \end{subfigure} 
    \end{tabular}
    \caption{(a): Vertical displacements of the point of interest A, highlighted in Fig.~\ref{f:ELS_BND_notch_dynamics}, under grid refinement. (b): Horizontal displacements of the point of interest A under grid refinement. $N = 8$ corresponds to $594$ solid DoF, $N=12$ corresponds to $1274$ solid DoF, $N=16$ corresponds to $2210$ solid DoF, and $N=20$ corresponds to $3402$ solid DoF.}
\end{figure}

\begin{figure}[]
\centering
\includegraphics[width=0.4\textwidth]{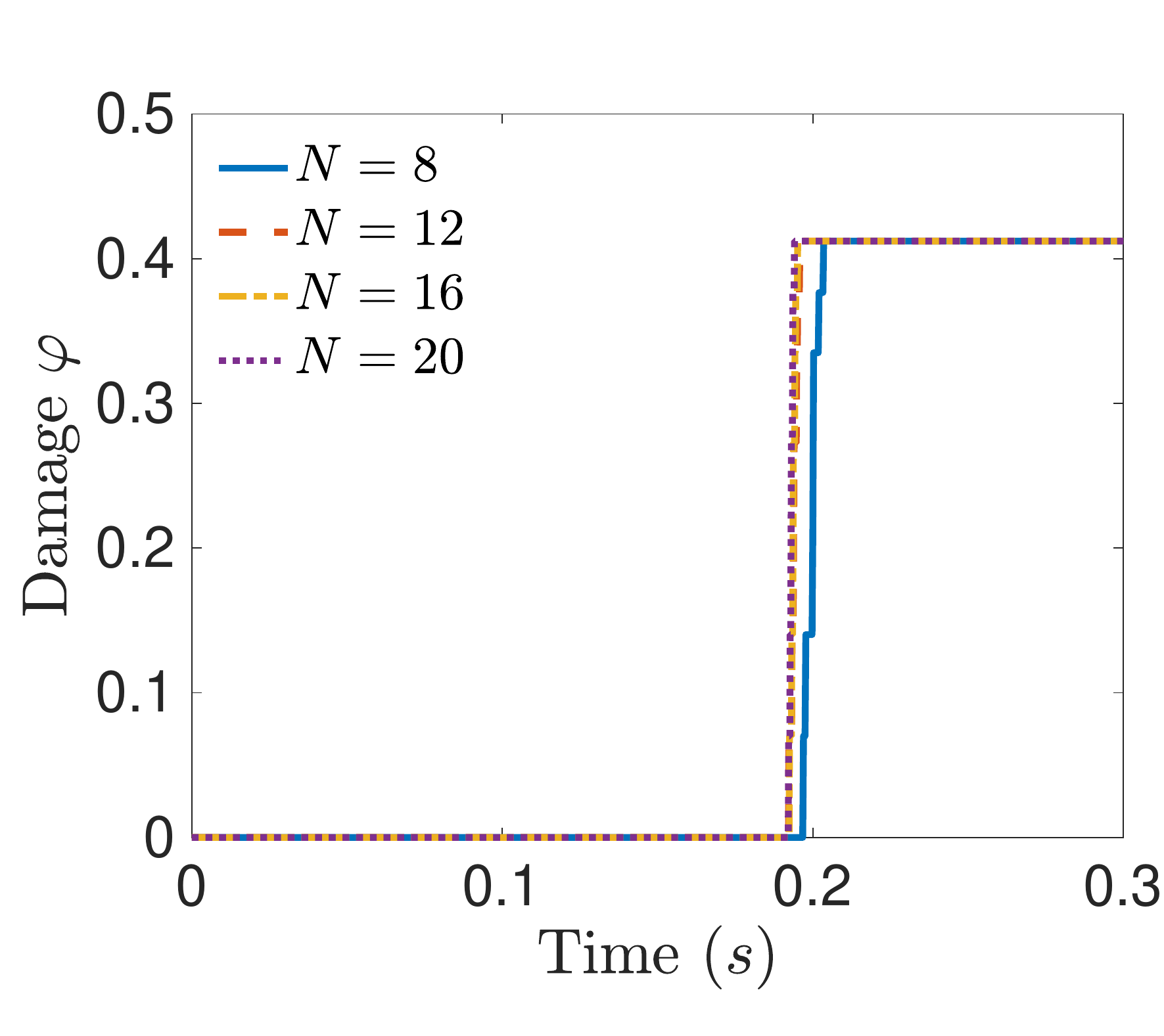}
\caption{Local damage growth at the point of interest B, highlighted in Fig.~\ref{f:ELS_BND_notch_dynamics}, during the failure process. $N = 8$ corresponds to $594$ solid DoF, $N=12$ corresponds to $1274$ solid DoF, $N=16$ corresponds to $2210$ solid DoF, and $N=20$ corresponds to $3402$ solid DoF. Note that $\varphi = 0$ implies that all initial bonds are connected and $\varphi = 1$ implies that all initial bonds are disconnected.}
\label{f:ELS_BND_notch_stress}
\end{figure}


Fig.~\ref{f:ELS_BND_notch_dynamics} shows the dynamic deformation and failure process of the elastic band along with the damage parameter $\varphi$ and Eulerian velocity vectors. 
The band undergoes a large deformation before the crack propagates, and it ultimately breaks the bond connectivity and completely ruptures. 
Fig.~\ref{f:ELS_BND_notch_disp_1} shows the vertical displacements of the point $A$ in Fig.~\ref{f:ELS_BND_notch_dynamics} under grid refinement, and Fig.~\ref{f:ELS_BND_notch_disp_2} shows horizontal displacements of the point $A$ under grid refinement. 
Fig.~\ref{f:ELS_BND_notch_stress} shows the local damage growth at the point of interest $B$ during the failure process under grid refinement. 


\section{Discussion and conclusion}

This paper develops an immersed peridynamics method to simulate FSI with material models that can experience material damage and failure. 
It uses non-ordinary state-based peridynamics to determine the internal body forces, which allows nonlinear material models of the immersed structural body that can accommodate discontinuities (i.e., crack formulation and propagation). 
Numerical tests consider both classical quasi-static benchmarks adapted from the solid mechanics literature \cite{reese1999new,cook1974improved,bonet1997nonlinear} and fully dynamic FSI benchmarks \cite{vadala2020stabilization}. 
Our numerical results demonstrate the constitutive correspondence of nonlinear hyperelastic material models for non-failure benchmarks, as detailed in Sec.~\ref{s:non-failure}, and show that the IPD method yields comparable accuracy under grid refinement that offered by a stabilized FE method \cite{reese1999new} and the IFED method \cite{vadala2020stabilization}. 
We obtain accurate results of overall deformations and volume conservation. 
We also test crack initiation, growth, and fracture in the immersed structural body by fluid stresses, as detailed in Sec.~\ref{s:failure}, and simulation results show $\horizonsize$-convergence as in the non-failure benchmarks. 
Moreover, we investigate the effect of the size of the peridynamic horizon using standard solid mechanics and FSI benchmark studies. 
When the IPD method is used with classical solid mechanics benchmarks, results are relatively insensitive to the size of the peridynamic horizons. 
In contrast, the elastic band benchmarks requires at least $\horizonsize = 2.015 \Delta X$ for non-failure tests and $\horizonsize = 3.015 \Delta X$ for failure tests. 

An interesting finding of this work is that the IPD method with a proper volume stabilization does not appear to suffer from spurious zero-energy modes.
Our numerical experiments suggest that the methodology developed herein reduces the effect of zero-energy modes by including the volumetric stabilization term to the strain energy functional; however, we currently lack a rigorous theoretical understanding of this empirical observation.
Zero-energy modes in deformed material bodies are observed in numerical discretizations of the PD correspondence model in the PD literature \cite{silling2017stability}.
These modes were initially regarded as discretization issuses, and various numerical treatments have been developed to reduce such instabilities \cite{silling2017stability, gu2017effective, luo2018stress, chowdhury2019modified}. 
However, Tupek and Radobitzky demonstrated that the instabilities are caused by the definition of the non-local deformation gradient tensor used in the mathematical formulation of the constitutive correspondence model \cite{TUPEK201482}.

One aspect of the present formulation, which is also shared by some other numerical approaches to failure mechanics \cite{BEHZADINASAB2021100045}, is that failure events can only occur within the discretized equations of motion, and not in the continuous equations.
For the present methodology, the reason is that the interpolated velocity field used within the IB framework to determine the motion of the material points is a continuous function of space.
Consequently, if we consider the limit in which two material points approach each other, the material velocities of those points will also converge.
In contrast, the spatially discrete equations can achieve discontinuous structural dynamics because the discrete material points have a non-zero lattice spacing, so that nearby material points can experience different velocities, independent of whether there are discontinuities in the Eulerian velocity field.
Although not considered here, modifications to the continuum formulation could allow for failure events for both the continuous and discrete equations of motion.
For instance, using different delta functions for different material points cloud allow nearby points to experience discontinuous dynamics. 
Despite these inconsistencies between the continuous and discrete equations of motion, we remark that an important finding of the current work is that the developed numerical methodology yields convergent and consistent failure predictions for a nontrivial range of grid spacings.


%

Our current IPD formulation is limited to a material that has uniformly distributed volumes along the structural body, which limits the fidelity of the methodology for complex structural geometries. 
This can be seen in the Cook's membrane benchmark Sec.~\ref{s:Benchmark_Cooks} with the stair-step geometry. 
To simulate the deformations of real hyperelastic materials under fluid traction, it is necessary to modify the volume terms in the discrete IPD formulation Eqs.~\eqref{discrete_nonlocal_deformation_gradient}--\eqref{discrete_pd_force}. 
Such modifications will ultimately enable the IPD method to simulate realistic material behaviors in more complex FSI problems. 
We also only consider isotropic material models. 
An important extension of this work will be to consider fiber-reinforced material models like those that have been developed to describe biomaterials.

%
%
%
%

\section*{Acknowledgement}
We gratefully acknowledge research support through NIH Award HL 117063 and NSF Awards OAC 1450327, OAC 1652541, OAC 1931516, and DMS 1929298. 
A.P.S.B acknowledges support through NSF Award OAC 1931368. 
We thank Jae H. Lee and Simone Rossi for their constructive feedback in improving the manuscript. 
We also thank Pablo Seleson for discussions on peridynamics that have helped to improve this manuscript.
Numerical simulations were performed using facilities provided by the University of North Carolina at Chapel Hill through the Research Computing division of UNC Information Technology Services.

\bibliography{references}

\end{document}